\documentclass{amsart}

\def\Z{\mathbb{Z}}

\newcommand{\stacks}{\Delta}
\newcommand{\toposX}{\mathfrak X}
\newcommand{\etale}{\'{e}tale\,}
\newcommand{\Cech}{\v{C}ech\,}
\newcommand{\bigdot}{\bullet}
\newcommand{\coproduct}{\coprod}
\DeclareMathOperator{\Sub}{\mathcal P}
\DeclareMathOperator{\coeq}{\mathfrak CE}
\DeclareMathOperator{\group}{\mathfrak GPD}
\DeclareMathOperator{\colim}{colim}
\DeclareMathOperator{\Top}{Top}
\DeclareMathOperator{\PTop}{Fun^{\text{pres}}}
\DeclareMathOperator{\holim}{holim}
\DeclareMathOperator{\cosk}{cosk}

\DeclareMathOperator{\bd}{\partial}

\DeclareMathOperator{\calU}{\mathcal{U}}
\DeclareMathOperator{\calW}{\mathcal{W}}
\DeclareMathOperator{\calE}{\mathcal{E}}
\DeclareMathOperator{\calB}{\mathcal{B}}
\DeclareMathOperator{\calK}{\mathcal{K}}
\DeclareMathOperator{\Sp}{Sp}
\DeclareMathOperator{\Simp}{\bold{CX}}
\DeclareMathOperator{\calF}{\mathcal{F}}
\DeclareMathOperator{\calG}{\mathcal{G}}
\DeclareMathOperator{\Hom}{Hom} 
\DeclareMathOperator{\HH}{H} \DeclareMathOperator{\BU}{BU}
\DeclareMathOperator{\id}{id} \DeclareMathOperator{\Fun}{Fun}
\DeclareMathOperator{\calC}{\mathcal{C}}
\DeclareMathOperator{\calI}{\mathcal{I}}
\DeclareMathOperator{\calJ}{\mathcal{J}}
\DeclareMathOperator{\SSet}{\mathcal{S}}
\DeclareMathOperator{\calS}{\mathfrak{S}}
\DeclareMathOperator{\calX}{\mathcal{X}}
\DeclareMathOperator{\red}{\text{hyp}}
\DeclareMathOperator{\calY}{\mathcal{Y}}
\DeclareMathOperator{\op}{op}
\DeclareMathOperator{\calD}{\mathcal{D}}
\DeclareMathOperator{\Ind}{Ind} \DeclareMathOperator{\Acc}{Acc}
\DeclareMathOperator{\Pre}{\PTop}
\DeclareMathOperator{\calP}{\mathcal{P}} \topmargin=0in
\newtheorem{theorem}{Theorem}[subsection]
\newtheorem{lemma}[theorem]{Lemma}
\newtheorem{proposition}[theorem]{Proposition}
\newtheorem{corollary}[theorem]{Corollary}
\newtheorem{fact}[theorem]{Fact}

\theoremstyle{definition}
\newtheorem{definition}[theorem]{Definition}
\newtheorem{example}[theorem]{Example}
\newtheorem{counterexample}[theorem]{Counterexample}

\newtheorem{remark}[theorem]{Remark}

\begin{document}

\title{On $\infty$-Topoi}
\author{Jacob Lurie}

\maketitle

Let $X$ be a topological space and $G$ an abelian group. There are
many different definitions for the cohomology group $\HH^n(X,G)$;
we will single out three of them for discussion here. First of
all, one has the singular cohomology $\HH^n_{\text{sing}}(X,G)$,
which is defined as the cohomology of a complex of $G$-valued
singular cochains. Alternatively, one may regard $\HH^n( \bigdot,
G)$ as a representable functor on the homotopy category of
topological spaces, and thereby define $\HH^n_{\text{rep}}(X,G)$
to be the set of homotopy classes of maps from $X$ into an
Eilenberg-MacLane space $K(G,n)$. A third possibility is to use
the sheaf cohomology $\HH^n_{\text{sheaf}}(X, \underline{G})$ of
$X$ with coefficients in the constant sheaf $\underline{G}$ on
$X$.

If $X$ is a sufficiently nice space (for example, a CW complex),
then all three of these definitions agree. In general, however,
all three give different answers. The singular cohomology of $X$
is constructed using continuous maps from simplices $\Delta^k$
into $X$. If there are not many maps {\em into} $X$ (for example
if every path in $X$ is constant), then we cannot expect
$\HH^n_{\text{sing}}(X,G)$ to tell us very much about $X$.
Similarly, the cohomology group $\HH^n_{\text{rep}}(X,G)$ is
defined using maps from $X$ into a simplicial complex, which
(ultimately) relies on the existence of continuous real-valued
functions on $X$. If $X$ does not admit many real-valued
functions, we should not expect $\HH^n_{\text{rep}}(X,G)$ to be a
useful invariant. However, the sheaf cohomology of $X$ seems to be
a good invariant for arbitrary spaces: it has excellent formal
properties in general and sometimes yields good results in
situations where the other approaches do not apply (such as the
{\etale} topology of algebraic varieties).

We shall take the point of view that the sheaf cohomology of a
space $X$ gives the right answer in all cases. We should then ask
for conditions under which the other definitions of cohomology
give the same answer. We should expect this to be true for
singular cohomology when there are many continuous functions {\em
into $X$}, and for Eilenberg-MacLane cohomology when there are
many continuous functions {\em out of} $X$. It seems that the
latter class of spaces is much larger than the former: it
includes, for example, all paracompact spaces, and consequently
for paracompact spaces one can show that the sheaf cohomology
$\HH^n_{\text{sheaf}}(X,G)$ coincides with the Eilenberg-MacLane
cohomology $\HH^n_{\text{rep}}(X,G)$. One of the main results of
this paper is a generalization of the preceding statement to
non-abelian cohomology, and to the case where the coefficient
system $G$ is not necessarily constant.

Classically, the non-abelian cohomology $\HH^1(X,G)$ of $X$ with
coefficients in a possibly non-abelian group $G$ is understood as
classifying $G$-torsors on $X$. When $X$ is paracompact, such
torsors are again classified by homotopy classes of maps from $X$
into an Eilenberg-MacLane space $K(G,1)$. Note that the group $G$
and the space $K(G,1)$ are essentially the same piece of data: $G$
determines $K(G,1)$ up to homotopy equivalence, and conversely $G$
may be recovered as the fundamental group of $K(G,1)$. To make
this canonical, we should say that specifying $G$ is equivalent to
specifying the space $K(G,1)$ {\em together with a base point};
the space $K(G,1)$ alone only determines $G$ up to inner
automorphisms. However, inner automorphisms of $G$ induce the
identity on $\HH^1(X,G)$, so that $\HH^1(X,G)$ is really a functor
which depends only on $K(G,1)$. This suggests the proper
coefficients for non-abelian cohomology are not groups, but
``homotopy types'' (which we regard as purely combinatorial
entities, represented perhaps by simplicial sets). We may define
the non-abelian cohomology $\HH_{\text{rep}}(X,K)$ of $X$ with
coefficients in any simplicial complex $K$ to be the collection of
homotopy classes of maps from $X$ into $K$. This leads to a good
notion whenever $X$ is paracompact. Moreover, we gain a great deal
by allowing the case where $K$ is not an Eilenberg-MacLane space.
For example, if $K = \BU \times \Z$ is the classifying space for
complex K-theory and $X$ is a compact Hausdorff space, then
$\HH_{\text{rep}}(X,K)$ is the usual complex K-theory of $X$,
defined as the Grothendieck group of the monoid of isomorphism
classes of complex vector bundles on $X$.

When $X$ is not paracompact, we are forced to seek a better way of
defining $\HH(X,K)$. Given the apparent power and flexibility of
sheaf-theoretic methods, it is natural to look for some
generalization of sheaf cohomology, using as coefficients
``sheaves of homotopy types on $X$.'' In other words, we want a
theory of $\infty$-stacks (in groupoids) on $X$, which we will
henceforth refer to simply as {\it stacks}. One approach to this
theory is provided by the Joyal-Jardine homotopy theory of
simplicial presheaves on $X$. According to this approach, if $K$
is a simplicial set, then the cohomology of $X$ with coefficients
in $K$ should be defined as $\HH_{JJ}(X,K) = \pi_0( \calF(X))$,
where $\calF$ is a fibrant replacement for the constant simplicial
presheaf with value $K$ on $X$. When $K$ is an Eilenberg-MacLane
space $K(G,n)$, then this agrees with the sheaf-cohomology group
(or set) $\HH^n_{\text{sheaf}}(X,G)$. It follows that if $X$ is
paracompact, then $\HH_{JJ}(X,K) = \HH_{\text{rep}}(X,K)$ whenever
$K$ is an Eilenberg-MacLane space.

However, it turns out that $\HH_{JJ}(X,K) \neq
\HH_{\text{rep}}(X,K)$ in general, even when $X$ is paracompact.
In fact, one can give an example of a compact Hausdorff space for
which $\HH_{JJ}(X, BU \times \Z)$ is not equal to the complex
$K$-theory of $X$. We shall proceed on the assumption that $K(X)$
is the ``correct'' answer in this case, and give an alternative to
the Joyal-Jardine theory which computes this answer. Our
alternative is distinguished from the Joyal-Jardine theory by the
fact that we require our stacks to satisfy a descent condition
only for coverings, rather than for arbitrary hypercoverings.
Aside from this point we proceed in the same way, setting
$\HH(X,K)  = \pi_0( \calF'(X))$, where $\calF'$ is the stack which
is obtained by forcing the ``constant prestack with value $K$'' to
satisfy this weaker form of descent. In general, $\calF'$ will not
satisfy descent for hypercoverings, and consequently it will not
be equivalent to the simplicial presheaf $\calF$ used in the
definition of $\HH_{JJ}$.

The resulting theory has the following properties:

\begin{itemize}
\item If $X$ is paracompact, $\HH(X,K)$ is the set of homotopy
classes from $X$ into $K$.

\item If $X$ is a paracompact space of finite covering dimension,
then our theory of stacks is equivalent to the Joyal-Jardine
theory. (This is also true for certain inductive limits of finite
dimensional spaces, and in particular for CW complexes.)

\item The cohomologies $\HH_{JJ}(X,K)$ and $\HH(X,K)$ always agree
when $K$ is ``truncated'', for example when $K$ is an
Eilenberg-MacLane space. In particular, $\HH(X,K(G,n))$ is equal
to the usual sheaf cohomology $\HH^n_{\text{sheaf}}(X,G)$.
\end{itemize}

In addition, our theory of $\infty$-stacks enjoys good formal
properties which are not always shared by the Joyal-Jardine
theory; we shall summarize the situation in \S \ref{versus}.
However, the good properties of our theory do not come without
their price. It turns out that the essential difference between
stacks (which are required to satisfy descent only for ordinary
coverings) and hyperstacks (which are required to satisfy descent
for arbitrary hypercoverings) is that the former can fail to
satisfy the Whitehead theorem: one can have, for example, a
pointed stack $(E,\eta)$ for which $\pi_i(E,\eta)$ is a trivial
sheaf for all $i \geq 0$, such that $E$ is not ``contractible''
(for the definition of these homotopy sheaves, see \S
\ref{homotopysheaves}).

In order to make a thorough comparison of our theory of stacks on
$X$ and the Joyal-Jardine theory of hyperstacks on $X$, it seems
desirable to fit both of them into some larger context. The proper
framework is provided by the notion of an $\infty$-topos, which is
intended to be an $\infty$-category that ``looks like'' the
$\infty$-category of $\infty$-stacks on a topological space, just
as an ordinary topos is supposed to be a category that ``looks
like'' the category of sheaves on a topological space. For any
topological space $X$ (or, more generally, any topos), the
$\infty$-stacks on $X$ comprise an $\infty$-topos, as do the
$\infty$-hyperstacks on $X$. However, it is the former
$\infty$-topos which enjoys the more universal position among
$\infty$-topoi related to $X$ (see Proposition
\ref{universality}).

Let us now outline the contents of this paper. We will begin in \S
\ref{categories} with an informal review of the theory of
$\infty$-categories. There are many approaches to the foundation
of this subject, each having its own particular merits and
demerits. Rather than single out one of those foundations here, we
shall attempt to explain the ideas involved and how to work with
them. The hope is that this will render this paper readable to a
wider audience, while experts will be able to fill in the details
missing from our exposition in whatever framework they happen to
prefer.

Section \ref{toposes} is devoted to the notion of an
$\infty$-topos. We will begin with an intrinsic characterization
(analogous to Giraud's axioms which characterize ordinary topoi:
see \cite{SGA}), and then argue that any $\infty$-category
satisfying our axioms actually arises as an $\infty$-category of
``stacks'' on {\em something}. We will then show that any
$\infty$-topos determines an ordinary topos in a natural way, and
vice versa. We will also explain the relationship between stacks
and hyperstacks.

Section \ref{paracompactness} relates the $\infty$-topos of stacks
on a paracompact space with some notions from classical homotopy
theory. In particular, we prove that the ``non-abelian
cohomology'' determined by the $\infty$-topos associated to a
paracompact space agrees, in the case of constant coefficients,
with the functor $H_{\text{rep}}(X,K)$ described above.

Section \ref{dimension} is devoted to the dimension theory of
$\infty$-topoi. In particular, we shall define the {\it homotopy
dimension} of an $\infty$-topos, which simultaneously generalizes
the covering dimension of paracompact topological spaces and the
Krull dimension of Noetherian topological spaces. We also prove a
generalization of Grothendieck's vanishing theorem for cohomology
on a Noetherian topological space of finite Krull dimension.

This paper is intended to be the first in a series devoted to the
notion of an $\infty$-topos. In the future we plan to discuss
proper maps of $\infty$-topoi, the proper base change theorem, and
a notion of ``elementary $\infty$-topos'' together with its
relationship to mathematical logic.

I would like to thank Max Lieblich, whose careful reading of this
document and many suggestions have improved the exposition
considerably. I would also like to thank Bertrand To\"{e}n, with
whom I shared several stimulating conversations about the subject
matter of this paper.

\section{$\infty$-Categories}\label{categories}

\subsection{General Remarks}
Throughout this paper, we will need to use the language of
$\infty$-categories. Since there are several approaches to this
subject in the literature, we devote this first section to
explaining our own point of view. We begin with some general
remarks directed at non-experts; seasoned homotopy-theorists and
category-theorists may want to skip ahead to the next section.

An ordinary category consists of a collection of objects, together
with sets of morphisms between those objects. The reader may also
be familiar with the notion of a $2$-category, in which there are
not only morphisms but also {\em morphisms between the morphisms},
which are called $2$-morphisms. The vision of higher categories is
that one should be able to discuss $n$-categories for any $n$, in
which one has not only objects, morphisms, and $2$-morphisms, but
$k$-morphisms for all $k \leq n$. Finally, in some sort of limit
one should obtain a theory of $\infty$-categories where one has
morphisms of all orders.

There are many approaches to realizing this vision. We might begin
by defining a $2$-category to be a ``category enriched over
categories.'' In other words, one considers a collection of
objects together with a {\em category} of morphisms $\Hom(A,B)$
for any two objects $A$ and $B$, and composition {\em functors}
$c_{ABC}: \Hom(A,B) \times \Hom(B,C) \rightarrow \Hom(A,C)$ (to
simplify the discussion, we shall ignore identity morphisms for a
moment). These functors are required to satisfy an associative
law, which asserts that for any quadruple $(A,B,C,D)$ of objects,
$c_{ACD} \circ (c_{ABC} \times 1) = c_{ABD} \circ (1 \times
c_{BCD})$ as functors $$\Hom(A,B) \times \Hom(B,C) \times
\Hom(C,D) \rightarrow \Hom(A,D)$$ This leads to the definition of
a {\it strict $2$-category}.

At this point, we should object that the definition of a strict
$2$-category violates one of the basic philosophical principles of
category theory: one should never demand that two functors be
``equal.'' Instead one should postulate the existence of a natural
transformation between two functors. This means that the
associative law should not take the form of an equation, but of
additional structure: a natural isomorphism $\gamma_{ABCD}:
c_{ACD} \circ (c_{ABC} \times 1) \simeq c_{ABD} \circ (1 \times
c_{BCD})$. We should also demand the natural transformations
$\gamma_{ABCD}$ be functorial in the quadruple $(A,B,C,D)$, and
that they satisfy certain higher associativity conditions. After
formulating the appropriate conditions, we arrive at the
definition of a {\it weak $2$-category}, or simply a {\it
$2$-category}.

Let us contrast the notions of ``strict $2$-category'' and ``weak
$2$-category.'' The former is easier to define, since we do not
have to worry about the higher associativity conditions satisfied
by the transformations $\gamma_{ABCD}$. On the other hand, the
latter notion seems more natural if we take the philosophy of
category theory seriously. In this case, we happen to be lucky:
the notions of ``strict $2$-category'' and ``weak $2$-category''
turn out to be equivalent. More precisely, any weak $2$-category
can be replaced by an ``equivalent'' strict $2$-category. The
question of which definition to adopt is therefore an issue of
aesthetics.

Let us plunge onward to $3$-categories. Following the above
program, we may define a {\it strict $3$-category} to consist of a
collection of objects together with strict $2$-categories
$\Hom(A,B)$ for any pair of objects $A$ and $B$, together with a
strictly associative composition law. Alternatively, we could seek
a definition of {\it weak $3$-category} by allowing $\Hom(A,B)$ to
be only a weak $2$-category, requiring associativity only up to
natural $2$- isomorphisms, which satisfy higher associativity laws
up to natural $3$-isomorphisms, which in turn satisfy still higher
associativity laws of their own. Unfortunately, it turns out that
these notions are {\em not} equivalent.

Both of these approaches have serious drawbacks. The obvious
problem with weak $3$-categories is that an explicit definition is
extremely complicated (see \cite{tricat}, where a definition is
worked out along these lines), to the point where it is
essentially unusable. On the other hand, strict $3$-categories
have the problem of not being the correct notion: most of the weak
$3$-categories which occur in nature (such as the fundamental
$3$-groupoid of a topological space) are not equivalent to strict
$3$-categories. The situation only gets worse (from either point
of view) as we pass to $4$-categories and beyond.

Fortunately, it turns out that major simplifications can be
introduced if we are willing to restrict our attention to
$\infty$-categories in which most of the higher morphisms are
invertible. Let us henceforth use the term $(\infty,n)$-category
to refer to $\infty$-categories in which all $k$-morphisms are
invertible for $k > n$. It turns out that $(\infty,0)$-categories
(that is, $\infty$-categories in which {\em all} morphisms are
invertible) have been extensively studied from another point of
view: they are the same thing as ``spaces'' in the sense of
homotopy theory, and there are many equivalent ways of describing
them (for example, using CW complexes or simplicial sets).

We can now proceed to define an $(\infty,1)$-category: it is
something which has a collection of objects, and between any two
objects $A$ and $B$ there is an $(\infty,0)$-category, or
``space,'' called $\Hom(A,B)$. We require that these morphism
spaces are equipped with an associative composition law. As
before, we are faced with two choices as to how to make this
precise: do we require associativity on the nose, or only ``up to
homotopy'' in some sense? Fortunately, it turns out not to matter:
as was the case with $2$-categories, any $(\infty,1)$-category
with a coherently associative multiplication can be replaced by an
equivalent $(\infty,1)$-category with a strictly associative
multiplication.

In this paper, we will deal almost exclusively with
$(\infty,1)$-categories (with a few exceptional appearances of
$(\infty,2)$-categories, for which we will not require any general
theory). {\em Unless otherwise specified, all $\infty$-categories
will be assumed to be $(\infty,1)$-categories.}

There are a number of models for these $\infty$-categories:
categories enriched over simplicial sets, categories enriched over
topological spaces, Segal categories (\cite{toen}), simplicial
sets satisfying a weak Kan condition (called quasi-categories by
Joyal, see \cite{quasicat} and \cite{weakKan}). All of these
approaches give essentially the same notion of $\infty$-category,
so it does not matter which approach we choose.

\subsection{$\infty$-Categories}

The first ingredient needed for the theory of $\infty$-categories
is the theory of $\infty$-groupoids, or ``spaces'': namely,
homotopy theory. Our point of view on this subject is that exists
some abstract and purely combinatorial notion of a ``homotopy
type'' which is $\infty$-categorical in nature, and therefore
slippery and difficult to define in terms of sets. However, there
are various ways to ``model'' homotopy types by classical
mathematical objects such as CW complexes or simplicial sets. Let
us therefore make the definition that a {\it space} is a CW
complex. We shall agree that when we discuss a space $X$, we shall
refer only to its homotopy-theoretic properties and never to
specific features of some particular representative of $X$.

We can now say what an $\infty$-category is:

\begin{definition}\label{ic}
An {\it $\infty$-category} $\calC$ consists of a collection of
objects, together with a space $\Hom_{\calC}(X,Y)$ for any pair of
objects $X,Y \in \calC$. These $\Hom$-spaces must be equipped with
an associative composition law
$$\Hom_{\calC}(X_0, X_1) \times \Hom_{\calC}(X_1, X_2) \times
\ldots \Hom_{\calC}(X_{n-1},X_n) \rightarrow
\Hom_{\calC}(X_0,X_n)$$ (defined for all $n \geq 0$).
\end{definition}

\begin{remark}
It is customary to use the compactly generated topology on
$$\Hom_{\calC}(X_0,X_1) \times \ldots \times \Hom_{\calC}(X_{n-1},X_n),$$
rather than the product topology, so that the product remains a CW
complex. This facilitates comparisons with more combinatorial
versions of the theory, but the issue is not very important.
\end{remark}

\begin{remark}
Of the numerous possible definitions of $\infty$-categories,
Definition \ref{ic} is the easiest to state and to understand.
However, it turns out to be one of the hardest to work with
explicitly, since many basic $\infty$-categorical constructions
are difficult to carry out at the level of topological categories.
We will not dwell on these technical points, which are usually
more easily addressed using the more sophisticated approaches.
\end{remark}

Let us now see how to work with the $\infty$-categories introduced
by Definition \ref{ic}. Note first that any $\infty$-category
$\calC$ determines an ordinary category $h\calC$ having the same
objects, but with $\Hom_{h\calC} (X,Y) = \pi_0 \Hom_{\calC}(X,Y)$.
The category $h\calC$ is called the {\it homotopy category} of
$\calC$ (or sometimes the {\it derived category} of $\calC$). To
some extent, working in the $\infty$-category $\calC$ is like
working in its homotopy category $h\calC$: up to equivalence,
$\calC$ and $h\calC$ have the same objects and morphisms. The
difference between $h\calC$ and $\calC$ is that in $\calC$, one
must not ask about whether or not morphisms are {\em equal};
instead one should ask whether or not one can find a path from one
to the other. One consequence of this difference is that the
notion of a commutative diagram in $h\calC$, which corresponds to
a {\it homotopy commutative} diagram in $\calC$, is quite
unnatural and usually needs to be replaced by the more refined
notion of a {\it homotopy coherent} diagram in $\calC$.

To understand the problem, let us suppose that $F: \calI
\rightarrow h\SSet$ is a functor from an ordinary category $\calI$
into the homotopy category of spaces $\SSet$. In other words, $F$
assigns to each object $x \in \calI$ a space $Fx$, and to each
morphism $\phi: x \rightarrow y$ in $\calI$ a continuous map of
spaces $F\phi: Fx \rightarrow Fy$ (well-defined up to homotopy),
such that $F(\phi \circ \psi)$ is homotopic to $F\phi \circ F\psi$
for any pair of composable morphisms $\phi, \psi$ in $\calI$. In
this situation, it may or may not be possible to {\em lift} $F$ to
an actual functor $\widetilde{F}$ from $\calI$ to the ordinary
category of topological spaces, such that $\widetilde{F}$ induces
a functor $\calI \rightarrow h \SSet$ which is equivalent to $F$.
In general there are obstructions to both the existence and the
uniqueness of the lifting $\widetilde{F}$, even up to homotopy. To
see this, we note that $\widetilde{F}$ determines extra data on
$F$: for every composable pair of morphisms $\phi$ and $\psi$,
$\widetilde{F}(\phi \circ \psi) = \widetilde{F} \phi \circ
\widetilde{F} \psi$, which means that $\widetilde{F}$ gives a {\em
specified} homotopy $h_{\phi, \psi}$ between $F(\phi \circ \psi)$
and $F\phi \circ F \psi$. We should imagine that the functor $F$
to the homotopy category $h \SSet$ is a first approximation to
$\widetilde{F}$; we obtain a second approximation when we take
into account the homotopies $h_{\phi, \psi}$. These homotopies are
not arbitrary: the associativity of composition gives a
relationship between $h_{\phi, \psi}, h_{\psi, \theta}, h_{\phi,
\psi \circ \theta}$ and $h_{\phi \circ \psi, \theta}$, for a
composable triple of morphisms $(\phi, \psi, \theta)$ in $\calI$.
This relationship may be formulated in terms of the existence of a
certain higher homotopy, which is once again canonically
determined by $\widetilde{F}$. To obtain the next approximation to
$\widetilde{F}$, we should take these higher homotopies into
account, and formulate the associativity properties that {\em
they} enjoy, and so forth. A {\it homotopy coherent} diagram in
$\calC$ is, roughly speaking, a functor $F: \calI \rightarrow
h\calC$, together with all of the extra data that would be
available if there we could lift $F$ to get some $\widetilde{F}$
which was functorial ``on the nose.''

An important consequence of the distinction between homotopy
commutativity and homotopy coherence is that the appropriate
notions of {\it limit} and {\it colimit} in $\calC$ do not
coincide with the notion of {\it limit} and {\it colimit} in
$h\calC$ (in which limits and colimits typically do not exist).
The appropriately defined limits and colimits in $\calC$ are
typically referred to as {\it homotopy limits} and {\it homotopy
colimits}, to avoid confusing them ordinary limits and colimits
inside of some category of models for $\calC$. We will try to
avoid using categories of models at all and work in an
``invariant'' fashion with an $\infty$-category $\calC$. In
particular, the terms {\it limit} and {\it colimit} in an
$\infty$-category $\calC$ will {\em always} mean homotopy limit
and homotopy colimit. We shall never speak of any other kind of
limit or colimit.

\begin{example}
Suppose given a collection of objects $\{X_{\alpha} \}$ in an
$\infty$-category $\calC$. The {\it product} $X = \prod_{\alpha}
X_{\alpha}$ (if it exists) is characterized by the following
universal property: $\Hom_{\calC}(Y, X) \simeq \prod_{\alpha}
\Hom_{\calC}(Y, X_{\alpha})$. Passing to connected components, we
see that we also have $\Hom_{h \calC}(Y,X) \simeq \prod_{\alpha}
\Hom_{h \calC}(Y, X_{\alpha})$. Consequently, any product in
$\calC$ is also a product in $h \calC$.
\end{example}

\begin{example}
Given two maps $\pi: X \rightarrow Z$ and $\psi: Y \rightarrow Z$
in an $\infty$-category $\calC$, the {\it fiber product} $X
\times_Z Y$ (if it exists) is characterized by the following
universal property: to specify a map $W \rightarrow X \times_Z Y$
is to specify maps $W \rightarrow X$, $W \rightarrow Y$, {\em
together with a homotopy between the induced composite maps $W
\rightarrow Z$}. In particular, there is a map from $X \times_Z Y$
to the analogous fiber product in $h \calC$ (if this fiber product
exists), which need not be an isomorphism.
\end{example}

We want to emphasize the point of view that $\infty$-categories
(with $n$-morphisms assumed invertible for $n
> 1$) are in many ways a less drastic generalization of categories
than ordinary $2$-categories (in which there can exist
non-invertible $2$-morphisms). Virtually all definitions and
constructions which make sense for ordinary categories admit
straightforward generalizations to $\infty$-categories, while for
$2$-categories even the most basic notions need to be rethought
(for example, one has a distinction between strict and lax
pullbacks).

We now summarize some of the basic points to keep in mind when
dealing with $\infty$-categories.

\begin{itemize}

\item Any ordinary category may be considered as an
$\infty$-category: one takes each of the morphism spaces
$\Hom_{\calC}(X,Y)$ to be discrete.

\item By a {\it morphism} from $X$ to $Y$ in $\calC$, we mean a point
of the space $\Hom_{\calC}(X,Y)$. Two morphisms are {\it
equivalent}, or {\it homotopic}, if they lie in the same path
component of $\Hom_{\calC}(X,Y)$.

\item If $f: C \rightarrow C'$ is a morphism in an
$\infty$-category $\calC'$, then we say that $f$ is an {\it
equivalence} if it becomes an isomorphism in the homotopy category
$h\calC$. In other words, $f$ is an equivalence if it has a
homotopy inverse.

\item Any $\infty$-category $\calC$ has an {\it opposite
$\infty$-category} $\calC^{op}$, having the same objects but with
$\Hom_{\calC^{op}}(X,Y) = \Hom_{\calC}(Y,X)$.

\item Given an $\infty$-category $\calC$ and an object $X \in
\calC$, we can define a {\it slice $\infty$-category}
$\calC_{/X}$. The objects of $\calC_{/X}$ are pairs $(Y,f)$ with
$Y \in \calC$ and $f \in \Hom_{\calC}(Y,X)$. The space of
morphisms from $\Hom_{\calC_{/X}}((Y,f), (Y',f'))$ is given by the
homotopy fiber of the map $\Hom_{\calC}(Y,Y') \stackrel{f' \circ
}{\rightarrow} \Hom_{\calC}(Y,X)$ over the point $f$.

If the category $\calC$ and the maps $f: Y \rightarrow X$, $f': Y'
\rightarrow X$ are clear from context, we will often write
$\Hom_{X}(Y,Y')$ for $\Hom_{\calC_{/X}}((Y,f),(Y',f'))$.

\begin{remark}
Let $\calC_0$ be a category which is enriched over topological
spaces (in other words, the morphism sets in $\calC$ are equipped
with topologies such that the composition laws are all
continuous), and let $\calC$ denote the $\infty$-category which it
represents. Given an object $X \in \calC_0$, one may form the
classical slice-category construction to obtain a category
${\calC_0}_{/ X}$, which is again enriched over topological
spaces. The associated $\infty$-category is not necessarily
equivalent to the slice $\infty$-category $\calC_{/X}$. The
difference is that the morphism spaces in $\calC_{/X}$ must be
formed using homotopy fibers of maps $\phi: \Hom_{\calC}(Y,Y')
\rightarrow \Hom_{\calC}(Y,X)$, while the morphism spaces in
${\calC_0}_{/X}$ are given by taking the ordinary fibers of the
maps $\phi$. If we want to construct an explicit model for
$\calC_{/X}$ as a topological category, then we need to choose an
explicit construction for the homotopy fiber of a continuous map
$\phi: E \rightarrow E'$. The standard construction for the
homotopy fiber of $\phi$ over a point $e' \in E'$ is the space $\{
(e, p): e \in E, p: [0,1] \rightarrow E', p(0)=\phi(e),
p(1)=e'\}$. However, if one uses this construction to define the
morphism spaces in $\calC_{/X}$, then it is not possible to define
a {\em strictly} associative composition law on morphisms.
However, there exists a composition law which is associative up to
{\em coherent} homotopy. This may then be replaced by a strictly
associative composition law after slightly altering the morphism
spaces.

The same difficulties arise with virtually every construction we
shall encounter. In order to perform some categorical construction
on an $\infty$-category, it is not enough to choose a topological
(or simplicial) category which models it and employ the same
construction in the naive sense. However, this is merely a
technical annoyance and not a serious problem: there are many
strategies for dealing with this, one (at least) for each of the
main definitions of $\infty$-categories. We will henceforth ignore
these issues; the reader may refer to the literature for more
detailed discussions and constructions.
\end{remark}

\item We call an $\infty$-category {\it small} if it has a set of
objects, just as with ordinary categories. We call an
$\infty$-category {\it essentially small} if the collection of
equivalence classes of objects forms a set.

\item Let $\calC$ and $\calC'$ be $\infty$-categories. One can
define a notion of functor $F: \calC \rightarrow \calC'$. A
functor $F$ carries objects of $C \in \calC$ to objects $FC \in
\calC'$, and induces maps between the morphism spaces
$\Hom_{\calC}(C,D) \rightarrow \Hom_{\calC'}(FC,FD)$. These maps
should be compatible with the composition operations, up to
coherent homotopy. See \cite{cordier} for an appropriate
definition, in the context of simplicial categories. We remark
that if we are given explicit models for $\calC$ and $\calC'$ as
simplicial or topological categories, then we do not expect every
functor $F$ to arise from an ordinary functor (in other words, we
do not expect to be able to arrange the situation so that $F$ is
compatible with composition ``on the nose''), although any
ordinary functor which is compatible with a simplicial or
topological enrichment does give rise to a functor between the
associated $\infty$-categories.

\item If $\calC$ and $\calC'$ are two small $\infty$-categories, then
there is an $\infty$-category $\Fun(\calC, \calC')$ of functors
from $\calC$ to $\calC'$. We will sometimes write
${\calC'}^{\calC}$ for $\Fun(\calC,\calC')$. We may summarize the
situation by saying that there is an $(\infty,2)$-category of
small $\infty$-categories.

\item We will also need to discuss functor categories
$\Fun(\calC,\calC')$ in cases when $\calC$ and $\calC'$ are not
small. So long as $\calC$ is small, this poses little difficulty:
we must simply bear in mind that if $\calC'$ is large, the functor
category is also large. Care must be taken if $\calC$ is also
large: in this case $\Fun(\calC, \calC')$ is ill-defined because
its collection of objects could be {\em very large} (possibly not
even a proper class), and the space of natural transformations
between two functors might also be large. In practice, we will
avoid this difficulty by only considering functors satisfying
certain continuity conditions. After imposing these restrictions,
$\Fun(\calC,\calC')$ will be an honest $\infty$-category (though
generally still a large $\infty$-category).

\item Given an $\infty$-category $\calC$ and some collection of
objects $S$ of $\calC$, we can form an $\infty$-category $\calC_0$
having the elements of $S$ as objects, and the same morphism
spaces as $\calC$. We shall refer to those $\infty$-categories
$\calC_0$ which arise in this way {\it full subcategories} of
$\calC$. (Our use of ``full subcategory" as opposed to ``full
sub-$\infty$-category" is to avoid awkward language and is not
intended to suggest that $\calC_0$ is an ordinary category).

\item A functor $F: \calC \rightarrow \calC'$ is said to be fully
faithful if for any objects $C,C' \in \calC$, the induced map
$\Hom_{\calC}(C,C') \rightarrow \Hom_{\calC'}(FC,FC')$ is a
homotopy equivalence. We say that $F$ is {\it essentially
surjective} if every object of $\calC'$ is equivalent to an object
of the form $FC$, $C \in \calC$. We say that $F$ is an {\it
equivalence} if it is fully faithful and essentially surjective.
In this case, one can construct a homotopy inverse for $F$ if one
assumes a sufficiently strong version of the axiom of choice.

We shall understand that all relevant properties of
$\infty$-categories are invariant under equivalence. For example,
the property of being small is not invariant under equivalence,
and is therefore not as natural as the property of being
essentially small. We note that an $\infty$-category is
essentially small if and only if it is equivalent to a small
$\infty$-category.

\item If a functor $F: \calC \rightarrow \calC'$ is fully
faithful, then it establishes an equivalence of $\calC$ with a
full subcategory of $\calC'$. In this situation, we may use $F$ to
identify $\calC$ with its image in $\calC'$.

\item If $X$ is an object in an $\infty$-category $\calC$, then we
sa that $X$ is {\it initial} if $\Hom_{\calC}(X, Y)$ is
contractible for all $Y \in \calC$. The dual notion of a {\it
final} object is defined in the evident way.

\item If $\calC$ is an $\infty$-category, then a {\it diagram} in
$\calC$ is a functor $F: \calI \rightarrow \calC$, where $\calI$
is a small $\infty$-category. Just as with diagrams in ordinary
categories, we may speak of limits and colimits of diagrams, which
may or may not exist. These may be characterized in the following
way (we restrict our attention to colimits; for limits, just
dualize everything): a colimit of $F$ is an object $C \in \calC$
such that there exists a natural transformation $F \rightarrow
F_C$ which induces equivalences $\Hom_{\calC^{\calI}}(F, F_{D})
\simeq \Hom_{\calC}(C,D)$, where $F_C$ and $F_D$ denote the
constant functors $\calI \rightarrow \calC$ having the values $C,D
\in \calC$. As we have explained, these are {\em
homotopy-theoretic} limits and colimits, and do not necessarily
enjoy any universal property in $h\calC$.

\item In the theory of ordinary categories (and in mathematics in
general), the category of sets plays a pivotal role. In the
$\infty$-categorical setting, the analogous role is filled by the
$\infty$-category $\SSet$ of spaces. One can take the objects of
$\SSet$ to be any suitable model for homotopy theory, such as CW
complexes or fibrant simplicial sets: in either case, there are
naturally associated ``spaces'' of morphisms, which may themselves
be interpreted as objects of $\SSet$.

\item An $\infty$-category is an {\it $\infty$-groupoid} if all of
its morphisms are equivalences. A small $\infty$-groupoid is
essentially the same thing as an object of $\SSet$.
\begin{remark}\label{content}
This assertion has some content: for $\infty$-groupoids with a
single object $\ast$, it reduces to Stasheff's theorem that the
coherently associative composition on $X = \Hom(\ast,\ast)$ is
precisely the data required to realize $X$ as a loop space (see
for example \cite{goerss}, where a version of this theorem is
proved by establishing the equivalence of the homotopy theory of
simplicial sets and simplicial groupoids). \end{remark}

\item For any $\infty$-category $\calC$, we define a {\it
prestack} on $\calC$ to be a functor $F: \calC^{\op} \rightarrow
\SSet$. We remark that this is not standard terminology: for us,
prestacks and stacks will always be valued in $\infty$-groupoids,
rather than in ordinary groupoids or in categories.

\item
Any object $C \in \calC$ gives rise to a prestack $\calF_C$ given
by the formula $\calF_C(C') = \Hom_{\calC}(C',C)$. This induces a
functor from $\calC$ to the $\infty$-category $\SSet^{\calC^{op}}$
of prestacks on $\calC$, which generalizes the classical Yoneda
embedding for ordinary categories. As with the classical Yoneda
embedding, this functor is fully faithful; prestacks which lie in
its essential image are called {\it representable}.

\item If $\calC_0$ is a full subcategory of an $\infty$-category
$\calC$, then one can also define a ``restricted Yoneda embedding"
$\calC \rightarrow \SSet^{\calC_0^{op}}$. This functor is fully
faithful whenever $\calC_0$ generates $\calC$, in the sense that
every object of $\calC$ can be obtained as a colimit of objects of
$\calC_0$.

\item The $(\infty,2)$-category of $\infty$-categories has all
$(\infty,2)$-categorical limits and colimits. We do not want to
formulate a precise statement here, but instead refer the reader
to Appendix \ref{appendixdiagram} for some discussion. Let us
remark here that forming {\em limits} is relatively easy: one can
work with them ``componentwise,'' just as with ordinary
categories.

\begin{example}
Given any family $\{ \calC_{\alpha} \}$ of $\infty$-categories,
their product $\calC = \prod_{\alpha} \calC_{\alpha}$ may be
constructed in the following way: the objects of $\calC$ consist
of a choice of one object $C_{\alpha}$ from each $\calC_{\alpha}$,
and the morphism spaces $\Hom_{\calC}( \{ C_{\alpha} \}, \{
C'_{\alpha} \}) = \prod_{\alpha} \Hom_{\calC_{\alpha}}(
C_{\alpha}, C'_{\alpha} )$.
\end{example}

\begin{example}
Suppose that $F: \calC' \rightarrow \calC$ and $G: \calC''
\rightarrow \calC$ are functors between $\infty$-categories. One
may form an $\infty$-category $\calC' \times_{\calC} \calC''$, the
{\it strict fiber product} of $\calC'$ and $\calC''$ over $\calC$.
An object of $\calC' \times_{\calC} \calC''$ consists of a pair $(
C', C'', \eta)$ where $C' \in \calC'$, $C'' \in \calC''$, and
$\eta: FC' \rightarrow FC''$ is an equivalence in $\calC$. If one
drops the requirement that $\eta$ be an equivalence, then one
obtains instead the {\it lax fiber product} of $\calC'$ and
$\calC''$ over $\calC$.
\end{example}

\begin{example}
Let $\calC$ and $\calI$ be $\infty$-categories. Then the functor
$\infty$-category $\calC^{\calI}$ may be regarded as a limit of
copies of the $\infty$-category $\calC$, where the limit is
indexed by the $\infty$-category $\calI$.
\end{example}

\begin{example}\label{slicesarelax}
Let $\calC$ be an $\infty$-category and $E \in \calC$ an object.
Then we may regard $E$ as a determining a functor $\ast
\rightarrow \calC$, where $\ast$ is an $\infty$-category with a
single object having a contractible space of endomorphisms. The
slice category $\calC_{/ E}$ may be regarded as a lax fiber
product of $\calC$ and $\ast$ over $\calC$ (via the functors $E$
and the identity).
\end{example}

\item One may describe $\infty$-categories by ``generators and
relations''. In particular, it makes sense to speak of a {\it
finitely presented} $\infty$-category. Such an $\infty$-category
has finitely many objects and its morphism spaces are determined
by specifying a finite number of generating morphisms, a finite
number of relations among these generating morphisms, a finite
number of relations among the relations, and so forth (a finite
number of relations in all).

\begin{example}\label{infinitemorphisms}
Let $\calC$ be the free $\infty$-category generated by a single
object $X$ and a single morphism $f: X \rightarrow X$. Then
$\calC$ is a finitely presented $\infty$-category with a single
object, and $\Hom_{\calC}(X,X) = \{ 1, f, f^2, \ldots \}$ is
infinite discrete. In particular, we note that the finite
presentation of $\calC$ does not guarantee finiteness properties
of the morphism spaces.
\end{example}

\begin{example}
If we identify $\infty$-groupoids with spaces, then writing down a
presentation for an $\infty$-groupoid corresponds to giving a cell
decomposition of the associated space. We therefore see that the
finitely presented $\infty$-groupoids correspond precisely to the
finite cell complexes.
\end{example}

\begin{example}
Suppose that $\calC$ is an $\infty$-category with only two objects
$X$ and $Y$, and that $X$ and $Y$ have contractible endomorphism
spaces and that $\Hom_{\calC}(X,Y)$ is empty. Then $\calC$ is
completely determined by the morphism space $\Hom_{\calC}(Y,X)$,
which may be arbitrary. The $\infty$-category $\calC$ is finitely
presented if and only if $\Hom_{\calC}(Y,X)$ is a finite cell
complex (up to homotopy equivalence).
\end{example}

\item
Given two $\infty$-categories $\calC$ and $\calC'$ and functors
$F: \calC \rightarrow \calC'$, $G: \calC' \rightarrow \calC$, we
have two functors
$$\Hom_{\calC}(X,GY), \Hom_{\calC'}(FX,Y): \calC^{op} \times \calC'
\rightarrow \SSet$$ An equivalence between these two functors is
called an {\it adjunction} between $F$ and $G$: we say that $F$ is
{\it left adjoint} to $G$ and that $G$ is {\it right adjoint} to
$F$. As with ordinary categories, adjoints are unique up to
canonical equivalence when they exist.

\item As with ordinary categories, many of the $\infty$-categories
which arise in nature are large. Working with these objects
creates foundational technicalities. We will generally ignore
these technicalities, which can be resolved in many different
ways.

\end{itemize}

\subsection{Accessibility}

Most of the categories which arise in nature (such as the category
of groups) have a proper class of objects, even when taken up to
isomorphism. However, there is a sense in which many of these
large categories are determined by a bounded amount of
information. For example, in the category of groups, every object
may be represented as a filtered colimit of finitely presented
groups, and the category of finitely presented groups is
essentially small.

This phenomenon is common to many examples, and may be formalized
in the notion of an {\it accessible category} (see for example
\cite{adamek} or \cite{makkai}). In this section, we shall adapt
the definition of accessible category and the basic facts
concerning them to the $\infty$-categorical setting. We shall not
assume that the reader is familiar with the classical notion of an
admissible category.

Our reasons for studying accessibility are twofold: first, we
would prefer to deal with $\infty$-categories as legitimate
mathematical objects, without appealing to a hierarchy of
universes (in the sense of Grothendieck) or some similar device.
Second, we sometimes wish to make arguments which play off the
difference between a ``large'' category and a ``small'' piece
which determines it, in order to make sense of constructions which
would otherwise be ill-defined. See, for example, the proof of
Theorem \ref{representable}.

Before we begin, we are going to need to introduce some
size-related concepts. Recall that a cardinal $\kappa$ is {\it
regular} if $\kappa$ cannot be represented as a sum of fewer than
$\kappa$ cardinals of size $< \kappa$.

\begin{definition}
Let $\kappa$ be a regular cardinal. If $\kappa > \omega$, then an
$\infty$-category $\calC$ is called {\it $\kappa$-small} it has
fewer than $\kappa$ objects, and for any pair of objects $X,Y$,
the space $\Hom_{\calC}(X,Y)$ may be presented using fewer than
$\kappa$ cells. If $\kappa = \omega$, then an $\infty$-category
$\calC$ is {\it $\kappa$-small} if it is finitely presented.
\end{definition}

\begin{remark}\label{brak}
One may unify the two clauses of the definition as follows: an
$\infty$-category $\calC$ is $\kappa$-small if it may be presented
by fewer than $\kappa$ ``generators and relations''. The reason
that the above definition needs a special case for $\kappa =
\omega$ is that a finitely presented category $\infty$-category
can have morphism spaces which are not finite (see Example
\ref{infinitemorphisms}).
\end{remark}

\begin{definition}
An $\infty$-category $\calC$ is {\it $\kappa$-filtered} if for any
diagram $F: \calI \rightarrow \calC$ where $\calI$ is
$\kappa$-small, there exists a natural transformation from $F$ to
a constant functor. We say that $\calC$ is {\it filtered} if it is
$\omega$-filtered. A diagram $\calI \rightarrow \calC$ is called
{\it $\kappa$-filtered} if $\calI$ is $\kappa$-filtered, and {\it
filtered} if $\calI$ is $\omega$-filtered.
\end{definition}

As the definition suggests, the most important case occurs when
$\kappa = \omega$.

\begin{remark}
Filtered $\infty$-categories are the natural generalization of
filtered categories, which are in turn a mild generalization of
directed partially ordered sets. Recall that a partially ordered
set $P$ is {\it directed} if every finite subset of $P$ has some
upper bound in $P$. One frequently encounters diagrams indexed by
directed partially ordered sets, for example by the directed
partially ordered set $\Z_{\geq 0} = \{0, 1, \ldots \}$ of natural
numbers: colimits $$\colim_{n \rightarrow \infty} X_n$$ of such
diagrams are among the most familiar constructions in mathematics.

In classical category theory, it is convenient to consider not
only diagrams indexed by directed partially ordered sets, but also
more generally diagram indexed by filtered categories. A filtered
category is defined to be a category $\calC$ satisfying the
following conditions:

\begin{itemize}
\item There exists at least one object in $\calC$.
\item Given any two objects $X,Y \in \calC$, there exists a third
object $Z$ and morphisms $X \rightarrow Z$, $Y \rightarrow Z$.
\item Given any two morphisms $f,g: X \rightarrow Y$ in $\calC$,
there exists a morphism $h: Y \rightarrow Z$ such that $h \circ f
= h \circ g$.
\end{itemize}

The first two conditions are analogous to the requirement that any
finite piece of $\calC$ has an ``upper bound'', while the last
condition guarantees that the upper bound is uniquely determined
in some asymptotic sense.

One can give a similar definition of filtered $\infty$-categories,
but one must assume a more general version of the third condition.
A pair of maps $f,g: X \rightarrow Y$ may be regarded as a map
$S^0 \rightarrow \Hom_{\calC}(X,Y)$, where $S^0$ denotes the
zero-sphere. The third condition asserts that there exists a map
$Y \rightarrow Z$ such that the induced map $S^0 \rightarrow
\Hom_{\calC}(X,Z)$ extends to a map $D^1 \rightarrow
\Hom_{\calC}(X,Z)$, where $D^1$ is the $1$-disk. To arrive at the
definition of a filtered $\infty$-category, we must require the
existence of $Y \rightarrow Z$ not only for maps $S^0 \rightarrow
\Hom_{\calC}(X,Y)$, but for any map $S^n \rightarrow
\Hom_{\calC}(X,Y)$ (and the requirement should be that the induced
map extends over the disk $D^{n+1}$). With this revision, the
above conditions are equivalent to our definition of
$\omega$-filtered $\infty$-categories.

The reader should not be intimidated by the apparent generality of
the notion of a $\kappa$-filtered colimit. Given any
$\kappa$-filtered diagram $\calI \rightarrow \calC$, one can
always find a diagram $\calJ \rightarrow \calC$ indexed by a
partially ordered set $\calJ$ which is ``equivalent'' in the sense
that they have the same colimit (provided that either colimit
exists), and $\calJ$ is $\kappa$-filtered in the sense that every
subset of $\calJ$ having cardinality $< \kappa$ has an upper bound
in $\calJ$. We sketch a proof for this in Appendix
\ref{appendixdiagram}.
\end{remark}

\begin{definition}
We will call an $\infty$-category $\calC$ {\it $\kappa$-closed} if
every $\kappa$-filtered diagram has a colimit in $\calC$. A
functor between $\kappa$-closed $\infty$-categories will be called
$\kappa$-continuous if it commutes with the formation of
$\kappa$-filtered colimits. Finally, we will call an object $C$ in
a $\kappa$-closed $\infty$-category $\calC$ {\it $\kappa$-compact}
if the functor $\Hom_{\calC}(C, \bigdot)$ is a $\kappa$-continuous
functor $\calC \rightarrow \SSet$. When $\kappa = \omega$, we will
abbreviate by simply referring to {\it continuous} functors and
{\it compact} objects.
\end{definition}

\begin{definition}
If $\calC$ is any $\kappa$-closed $\infty$-category, we will let
$\calC_{\kappa}$ denote the full subcategory consisting of
$\kappa$-compact objects.
\end{definition}

\begin{remark}
Any $\kappa$-closed category ($\kappa$-continuous morphism,
$\kappa$-compact object) is also $\kappa'$-closed
($\kappa'$-continuous, $\kappa'$-compact) for any $\kappa' \geq
\kappa$.
\end{remark}

\begin{remark}
There are a number of reasons to place emphasis on
$\kappa$-filtered colimits (particularly in the case where $\kappa
= \omega$). They exist more often and tend to be more readily
computable than colimits in general. For example, the category of
groups admits all colimits, but it is usually difficult to give an
explicit description of a colimit (typically this involves solving
a word problem). However, filtered colimits of groups are easy to
describe, because the formation of filtered colimits is compatible
with passage to the underlying set of a group. In other words, the
forgetful functor $F$, which assigns to each group its underlying
set, is continuous.
\end{remark}

We now come to a key construction. Let $\calC$ be a
$\infty$-category and $\kappa$ a regular cardinal. We will define
a new $\infty$-category $\Ind_{\kappa}(\calC)$ as follows. The
objects of $\Ind_{\kappa}(\calC)$ are small, $\kappa$-filtered
diagrams $\calI \rightarrow \calC$. We imagine that an object of
$\Ind_{\kappa}(\calC)$ is a {\em formal} $\kappa$-filtered colimit
of objects in $\calC$, and use the notation
$$`` \colim_{D \in \calI} F(D)"$$ to denote the object
corresponding to a diagram $F: \calI \rightarrow \calC$. The
morphism spaces are given by

$$\Hom_{\Ind_{\kappa}(\calC)}( `` \colim_{D \in \calI} F(D) ", ``
\colim_{D' \in \calI'} F'(D')") = \lim_{D \in \calI} \colim_{D'
\in \calI'} \Hom_{\calC}(FD, F'D')$$

When $\kappa=\omega$, we will sometimes write $\Ind(\calC)$
instead of $\Ind_{\kappa}(\calC)$. For a more detailed discussion
of the homotopy theory of filtered diagrams, from the point of
view of model categories, we refer the reader to
\cite{filteredhomotopy}.

\begin{remark}
If $\calC$ is an ordinary category, then so is
$\Ind_{\kappa}(\calC)$. Many categories have the form
$\Ind(\calC)$. For example, the category of groups is equivalent
to $\Ind(\calC)$, where $\calC$ is the category of finitely
presented groups. A similar remark applies to any other sort of
(finitary) algebraic structure.
\end{remark}

There is a functor $\calC \rightarrow \Ind_{\kappa}(\calC)$, which
carries an object $C \in \calC$ to the constant diagram $\ast
\rightarrow \calC$ having value $C$, where $\ast$ denotes a
category having a single object with a contractible space of
endomorphisms. This functor is fully faithful, and identifies
$\calC$ with a full subcategory of $\Ind_{\kappa}(\calC)$. One can
easily verify that $\Ind_{\kappa}(\calC)$ is $\kappa$-closed, and
that the essential image of $\calC$ consists of $\kappa$-compact
objects. In fact, $\Ind_{\kappa}(\calC)$ is universal with respect
to these properties:

\begin{proposition}
Let $\calC$ and $\calC'$ be $\infty$-categories, where $\calC'$ is
$\kappa$-closed. Then restriction induces an equivalence of
$\infty$-categories
$$\Fun^{\kappa}( \Ind_{\kappa}(\calC), \calC') \rightarrow
\Fun(\calC, \calC')$$ where the left hand side denotes the
$\infty$-category of $\kappa$-continuous functors. In other words,
any functor $F: \calC \rightarrow \calC'$ extends uniquely to a
$\kappa$-continuous functor $\widetilde{F}: \Ind_{\kappa}(\calC)
\rightarrow \calC'$. If $F$ is fully faithful and its image
consists of $\kappa$-compact objects, then $\widetilde{F}$ is also
fully faithful.
\end{proposition}

\begin{proof}
This follows more or less immediately from the construction.
\end{proof}

From this, we can easily deduce the following:

\begin{proposition}\label{clear}
Let $\calC$ be an $\infty$-category. The following conditions are
equivalent:

\begin{itemize}
\item The $\infty$-category $\calC$ is equivalent to $\Ind_{\kappa} \calC_0$ for some
small $\infty$-category $\calC_0$.
\item The $\infty$-category $\calC$ is $\kappa$-closed,
$\calC_{\kappa}$ is essentially small, and the functor
$\Ind_{\kappa}(\calC_{\kappa}) \rightarrow \calC$ is an
equivalence of categories.
\item The $\infty$-category $\calC$ is $\kappa$-closed, and has a
small subcategory $\calC_0$ consisting of $\kappa$-compact objects
such that every object of $\calC$ can be obtained as a
$\kappa$-filtered colimit of objects in $\calC_0$.
\end{itemize}
\end{proposition}

An $\infty$-category $\calC$ satisfying the equivalent conditions
of Proposition \ref{clear} will be called {\it
$\kappa$-accessible}.

The only difficulty in proving Proposition \ref{clear} is in
verifying that if $\calC_0$ is small, then $\Ind_\kappa(\calC_0)$
has only a set of $\kappa$-compact objects, up to equivalence. It
is tempting to guess that any such object must be equivalent to an
object of $\calC_0$. The following example shows that this is not
necessarily the case.

\begin{example}
Let $R$ be a ring, and let $\calC_0$ denote the (ordinary)
category of finitely generated free $R$-modules. Then $\calC =
\Ind(\calC_0)$ is equivalent to the category of flat $R$-modules
(by Lazard's theorem; see for example the appendix of
\cite{lazard}). The compact objects of $\calC$ are precisely the
finitely generated projective $R$-modules, which need not be free.
\end{example}

However, a $\kappa$-compact object $C$ in $\Ind_{\kappa}(\calC_0)$
is not far from being an object of $\calC_0$: it follows readily
from the definition that $C$ is a {\it retract} of an object of
$\calC_0$, meaning that there are maps
$$ C \stackrel{f}{\rightarrow} D \stackrel{g}{\rightarrow} C$$
with $D \in \calC_0$ and $g \circ f$ homotopic to the identity of
$C$. In this case, we may recover $C$ given the object $D$ and
appropriate additional data. Since $g \circ f$ is homotopic to the
identity, we deduce easily that the morphism $r = f \circ g$ is
idempotent in the homotopy category $h\calC$. Moreover, in
$h\calC$, the object $C$ is the equalizer of $r$ and the identity
morphism of $D$. This shows that there are only a bounded number
of possibilities for $C$, and proves Proposition \ref{clear}.

\begin{definition}
An $\infty$-category $\calC$ is {\it accessible} if it is
$\kappa$-accessible for some regular cardinal $\kappa$. A functor
between accessible $\infty$-categories is {\it accessible} if it
is $\kappa$-continuous for some $\kappa$ (and therefore for all
sufficiently large $\kappa$).
\end{definition}

It is not necessarily true that a $\kappa$-accessible
$\infty$-category $\calC$ is $\kappa'$ accessible for all $\kappa'
> \kappa$. However, it is true that $\calC$ is
$\kappa'$-accessible for many other cardinals $\kappa'$. Let us
write $\kappa' \gg \kappa$ if $\tau^{\kappa} < \kappa'$ for any
$\tau < \kappa'$. Note that there exist arbitrarily large regular
cardinals $\kappa'$ with $\kappa' \gg \kappa$: for example, one
may take $\kappa'$ to be the successor of any cardinal having the
form $\tau^{\kappa}$.

\begin{lemma}\label{estimate}
If $\kappa' \gg \kappa$, then any $\kappa'$-filtered partially
ordered set $\calI$ may be written as a union of $\kappa$-filtered
subsets having size $< \kappa'$. Moreover, the family of such
subsets is $\kappa'$-filtered.
\end{lemma}

\begin{proof}
It will suffice to show that every subset of $S \subseteq
\calI$ having cardinality $< \kappa'$ can be included in a larger
subset $S'$, such that $|S'| < \kappa'$, but $S'$ is
$\kappa$-filtered.

We define a transfinite sequence of subsets $S_{\alpha} \subseteq
\calI$ by induction. Let $S_{0} = S$, and when $\lambda$ is a
limit ordinal we let $S_{\lambda} = \bigcup_{\alpha < \lambda}
S_{\alpha}$. Finally, we let $S_{\alpha + 1}$ denote a set which
is obtained from $S_{\alpha}$ by adjoining an upper bound for
every subset of $S_{\alpha}$ having size $< \kappa$ (which exists
because $\calI$ is $\kappa'$-filtered). It follows from the
assumption $\kappa' \gg \kappa$ that if $S_{\alpha}$ has size $<
\kappa'$, then so does $S_{\alpha + 1}$. Since $\kappa'$ is
regular, we deduce easily by induction that $|S_{\alpha}| <
\kappa'$ for all $\alpha < \kappa'$. It is easy to check that the
set $S' = S_{\kappa}$ has the desired properties.
\end{proof}

We now show that accessible categories enjoy several pleasant
boundedness properties:

\begin{proposition}\label{firss}
Let $\calC$ be an accessible category. Then every object of
$\calC$ is $\kappa$-compact for some regular cardinal $\kappa$.
For any regular cardinal $\kappa$, the full subcategory of $\calC$
consisting of $\kappa$-compact cardinals is essentially small.
\end{proposition}

\begin{proof}
Suppose $\calC$ is $\tau$-accessible. Then every object of $\calC$
may be written as a $\tau$-filtered colimit of $\tau$-compact
objects. But it is easy to see that for $\kappa \geq \tau$, the
colimit of a $\kappa$-small diagram of $\tau$-compact object is
$\kappa$-compact.

Now fix a regular cardinal $\kappa \gg \tau$, and let $C$ be a
$\kappa$-compact object of $\calC$. Then $C$ may be written as a
$\tau$-filtered colimit of a diagram $\calI \rightarrow
\calC_{\tau}$. Without loss of generality, we may assume that
$\calI$ is actually the $\infty$-category associated to a
partially ordered set. By Lemma \ref{estimate}, we may write
$\calI$ as a $\kappa$-filtered union of $\tau$-filtered,
$\kappa$-small subsets $\calI_{\alpha}$. Let $C_{\alpha}$ denote
the colimit of the diagram indexed by $\calI_{\alpha}$; then $C =
\colim_{\alpha} C_{\alpha}$ where {\em this} colimit is
$\kappa$-filtered. Consequently, $C$ is a retract of some
$C_{\alpha}$. Since there are only a bounded number of
possibilites for $\calI_{\alpha}$ and its functor to
$\calC_{\tau}$, we see that there are only a bounded number of
possibilities for $C$, up to equivalence.
\end{proof}

\begin{remark}
It is not really necessary to reduce to the case where $\calI$ is
a partially ordered set. We could instead extend Lemma
\ref{estimate} to $\infty$-categories: the difficulties are
primarily notational.
\end{remark}

\begin{proposition}
Let $\calC$ be a $\kappa$-accessible $\infty$-category. Then
$\calC$ is $\kappa'$-accessible for any $\kappa' \gg \kappa$.
\end{proposition}

\begin{proof}
It is clear that $\calC$ is $\kappa'$-closed. Proposition
\ref{firss} implies that $\calC_{\kappa'}$ is essentially small.
To complete the proof, it suffices to show that any object $C \in
\calC$ can be obtained as the colimit of a $\kappa'$-filtered
diagram in $\calC_{\kappa'}$. By assumption, $C$ can be
represented as the colimit of a $\kappa$-filtered diagram $F:
\calI \rightarrow \calC_{\kappa}$. Without loss of generality, we
may assume that $\calI$ is a partially ordered set.

By Lemma \ref{estimate}, we may write $\calI$ as a
$\kappa'$-filtered union of $\kappa'$-small, $\kappa$-filtered
partially ordered sets $\calI_{\alpha}$. Let $F_{\alpha} =
F|\calI_{\alpha}$. Then we get
$$C = \colim F = \colim_{\alpha} (\colim F_{\alpha})$$

Since each $\calI_{\alpha}$ is $\kappa'$-small, the colimit of
each $F_{\alpha}$ is $\kappa'$-compact. It follows that $C$ is a
$\kappa'$-filtered colimit of $\kappa'$-small objects, as desired.
\end{proof}

We now show that accessible functors between accessible
$\infty$-categories are themselves determined by a bounded amount
of data.

\begin{proposition}
Let $F: \calC \rightarrow \calC'$ be an accessible functor between
accessible $\infty$-categories. Then there exists a regular
cardinal $\kappa$ such that $\calC$ and $\calC'$ are
$\kappa$-accessible, and $F$ is the unique $\kappa$-continuous
extension of a functor $F_{\kappa}: \calC_{\kappa} \rightarrow
\calC'_{\kappa}$.
\end{proposition}

\begin{proof}
Choose $\tau$ so large that $F$ is $\tau$-continuous. Enlarging
$\tau$ if necessary, we may assume that $\calC$ and $\calC'$ are
$\tau$-accessible. Then $\calC$ has only a bounded number of
$\tau$-compact objects. By Proposition \ref{firss}, the images of
these objects under $F$ are all $\tau'$-compact for some regular
cardinal $\tau' \geq \tau$. We may now choose $\kappa$ to be any
regular cardinal such that $\kappa \gg \tau'$.

Any $\kappa$-compact object of $\calC$ can be written as a
$\kappa$-small, $\tau$-filtered colimit of $\tau$-compact objects.
Since $F$ preserves $\tau$-filtered colimits, the image of a
$\kappa$-compact of $\calC$ is a $\kappa$-small, $\tau$-filtered
colimit of $\tau'$-compact objects, hence $\kappa$-compact. Thus
the restriction of $F$ maps $\calC_{\kappa}$ into
$\calC'_{\kappa}$. To complete the proof, it suffices to note that
the $\tau$-continuity of $F$ implies $\kappa$-continuity since
$\kappa > \tau$.
\end{proof}

The above proposition shows that there is an honest
$\infty$-category of accessible functors from $\calC$ to $\calC'$,
for any accessible $\infty$-categories $\calC$ and $\calC'$.
Indeed, this $\infty$-category is just a (filtered) union of the
essentially small $\infty$-categories $\Hom(\calC_{\kappa},
\calC'_{\kappa})$, where $\kappa$ ranges over regular cardinals
such that $\calC$ and $\calC'$ are $\kappa$-accessible. We will
denote this $\infty$-category by $\Acc(\calC, \calC')$. Any
composition of accessible functors is accessible, so that we
obtain a (large) $(\infty,2)$-category $\Acc$ consisting of
accessible $\infty$-categories and accessible functors.

\begin{proposition}\label{access}
The $(\infty,2)$-category $\Acc$ has all $(\infty,2)$-categorical
limits, which are simply given by taking the limits of the
underlying $\infty$-categories.
\end{proposition}

In other words, a limit of accessible $\infty$-categories (and
accessible functors) is accessible.

\begin{proof}
We may choose a regular cardinal $\kappa$ such that all of the
accessible categories are $\kappa$-accessible and all of the
relevant functors between them are $\kappa$-continuous. Let
$\calC$ denote the limit in question. It is clear that $\calC$ is
$\kappa$-closed: one may compute $\kappa$-filtered colimits in
$\calC$ componentwise (since all of the relevant functors are
$\kappa$-continuous). To complete the proof, we must show that
(possibly after enlarging $\kappa$), the limit is generated by a
set of $\kappa$-compact objects. The proof of this point involves
tedious cardinality estimates, and will be omitted. The reader may
consult \cite{adamek} for a proof in the $1$-categorical case,
which contains all of the relevant ideas.
\end{proof}

\begin{example}\label{bbqq}
Let $\calC$ be an accessible $\infty$-category, and let $X$ be an
object of $\calC$. Then $\calC_{/X}$ is an accessible
$\infty$-category, since it is a lax pullback of accessible
$\infty$-categories (see Example \ref{slicesarelax}).
\end{example}

We conclude this section with the following useful observation:

\begin{proposition}\label{adjoints}
Let $G: \calC \rightarrow \calC'$ be a functor between accessible
$\infty$-categories. If $G$ admits an adjoint, then $G$ is
accessible.
\end{proposition}

\begin{proof}
If $G$ is a left adjoint, then $G$ commutes with all colimits and
is therefore $\omega$-continuous. Let us therefore assume that $G$
is right adjoint to some functor $F$.

Choose $\kappa$ so that $\calC$ and $\calC'$ are
$\kappa$-accessible. Now choose $\kappa' \gg \kappa$ large enough
that $FC'$ is $\kappa'$-compact for all $C' \in \calC_{\kappa}$.
We claim that $G$ is $\kappa'$-continuous.

Let $\pi: \calI \rightarrow \calC$ be any $\kappa'$-filtered
diagram having colimit $C$. We must show that $GC$ is the colimit
of the diagram $G \circ \pi$.

We need to show that $\Hom_{\calC'}(C',GC)$ is equivalent to
$\Hom_{\calC'}(C', \colim_{C_\alpha \in \calI} \{GC_\alpha \})$
for all objects $C' \in \calC'$. It suffices to check this for
$C'$ in $\calC'_{\kappa}$, since every object in $\calC'$ is a
colimit of objects in $\calC'_{\kappa}$. Then

$$\begin{array}{ccl}
\Hom_{\calC'}(C', \colim_{C_{\alpha} \in \calI} \{GC_{\alpha} \})
& = & \colim_{C_{\alpha} \in \calI} \Hom_{\calC'}(C', GC_{\alpha})
\\
& = &\colim_{C_{\alpha} \in \calI} \Hom_{\calC}(FC',C_{\alpha})
\\
&= &\Hom_{\calC}(FC', C) \\
&= &\Hom_{\calC'}(C',GC),\\
\end{array}$$ as required.
\end{proof}

\subsection{Presentable $\infty$-Categories}

In this section, we discuss the notion of a {\it presentable
$\infty$-category}. Roughly speaking, this is an $\infty$-category
with a set of small generators that admits arbitrary colimits.
This includes, for example, any ordinary category with a set of
small generators that admits arbitrary colimits: such categories
are called locally presentable in \cite{adamek}. The definition,
in the $\infty$-categorical setting, was given by Carlos Simpson
in \cite{simpson}, who calls them $\infty$-pretopoi.

\begin{proposition}\label{pretop}
Let $\calC$ be an $\infty$-category. Then the following are
equivalent:

\begin{enumerate}
\item The $\infty$-category $\calC$ is accessible and admits
arbitrary colimits.

\item For all sufficiently large regular cardinals $\kappa$, the
$\infty$-category $\calC_{\kappa}$ is essentially small, admits
colimits for all $\kappa$-small diagrams, and the induced
$\Ind_{\kappa}(\calC_{\kappa}) \rightarrow \calC$ is an
equivalence.

\item There exists a regular cardinal $\kappa$ such that the
$\infty$-category $\calC_{\kappa}$ is essentially small, admits
$\kappa$-small colimits, and the induced functor
$\Ind_{\kappa}(\calC_{\kappa}) \rightarrow \calC$ is an
equivalence.

\item There is a small $\infty$-category $\calC'$ which
admits all $\kappa$-small colimits such that $\calC$ is equivalent
to $\Ind_{\kappa}(\calC')$.

\item The $\infty$-category $\calC$ admits all colimits, and there exists
a cardinal $\kappa$ and a set $S$ of $\kappa$-compact objects of
$\calC$ such that every object of $\calC$ is a colimit of objects
in $S$.
\end{enumerate}

\end{proposition}

\begin{proof}
Since $\kappa$-small colimits of $\kappa$-compact objects are
$\kappa$-compact when they exist, it is clear that the $(1)$
implies $(2)$ for every cardinal $\kappa$ such that $\calC$ is
$\kappa$-accessible. To complete the proof of this step, it
suffices to show that under these hypotheses, the
$\kappa$-accessibility of $\calC$ implies $\kappa'$-accessibility
for all $\kappa' > \kappa$. This follows from the following
observation: if $C \in \calC$ is the colimit of a
$\kappa$-filtered diagram $\calI \rightarrow \calC$, then we may
replace $\calI$ by a partially ordered set, which is the union of
its subsets $\calI_{\alpha}$ having size $< \kappa'$. Then $C$ is
the $\kappa'$-filtered colimit of objects $C_{\alpha}$, where
$C_{\alpha}$ denotes the colimit of the induced diagram $\calI
\rightarrow \calC$ (which might not be filtered). Since
$C_{\alpha}$ is a $\kappa'$-small colimit of $\kappa$-compact
objects, it is $\kappa'$-compact.

The implications $(2) \Rightarrow (3) \Rightarrow (4) \Rightarrow
(5)$ are obvious. We will complete the proof by showing that $(5)
\Rightarrow (1)$.  Assume that there exists a regular cardinal
$\kappa$ and a set $S$ of $\kappa$-compact objects of $\calC$ such
that every object of $\calC$ is a colimit of objects in $S$. We
first claim that $\calC_{\kappa}$ is essentially small: any
$\kappa$-compact object $C \in \calC$ is a colimit of a diagram
involving only objects of $S$, hence a $\kappa$-filtered colimit
of the colimits of $\kappa$-small diagrams in $S$. The
$\kappa$-compactness of $C$ implies that $C$ is a retract of a
$\kappa$-small colimit of objects in $S$, and thus there are only
a bounded number of possibilities for $C$ up to equivalence. To
complete the proof, it suffices to show that every object of
$\calC$ is a $\kappa$-filtered colimit of objects of
$\calC_{\kappa}$. But this follows immediately, since the colimit
of any $\kappa$-small diagram in $S$ belongs to $\calC_{\kappa}$.
\end{proof}

An $\infty$-category $\calC$ satisfying the equivalent conditions
of Proposition \ref{pretop} is called {\it presentable}.

\begin{remark}\label{tensored}
Let $\calC$ be a presentable $\infty$-category. Then $\calC$ is
``tensored over $\SSet$'' in the following sense: there exists a
functor $\calC \times \SSet \rightarrow \calC$, which we shall
denote by $\otimes$, having the property that
$$\Hom_{\calC}( C \otimes S, C') = \Hom_{\SSet}(S,
\Hom_{\calC}(C,C'))$$ In order to prove this, we note that the $C
\otimes S$ is canonically determined by its universal property and
that its formation is compatible with colimits in $S$. Any space
$S \in \SSet$ may be obtained as a colimit of points (if we regard
$S$ as an $\infty$-groupoid, then $S$ is the colimit of the
constant diagram $S \rightarrow \SSet$ having the value $\ast$).
Thus, the existence of $C \otimes S$ can be deduced from the
existence when $S$ is a single point, in which case we may take $C
\otimes S = C$.
\end{remark}

The key fact about presentable $\infty$-categories which makes
them very pleasant to work with is the following representability
criterion:

\begin{proposition}\label{representable}

Let $\calC$ be a presentable $\infty$-category, and let $\calF$ be
a prestack on $\calC$. Suppose that $\calF$ carries colimits in
$\calC$ into limits in $\SSet$. Then $\calF$ is representable.
\end{proposition}

\begin{proof}
Choose a regular cardinal $\kappa$ such that $\calC$ is
$\kappa$-accessible. Consider the restriction of $\calF$ to
$\calC_{\kappa}$: we may view this as an $\infty$-category
``fibered in spaces'' over $\calC_{\kappa}$. More precisely, we
may form a new $\infty$-category $\widetilde{\calC}$ whose objects
consist of pairs $C \in \calC_{\kappa}$, $\eta \in \calF(C)$, with
$\Hom_{\widetilde{\calC}}( (C,\eta), (C',\eta') )$ given by the
fiber of $p: \Hom_{\calC}(C,C') \rightarrow \calF(C)$ over the
point $\eta$, where $p$ is induced by pulling back $\eta' \in
\calF(C')$. Note that $\widetilde{\calC}$ admits all
$\kappa$-small colimits, and in particular is $\kappa$-filtered.

There exists a functor $\pi: \widetilde{\calC} \rightarrow
\calC_{\kappa}$. Regard $\pi$ as a diagram in $\calC$, and let
$F_{\kappa}$ denote its colimit. Then, by the compatibility of
$\calF$ with colimits, we obtain a natural ``universal element''
$\eta_{\kappa} \in \calF(F_{\kappa})$. We would like to show that
this exhibits $F_{\kappa}$ as representing the presheaf $\calF$.

We first prove the $F_{\kappa}$ is a good approximation to $\calF$
in the following sense: given any morphism $f: X \rightarrow Y$
between $\kappa$-compact objects of $\calC$, the natural map
$\Hom_{\calC}(Y, F_{\kappa}) \rightarrow \calF(Y)
\times_{\calF(X)} \Hom_{\calC}(X, F_{\kappa})$ is surjective on
$\pi_0$. Indeed, suppose we are given a compatible triple
$(X,\eta_X)$, $(Y,\eta_Y)$, and $p: X \rightarrow F_{\kappa}$,
where $\eta_Y \in \calF(Y)$, $\eta_X = f^{\ast} \eta_Y \in
\calF(X)$, and we are given a homotopy between
$p^{\ast}(\eta_{\kappa})$ and $\eta_X$. Since $X$ is
$\kappa$-compact, $p$ factors through some pair $(Z,\eta_Z) \in
\widetilde{\calC}$. Let $W = Y \coprod_X Z$. Since $\calF$ is
compatible with pushouts, the homotopy between $\eta_Z|X$ and
$\eta_Y|X$ gives us $\eta_W \in \calF(W)$ with $\eta_W|Z=\eta_Z$
and $\eta_W|Y = \eta_Y$. Then $(W,\eta_W) \in \widetilde{\calC}$,
and the natural map $Y \rightarrow W$ gives rise to a point of
$\Hom_{\calC}(Y, F_{\kappa})$. One readily checks that this point
has the desired properties.

Now let us show that $F_{\kappa}$ represents $\calF$. In other
words, we want to show that the natural map
$$\Hom_{\calC}(Y,F_{\kappa}) \rightarrow \calF(Y)$$
is a homotopy equivalence for every $Y$. Since both sides are
compatible with colimits in the variable $Y$, we may assume that
$Y$ is $\kappa$-compact. By Whitehead's theorem, it will suffice
to show that the natural map
$$ \Hom_{\calC}(Y,F_{\kappa}) \rightarrow \calF(Y) \times_{ [S^n,
\calF(Y)] } [S^n, \Hom_{\calC}(Y, F_{\kappa})]$$ is surjective on
$\pi_0$, where $[S^n, Z]$ denotes the ``free loop space''
consisting of maps from an $n$-sphere $S^n$ into $Z$. This reduces
to a special case of what we proved above if we take $X = Y
\otimes S^n$.
\end{proof}

The representability criterion of Proposition \ref{representable}
has many consequences, as we now demonstrate.

\begin{corollary}
A presentable $\infty$-category admits arbitrary limits.
\end{corollary}

\begin{proof}
Any limit of representable prestacks carries colimits into limits,
and is therefore representable by Proposition \ref{representable}.
\end{proof}

\begin{remark}
Now that we know that $\calC$ has arbitrary limits, we can apply
an argument dual to that of Remark \ref{tensored} to show that
$\calC$ is {\it cotensored over $\SSet$}. In other words, for any
$C \in \calC$ and any $X \in \calS$, there exists an object $C^X
\in \calC$ and a natural equivalence $\Hom_{\calC}(\bigdot,C^X) =
\Hom_{\calS}(X, \Hom_{\calC}(\bigdot,C)$.
\end{remark}

From Proposition \ref{representable} we also get a version of the
adjoint functor theorem:

\begin{corollary}\label{adjointfunctor}
Let $F: \calC \rightarrow \calC'$ be a functor between presentable
$\infty$-categories. Then $F$ has a right adjoint if and only if
$F$ preserves all colimits.
\end{corollary}

\begin{proof}
Suppose we have a small diagram $\{X_i\}$ in $\calC$ with colimit
$X$. If $F$ has a right adjoint $G$, then $\Hom_{\calC'}(FX,Y) =
\Hom_{\calC}(X,GY) = \lim_{i} \Hom_{\calC}(X_i,GY) = \lim_i
\Hom_{\calC'}(FX_i,Y)$, so that $FX$ is a colimit of the induced
diagram $\{FX_i\}$ in $\calC'$. Conversely, if $F$ preserves all
colimits, then for each $Y \in \calC'$, we may define a prestack
$\calF_Y$ on $\calC$ by the equation $\calF_Y(X) =
\Hom_{\calC'}(FX,Y)$. The hypothesis on $F$ implies that
$\calF_{Y}$ carries colimits into limits, so $\calF_Y$ is
representable by an object $GY \in \calC$. The universal property
enjoyed by $GY$ automatically ensures that $G$ is a functor
$\calC' \rightarrow \calC$, left adjoint to $F$.
\end{proof}

If $\calC$ and $\calC'$ are presentable, we let $\Pre(\calC,
\calC')$ denote the $\infty$-category of colimit-preserving
functors from $\calC$ to $\calC'$. These give the appropriate
notion of ``morphism'' between presentable $\infty$-categories.

The collection of colimit functors between presentable
$\infty$-categories is stable under compositions. We may summarize
the situation informally by saying that there is an
$(\infty,2)$-category $\Pre$ of presentable $\infty$-categories,
with colimit-preserving functors as morphisms. This
$(\infty,2)$-category admits arbitrary $(\infty,2)$-categorical
limits, which are just $(\infty,2)$-categorical limits of the
underlying $\infty$-categories (such limits are accessible by
Proposition \ref{access}, and colimits may be computed levelwise).
This $(\infty,2)$-category also admits arbitrary
$(\infty,2)$-categorical colimits, which may be constructed using
generators and relations (see \cite{adamek} for a discussion in
the $1$-categorical case), but we shall not need this.

\begin{remark}\label{slicepresentable}
If $\calC$ is a presentable $\infty$-category, then any slice
$\calC_{/X}$ is also presentable. Indeed, $\calC_{/X}$ is
accessible (see Example \ref{bbqq}) and has all colimits (which
are also colimits in $\calC$).
\end{remark}

\begin{remark}
Let $\calC$ be a small $\infty$-category, and let $\calP =
\SSet^{\calC^{op}}$ denote the $\infty$-category of stacks on
$\calC$. Then $\calP$ is presentable: in fact, it is a limit of
copies of $\SSet$ indexed by $\calC^{op}$ (one can also argue
directly). But $\calP$ also has a dual description: it is a {\em
colimit} of copies of $\SSet$ indexed by $\calC$. In other words,
$\calP$ is the {\it free presentable $\infty$-category generated
by $\calC$}, in the sense that for any $\infty$-pretopos $\calC'$
we have an equivalence of $\infty$-categories
$$\pi: \Pre(\calP, \calC') \rightarrow \Fun(\calC, \calC')$$
where the left hand side denotes the $\infty$-category of
right-exact continuous functors from $\calP$ to $\calC'$. The
functor $\pi$ is given by restriction, and it has a homotopy
inverse which carries a functor $F: \calC \rightarrow \calC'$ to
$\widetilde{F}: \calP \rightarrow \calC'$, where
$\widetilde{F}(\calF)$ is the homotopy colimit of the diagram $F
\circ p: \calI \rightarrow \calC'$, where $\calI$ is the
$\infty$-category fibered over $\calC$ in groupoids defined by
$\calF$ and $p: \calI \rightarrow \calC$ is the canonical
projection.
\end{remark}

We conclude by showing that the $(\infty,2)$-category of
presentable $\infty$-categories admits an internal $\Hom$-functor.
The proof uses material from the next section.

\begin{proposition}
Let $\calC$ and $\calC'$ be presentable $\infty$-categories. The
$\infty$-category $\Pre(\calC,\calC')$ of colimit-preserving
functors from $\calC$ to $\calC'$ is presentable.
\end{proposition}

\begin{proof}
We use Proposition \ref{simpsongiraud}, to be proved in the next
section, which asserts that $\calC$ is a localization of
$\SSet^{\calC_0^{op}}$ for some small $\infty$-category $\calC_0$.
This immediately implies that $\Pre(\calC,\calC')$ is a
localization of $\Pre(\SSet^{\calC_0^{op}}, \calC')$. Using
Proposition \ref{simpsongiraud} again, it suffices to consider the
case where $\calC= \SSet^{\calC_0^{op}}$. In this case, we deduce
that $\Pre(\calC,\calC') = \Fun( \calC_0, \calC')$, which is a
limit (indexed by $\calC_0$) of $\infty$-categories having the
form $\Hom( \ast, \calC') = \calC'$. Since $\calC'$ is an
presentable, the limit is also presentable.
\end{proof}

\subsection{Localization}

In this section, we discuss {\it localizations} of a presentable
$\infty$-category $\calC$. These ideas are due to Bousfield (see
for example \cite{bousfield}). Suppose we have a collection $S$ of
morphisms in $\calC$ which we would like to invert. In other
words, we seek an presentable $\infty$-category $S^{-1}\calC$
equipped with a colimit-preserving functor $\pi: \calC \rightarrow
S^{-1} \calC$ which is ``universal'' with respect to the fact that
every morphism in $S$ becomes an equivalence in $S^{-1} \calC$. It
is an observation of Bousfield that we can often find $S^{-1}
\calC$ {\it inside of $\calC$}, as the set of ``local'' objects of
$\calC$.

\begin{example}\label{excom}
Let $\calC$ be the (ordinary) category of abelian groups, $p$ a
prime number, and let $S$ denote the collection of morphisms $f$
whose kernel and cokernel consist entirely of $p$-power torsion. A
morphism $f$ lies in $S$ if and only if it induces an isomorphism
after inverting the prime number $p$. In this case, we may
identify $S^{-1} \calC$ with the category of abelian groups which
are {\em uniquely $p$-divisible}.
\end{example}

In Example \ref{excom}, the natural functor $\calC \rightarrow
S^{-1} \calC$ is actually left adjoint to an inclusion functor.
Let us therefore begin by examining functors which are left
adjoint to inclusions.

If $\calC_0 \subseteq \calC$ is a full subcategory of an
$\infty$-category, and $L$ is a left adjoint to the inclusion
$\calC_0 \rightarrow \calC$, then we may also regard $L$ as a
functor from $\calC$ to itself. The adjunction gives us a natural
transformation $\alpha: \id_{\calC} \rightarrow L$. From the pair
$(L, \alpha)$ we can recover $\calC_0$ as the full subcategory of
$\calC$ consisting of those objects $C \in \calC$ such that
$\alpha(C): C \rightarrow LC$ is an equivalence. Conversely, if we
begin with a functor $L: \calC \rightarrow \calC$ and a natural
transformation $\alpha: \id_{\calC} \rightarrow L$, and {\em
define} $\calC_0$ as above, then $L$ is left adjoint to the
inclusion $\calC_0 \rightarrow \calC$ if and only if $(L,\alpha)$
is {\em idempotent} in the following sense:

\begin{itemize}
\item For any $C \in \calC$, the morphisms $$L(\alpha(C)),
\alpha(LC): LC \rightarrow LLC$$ are equivalences which are
homotopic to one another.
\end{itemize}

In this situation, we shall say that $L$ is an {\it idempotent}
functor $\calC \rightarrow \calC$ (the transformation $\alpha$
being understood).

\begin{proposition}
Suppose that $\calC$ is an accessible category and that $L: \calC
\rightarrow \calC$ is idempotent. Then the induced functor $L_0:
\calC \rightarrow L\calC$ is accessible if and only if $L\calC$ is
an accessible category.
\end{proposition}

\begin{proof}
If $L \calC$ is accessible, then $L_0$ can be described as the
left adjoint of the inclusion $L \calC \subseteq \calC$ which is
accessible by Proposition \ref{adjoints}. For the converse,
suppose that $L_0$ is $\kappa$-continuous for some regular
cardinal $\kappa$. Without loss of generality we may enlarge
$\kappa$ so that $\calC$ is $\kappa$-accessible. It follows
immediately that $L \calC$ is stable under the formation of
$\kappa$-filtered colimits. Since $L_0$ is a left adjoint, it
preserves all colimits. Every object of $\calC$ is a
$\kappa$-filtered colimit of $\kappa$-compact objects of $\calC$;
hence every object of $L \calC$ is a $\kappa$-filtered colimit of
objects of the form $LC$, $C \in \calC_{\kappa}$. To complete the
proof, it suffices to note that every such $LC$ is
$\kappa$-compact in $L \calC$.
\end{proof}

An accessible idempotent functor in an accessible category will be
referred to as a {\it localization functor}. We will say that an
$\infty$-category $\calC_0$ is a {\it localization} of $\calC$ if
$\calC_0$ is equivalent to the essential image of some
localization functor $L: \calC \rightarrow \calC$.

We now give Simpson's characterization of presentable
$\infty$-categories. The proof can also be found in
\cite{simpson}, but we include it here for completeness.

\begin{proposition}\label{simpsongiraud}
Let $\calC$ be an $\infty$-category. The following conditions are
equivalent:
\begin{enumerate}
\item There exists a regular cardinal $\kappa$ such that the $\infty$-category $\calC$ is
$\kappa$-accessible and the Yoneda embedding $\calC \rightarrow
\calP =  \SSet^{\calC_{\kappa}^{op}}$ has a left adjoint.
\item The $\infty$-category $\calC$ is a localization of an
$\infty$-category of prestacks.
\item The $\infty$-category $\calC$ is a localization of a presentable $\infty$-category.
\item The $\infty$-category $\calC$ is presentable.
\end{enumerate}
\end{proposition}

\begin{proof}
If $\calC$ satisfies $(1)$, then we may regard it as a full
subcategory of $\calP$ via the Yoneda embedding. The left adjoint
to the Yoneda embedding is then viewed as an localization functor
$\calP \rightarrow \calP$ having $\calC$ as its essential image,
which proves $(2)$.

It is obvious that $(2)$ implies $(3)$. To see that $(3)$ implies
$(4)$, note that if $\calC$ is the essential image of some
localization functor defined on a larger presentable
$\infty$-category $\calC'$, then $\calC$ is accessible by
Proposition \ref{adjoints} and arbitrary colimits may be formed in
$\calC$ by forming those colimits in $\calC'$ and then applying
the functor $L$.

Finally, suppose that $\calC$ is presentable, and choose a regular
cardinal $\kappa$ such that $\calC$ is $\kappa$-accessible. Let
$\calP = \SSet^{\calC_{\kappa}^{op}}$ denote the $\infty$-category
of prestacks on $\calC_{\kappa}$. The functor $\calC_{\kappa}
\rightarrow \calC$ extends uniquely to a right-exact continuous
functor $\calP \rightarrow \calC$. One can easily check that this
functor is left adjoint to the Yoneda embedding, which proves
$(1)$.
\end{proof}

\begin{remark}
Carlos Simpson has also shown that the theory of presentable
$\infty$-categories is equivalent to that of {\it cofibrantly
generated model categories} (see \cite{simpson}). Since most of
the $\infty$-categories we shall meet are presentable, the subject
could also be phrased in the language of model categories. But we
shall not do this, since the restriction to presentable
$\infty$-categories seems unnatural and is often technically
inconvenient.
\end{remark}

We now proceed to discuss localizations of presentable
$\infty$-categories in more detail.

\begin{definition}
Let $\calC$ be an $\infty$-category and $S$ a collection of
morphisms of $\calC$. An object $C \in \calC$ is {\it $S$-local}
if, for any morphism $s: X \rightarrow Y$ in $S$, the induced
morphism $\Hom_{\calC}(Y,C) \rightarrow \Hom_{\calC}(X,C)$ is an
equivalence.

A morphism $f: X \rightarrow Y$ in $\calC$ is an {\it
$S$-equivalence} if the natural map $\Hom_{\calC}(Y,C) \rightarrow
\Hom_{\calC}(X,C)$ is a homotopy equivalence whenever $C$ is
$S$-local.
\end{definition}

Let $\rightarrow$ denote the (ordinary) category consisting of two
objects and a single morphism between them. For any
$\infty$-category $\calC$, we may form a functor $\infty$-category
$\calC^{\rightarrow}$ whose objects are the morphisms of $\calC$.
We note that $\calC^{\rightarrow}$ is presentable whenever $\calC$
is presentable.

\begin{definition}
Let $\calC$ be a presentable $\infty$-category. A collection $S$
of morphisms of $\calC$ will be called {\it saturated} if it
satisfies the following conditions:
\begin{itemize}
\item Every equivalence belongs to $S$.
\item If $f$ and $g$ are homotopic, then $f \in S$ if and only if
$g \in S$.
\item Suppose that $f: X \rightarrow Y$ and $g: Y \rightarrow Z$
are a composable pair of morphisms in $\calC$. If any two of $f$,
$g$, and $g \circ f$ belong to $S$, then so does the third.
\item If $f: X \rightarrow Y$ and $g: X \rightarrow Z$ are
morphisms of $S$ and $f$ belongs to $S$, then the induced morphism
$f': Z \rightarrow Z \coproduct_{X} Y$ belongs to $S$.
\item The full subcategory of $\calC^{\rightarrow}$ consisting of
morphisms which lie in $S$ is stable under the formation of
colimits.
\end{itemize}
\end{definition}

\begin{remark}
Let $S$ be a collection of morphisms in $\calC$. It is clear that
every morphism in $S$ is an $S$-equivalence. If $\calC$ is
presentable, then the collection of $S$-equivalences is saturated.
\end{remark}

\begin{remark}
If $S$ is saturated, and $f: X \rightarrow Y$ lies in $S$, then
the induced map $X \otimes S^n \rightarrow Y \otimes S^n$ lies in
$S$.
\end{remark}

\begin{remark}
If $\calC$ is any presentable $\infty$-category, then the
collection of all equivalences in $\calC$ is saturated.
\end{remark}

\begin{remark}
If $F: \calC \rightarrow \calC'$ is a colimit preserving functor
and $S'$ is a saturated collection of morphisms of $\calC'$, then
$S = \{ f \in \Hom_{\calC}(X,Y): Ff \in S'\}$ is saturated. In
particular, the collection of morphisms $f$ such that $Ff$ is an
equivalence is saturated.
\end{remark}

\begin{remark}
Let $\calC$ be a presentable $\infty$-category. Then any
intersection of saturated collections of morphisms in $\calC$ is
saturated. Consequently, for any collection $S$ of morphisms of
$\calC$ there is a smallest saturated collection $\overline{S}$
containing $S$, which we call the {\it saturation} of $S$. It
follows immediately that every morphism in $\overline{S}$ is an
$S$-equivalence.

We shall say that a saturated collection $S$ of morphisms is {\it
setwise generated} if $S = \overline{S_0}$ for some {\em set} $S_0
\subseteq S$.
\end{remark}

\begin{proposition}\label{local}
Let $\calC$ be a presentable $\infty$-category, and let $S$ be a
set of morphisms in $\calC$. Let $\calC_{S}$ denote the full
subcategory of $\calC$ consisting of $S$-local objects. Then the
inclusion $\calC_{S} \rightarrow \calC$ has an accessible left
adjoint $L$. Furthermore, the following are equivalent for a
morphism $f: X \rightarrow Y$ in $\calC$:
\begin{enumerate}
\item The morphism $f$ lies in the saturation $\overline{S}$.
\item The morphism $f$ is an $S$-equivalence.
\item The morphism $Lf$ is an equivalence.
\end{enumerate}
\end{proposition}

\begin{proof}
We may enlarge $S$ by adjoining elements of $\overline{S}$ without
changing the collection of $S$-local objects. Thus, we may assume
that for any morphism $s: X \rightarrow Y$ in $S$, the morphism $X
\coproduct_{X \otimes S^n} (Y \otimes S^n) \rightarrow Y$ also
belongs to $S$. This enlargement allows us to simplify the
definition of an $S$-local object: by Whitehead's theorem, an
object $C$ is $S$-local if and only if
$$\Hom_{\calC}(Y,C) \rightarrow \Hom_{\calC}(X,C)$$
is surjective on $\pi_0$ for every morphism $s: X \rightarrow Y$
in $S$.

Choose a regular cardinal $\kappa$ such that $\calC$ is
$\kappa$-accessible, and such that the source and target of every
morphism in $S$ is $\kappa$-compact. We will show that for every
object $C \in \calC_{\kappa}$, there exists a morphism $f: C
\rightarrow C_S$ such that $f \in \overline{S}$ and $C_S$ is
$S$-local. This will prove the existence of the left adjoint $L$
(carrying $C$ to $C_S$), at least on the subcategory
$\calC_{\kappa} \subseteq \calC$. This functor $L$ has a unique
$\kappa$-continuous extension to $\calC$. Since any
$\kappa$-filtered colimit of $S$-local objects is $S$-local, this
$\kappa$-continuous extension will be a left adjoint defined on
all of $\calC$. The closure of $\overline{S}$ under
$\kappa$-filtered colimits will also show that the adjunction
morphism $C \rightarrow LC$ lies in $\overline{S}$, for any $C \in
\calC$.

Let us therefore assume that $C \in \calC_{\kappa}$. Let $\calJ$
denote the $\infty$-category whose objects are morphisms $C
\rightarrow C'$ in $\calC_{\kappa}$, which lie in $\overline{S}$.
One readily verifies (using the fact that $\overline{S}$ is
saturated) that $\calJ$ is a $\kappa$-filtered category, equipped
with a functor $\pi: \calJ \rightarrow \calC_{\kappa}$ carrying
the diagram $C \rightarrow C'$ to the object $C' \in
\calC_{\kappa}$. We may therefore regard $\pi$ as a
$\kappa$-filtered diagram in $\calC_{\kappa}$, which has a colimit
$C_S$ in $\calC$. Since $\pi$ is a diagram of morphisms in
$\overline{S}$, the natural map $C \rightarrow C_S$ lies in
$\overline{S}$. We must now show that $C_S$ is $S$-local.

Suppose $s: X \rightarrow Y$ is a morphism in $S$. We need to show
that $\pi_0 \Hom_{\calC}(Y,C_S) \rightarrow \pi_0
\Hom_{\calC}(X,C_S)$ is surjective. To this end, choose any
morphism $f: X \rightarrow C_S$. Since $X$ is $\kappa$-compact,
the morphism $f$ factors through some morphism $X \rightarrow C'$,
where $C \rightarrow C'$ appears in the $\kappa$-filtered diagram
$\pi$. Then the pushout morphism $C \rightarrow C' \coprod_{X} Y$
lies in $\overline{S}$, and so also occurs in the diagram $\pi$.
Thus we obtain a map $f': Y \rightarrow (C' \coprod_{X} Y)
\rightarrow C_S$. It is clear from the construction that $f$
factors through $f'$.

It remains to verify the equivalence of the three conditions
listed above. We have already seen that $(1) \implies (2)$, and
the equivalence of $(2)$ and $(3)$ is formal. To complete the
proof, it suffices to show that $(3)$ implies $(1)$. If $f: X
\rightarrow Y$ induces an equivalence after applying $L$, then the
morphisms $X \rightarrow Y_S$ and $Y \rightarrow Y_S$ both lie in
$\overline{S}$, so that $f$ lies in $\overline{S}$ since
$\overline{S}$ is saturated.
\end{proof}

We note a universal property enjoyed by the localization
$\calC_S$:

\begin{proposition}
Let $\calC$ and $\calC'$ be presentable $\infty$-categories, let
$S$ be a set of morphisms of $\calC$, and let $L: \calC
\rightarrow \calC_S$ be left adjoint to the inclusion. Composition
with $L$ induces a fully faithful functor $\pi: \Pre(\calC_{S},
\calC') \rightarrow \Pre(\calC, \calC')$ whose essential image
consists of those functors $F: \calC \rightarrow \calC'$ for which
$Fs$ is an equivalence for all $s \in S$.
\end{proposition}

\begin{proof}
Restriction induces a functor $\psi: \Pre(\calC, \calC')
\rightarrow \Pre(\calC_{S}, \calC')$ which is right adjoint to
$\pi$. It is clear that the adjunction $\psi \circ \pi \rightarrow
1$ is an equivalence, which shows that $\pi$ is fully faithful. By
construction, the essential image consists of all functors for
which $F$ induces equivalences $FC \simeq FLC$ for all $C \in
\calC$. Since every morphism $C \rightarrow LC$ lies in
$\overline{S}$, it will suffice to show that if $F$ inverts every
morphism in $S$, then it inverts every morphism in $\overline{S}$.
This follows immediately from the observation that the set of
morphisms which become invertible after application of $F$ is
saturated.
\end{proof}

We conclude this section by showing that {\em every} localization
of an presentable $\infty$-category $\calC$ has the form
$\calC_{S}$, where $S$ is some set of morphisms of $\calC$.

We begin with the following observation: let $F: \calC \rightarrow
\calC'$ be an accessible functor between accessible
$\infty$-categories, and let $S$ denote the full subcategory of
$\calC^{\rightarrow}$ consisting of morphisms $s$ such that $Fs$
is an equivalence. Then $S$ may be constructed as a limit of
accessible $\infty$-categories, so that $S$ is accessible: in
particular, there exists a {\em set} of objects $S_0 \subseteq S$
such that every object in $S$ is a filtered colimit of objects in
$S_0$. It follows that each element of $S$ is an
$S_0$-equivalence, so that if $\calC$ is presentable then $S
\subseteq \overline{S_0}$. If in addition $F$ is a
colimit-preserving functor between presentable
$\infty$-categories, then $S$ is saturated so we obtain $S =
\overline{S_0}$. Consequently, we obtain the following:

\begin{theorem}
Let $\calC$ be a presentable $\infty$-category. Let $S$ be a
collection of morphisms in $\calC$. The following conditions are
equivalent:

\begin{enumerate}
\item The collection $S$ is the saturation of some {\em set} of
morphisms in $\calC$.
\item There exists a colimit preserving functor $F: \calC
\rightarrow \calC'$ such that $S$ consists of all morphisms $s$
for which $Fs$ is an equivalence.
\item There exists a localization functor $L: \calC \rightarrow
\calC$ such that $S$ consists of all morphisms $s$ for which $Ls$
is an equivalence.
\end{enumerate}

\end{theorem}

\begin{proof}
The above discussion shows that $(2)$ implies $(1)$, and it is
clear that $(3)$ implies $(2)$. Proposition \ref{local} shows that
$(1)$ implies $(3)$.
\end{proof}

It follows that {\em any} localization functor defined on a
presentable $\infty$-category $\calC$ may be obtained as a
localization $\calC \rightarrow \calC_{S}$ for some {\em set} of
morphisms $S$ of $\calC$.

\section{$\infty$-Topoi}\label{toposes}

In this section, we will define the notion of an $\infty$-topos.
These were defined by Simpson in \cite{simpson} as
$\infty$-categories which occur as left-exact localizations of
$\infty$-categories of prestacks. In \S \ref{toposdef}, we give
our own definition, which we will show to be equivalent to
Simpson's definition in \S \ref{leftexactloc}. We will then
explain how to extract topoi from $\infty$-topoi and vice versa.
Finally, in \S \ref{hyperstack}, we explain how to obtain the
Joyal-Jardine $\infty$-topos of hyperstacks on topos $X$ as a
localization of our $\infty$-topos of stacks on $X$.

There are several papers on higher topoi in the literature. The
papers \cite{street} and \cite{ditopoi} both discuss a notion of
$2$-topos (the second from an elementary point of view). However,
the basic model for these $2$-topoi was the $2$-category of
(small) categories, rather than the $2$-category of (small)
groupoids. Jardine (\cite{jardine}) has exhibited a model
structure on the simplicial presheaves on a Grothendieck site, and
the $\infty$-category associated to this model category is an
$\infty$-topos in our sense. This construction was generalized
from ordinary categories with a Grothendieck topology to
simplicial categories with a Grothendieck topology (appropriately
defined) in \cite{toen}; this also produces $\infty$-topoi.
However, not every $\infty$-topos arises in this way: one can
construct only $\infty$-topoi which are $t$-complete in the sense
of \cite{toen}; we will summarize the situation in Section
\ref{hyperstack}. We should remark here that our notion of
$\infty$-topos is equivalent to the notion of a {\it Segal topos}
defined in \cite{toen}.

In our terminology, ``topos'' shall always mean ``Grothendieck
topos'': we will discuss an ``elementary'' version of these ideas
in a sequel to this paper.

\subsection{Equivalence Relations and Groupoid Objects}

There are several equivalent ways to define a (Grothendieck)
topos. The following is proved in \cite{SGA}:

\begin{proposition}\label{toposdefined}
Let $\calC$ be a category. The following are equivalent:
\begin{enumerate}
\item The category $\calC$ is $($equivalent to$)$ the category of
sheaves on some Grothendieck site.
\item The category $\calC$ is $($equivalent to$)$ a left-exact
localization of the category of presheaves on some small category.
\item The category $\calC$ satisfies Giraud's
axioms:
\begin{itemize}
\item $\calC$ has a set of small generators.
\item $\calC$ admits arbitrary colimits.
\item Colimits in $\calC$ are universal (that is, the formation of
colimits commutes with base change).
\item Sums in $\calC$ are disjoint.
\item All equivalence relations in $\calC$ are effective.
\end{itemize}
\end{enumerate}
\end{proposition}

We will prove an analogue of this result (Theorem \ref{giraud}) in
the $\infty$-categorical context. However, our result is quite so
satisfying: it establishes the equivalence of analogues of
conditions $(2)$ and $(3)$, but does not give an explicit method
for constructing all of the left-exact localizations analogous to
$(1)$.

In the next section, we will define an $\infty$-topos to be an
$\infty$-category which satisfies $\infty$-categorical versions of
Giraud's axioms. The first four of these axioms generalize in a
very straightforward way to $\infty$-categories, but the last (and
most interesting axiom) is more subtle.

Recall that if $X$ is an object in an (ordinary) category $\calC$,
then an {\it equivalence relation} $R$ on $X$ is an object of
$\calC$ equipped with a map $p: R \rightarrow X \times X$ such
that for any $S$, the induced map $\Hom_{\calC}(S,R) \rightarrow
\Hom_{\calC}(S,X) \times \Hom_{\calC}(S,X)$ exhibits
$\Hom_{\calC}(S,R)$ as an equivalence relation on
$\Hom_{\calC}(S,X)$.

If $\calC$ admits finite limits, then it is easy to construct
equivalence relations in $\calC$: given any map $X \rightarrow Y$
in $\calC$, the induced map $X \times_Y X \rightarrow X \times X$
is an equivalence relation on $X$. If the category $\calC$ admits
finite colimits, then one can attempt to invert this process:
given an equivalence relation $R$ on $X$, one can form the
coequalizer of the two projections $R \rightarrow X$ to obtain an
object which we will denote by $X/R$. In the category of sets, one
can recover $R$ as the fiber product $X \times_{X/R} X$. In
general, this need not occur: one always has $R \subseteq X
\times_{X/R} X$, but the inclusion may be strict (as subobjects of
$X \times X$). If equality holds, then $R$ is said to be an {\it
effective equivalence relation}, and the map $X \rightarrow X/R$
is said to be an {\it effective epimorphism}. In this terminology,
we have the following:

\begin{fact}\label{factoid}
In the category $\calC$ of sets, every equivalence relation is
effective and the effective epimorphisms are precisely the
surjective maps.
\end{fact}

The first assertion of Fact \ref{factoid} remains valid in any
topos, and according to the axiomatic point of view it is one of
the defining features of topoi.

If $\calC$ is a category with finite limits and colimits in which
all equivalence relations are effective, then we obtain a
one-to-one correspondence between equivalence relations on an
object $X$ and ``quotients'' of $X$ (that is, isomorphism classes
of effective epimorphisms $X \rightarrow Y$). This correspondence
is extremely useful because it allows us to make elementary
descent arguments: one can deduce statements about quotients of
$X$ from statements about $X$ and about equivalence relations on
$X$ (which live over $X$). We would like to have a similar
correspondence in certain $\infty$-categories.

In order to formulate the right notions, let us begin by
considering the $\infty$-category $\SSet$ of spaces. The correct
notion of surjection of spaces $X \rightarrow Y$ is a map which
induces a surjection on components $\pi_0 X \rightarrow \pi_0 Y$.
However, in this case, the fiber product $R= X \times_Y X$ is does
not give an equivalence relation on $X$, because the map $R
\rightarrow X \times X$ is not necessarily ``injective'' in any
reasonable sense. However, it does retain some of the pleasant
features of an equivalence relation: instead of transitivity, we
have a {\it coherently associative} composition map $R \times_X R
\rightarrow R$. In order to formalize this observation, it is
convenient to introduce some simplicial terminology.

Let $\Delta$ denote the (ordinary) category of finite, nonempty,
linearly ordered sets. If $\calC$ is an $\infty$-category, then a
{\it simplicial object} of $\calC$ is a functor $C_{\bigdot}:
\Delta^{op} \rightarrow \calC$. Since every object of $\Delta$ is
isomorphic to $\{0, \ldots, n\}$ for some $n \geq 0$, we may think
of $C_{\bigdot}$ as being given by objects $\{ C_{n} \}_{n \geq
0}$ in $\calC$, together with various maps between the $C_{n}$
which are compatible up to coherent homotopy. In particular, we
have two maps $\pi_0, \pi_1: C_1 \rightarrow C_0$, and we have
natural maps $p_n: C_{n} \rightarrow C_1 \times_{C_0} \ldots
\times_{C_0} C_1$, where there are $n$ factors in the product and
for each fiber product, the map from the left copy of $C_1$ to
$C_0$ is given by $\pi_1$ and the map from the right copy of $C_1$
to $C_0$ is given by $\pi_0$. If $p_n$ is an equivalence for each
$n \geq 0$, then we shall say that $C_{n}$ is a {\it category
object of $\calC$}. (Note that this definition does not really
require that $\calC$ has fiber products: a priori, the target of
$p_n$ may be viewed as a prestack on $\calC$.)

\begin{example}
Suppose that $\calC$ is the category of sets. Then a category
object $C_{\bigdot}$ of $\calC$ consists of a set $C_0$ of {\it
objects}, a set $C_1$ of {\it morphisms}, together with source and
target maps $\pi_0, \pi_1: C_1 \rightarrow C_0$ and various other
maps between iterated fiber powers of $C_1$ over $C_0$. In
particular, one has a map $C_1 \times_{C_0} C_1 \rightarrow C_1$
which gives rise to an multiplication on ``composable'' pairs of
morphisms. One can show that this recovers the usual definition of
a category.
\end{example}

\begin{example}
More generally, suppose that $\calC$ is an ordinary category. Then
what we have called a category object of $\calC$ is equivalent to
the usual notion of a category object of $\calC$: namely, data
which represents a ``category-valued functor'' on $\calC$.
\end{example}

In the general case, the requirement that $p_n$ be an equivalence
shows that $C_n$ is determined by $C_1$ and $C_0$. Thus, a
category object in $\calC$ may be viewed as a pair of morphisms
$\pi_0,\pi_1: C_1 \rightarrow C_0$, together with some additional
data which formalizes the idea that there should be a coherently
associative composition law on composable elements of $C_1$.

If an $\infty$-category $\calC$ has fiber products, then any
morphism $U \rightarrow X$ gives rise to a category object
$U_{\bigdot}$ of $\calC$, with $U_{n}$ given by the $(n+1)$-fold
fiber power of $U$ over $X$. However, this category object has a
further important property: it is a {\it groupoid object}. This
can be formalized in many ways: for example, one can assert the
existence of an ``inverse map'' $i: U_1 \rightarrow U_1$ such that
$\pi_0 \circ i \simeq \pi_1$, $\pi_1 \circ i \simeq \pi_0$, and
various compositions involving the inverse map are homotopic to
the identity. It is important to note that the inverse map is
canonically determined by the category structure of $\{
U_{\bigdot} \}$ {\em if it exists}, and therefore we will not need
to axiomatize its properties. In other words, a groupoid object in
$\calC$ is a category object of $\calC$ which {\em satisfies
conditions}, rather than a category object with {\em extra
structure}.

If an $\infty$-category $\calC$ admits all limits and colimits,
then we may view any simplicial object $C_{\bigdot}: \Delta^{op}
\rightarrow \calC$ as a diagram in $\calC$ and form its colimit $|
C_{\bigdot} |$, which is also called the {\it geometric
realization} of $C_{\bigdot}$. There is a natural map $C_0
\rightarrow | C_{\bigdot} |$, which induces a map of simplicial
objects $\psi: C_{\bigdot} \rightarrow U_{\bigdot}$, where $U_n$
denotes the $(n+1)$-fold fiber power of $C_0$ over $| C_{\bigdot}
|$. If the map $\psi$ is an equivalence, then we shall say that
$C_{\bigdot}$ is an {\it effective groupoid}, and that $C_0
\rightarrow | C_{\bigdot} |$ is an {\it effective epimorphism}.

\begin{remark}
Since $C_0 \simeq U_0$ and the objects $C_n$ ($U_n$) are
determined by taking fiber products of $C_1$ over $C_0$ ($U_1$
over $U_0$), we see that $\psi$ is an equivalence if and only if
it induces an equivalence $C_1 \simeq U_1$. In other words,
$C_{\bigdot}$ is effective if and only if $C_1 \simeq C_0
\times_{| C_{\bigdot} |} C_0$.
\end{remark}

The notion of a groupoid object in $\calC$ will be our
$\infty$-categorical replacement for the $1$-categorical notion of
an object with an equivalence relation. This new notion is useful
thanks to the following analogue of Fact \ref{factoid}:

\begin{proposition}\label{steep}
In the $\infty$-category $\SSet$ of spaces, all groupoid objects
are effective and a map $X \rightarrow Y$ is an effective
epimorphism if and only if $\pi_0 X \rightarrow \pi_0 Y$ is
surjective.
\end{proposition}

In contrast to Fact \ref{factoid}, Proposition \ref{steep} is
nontrivial. For example, a category object $U_{\bigdot}$ in
$\SSet$ with $U_{0} = \ast$ is more-or-less the same thing as an
$A_{\infty}$-space in the terminology of \cite{stasheff} (in other
words, the space $U_1$ has a coherently associative multiplication
operation). The first part of Proposition \ref{steep} asserts, in
this case, that an $A_{\infty}$-space whose set of connected
components forms a group can be realized as a loop space. We refer
the reader to our earlier discussion in Remark \ref{content}.

\begin{remark}
In some sense, the notion of a groupoid object is {\em simpler
than} the notion of an equivalence relation. For example, let
$\calC$ be the category of sets. Then groupoid objects are simply
groupoids in the usual sense, while equivalence relations
correspond to groupoids satisfying an extra discreteness condition
(no objects can have nontrivial automorphisms). The removal of
this discreteness condition is what distinguishes the theory of
$\infty$-topoi from the theory of ordinary topoi. It permits us to
form quotients of objects using more general kinds of ``gluing
data''. In geometric contexts, this extra flexibility allows the
construction of useful objects such as orbifolds and algebraic
stacks, which are useful in a variety of mathematical situations.

One can imagine weakening the gluing conditions even further, and
considering axioms having the form ``every category object is
effective''. This seems to be a very natural approach to a theory
of topos-like $(\infty,\infty)$-categories. We will expound
further on this idea in a sequel to this paper.
\end{remark}

\subsection{The Definition of an $\infty$-Topos}\label{toposdef}

Just as a topos is a category which resembles the category of
sets, we shall define an $\infty$-topos to be an $\infty$-category
which resembles the $\infty$-category $\SSet$ of spaces. We will
begin by extracting out certain properties enjoyed by $\SSet$. We
will then take these properties as the axioms for an
$\infty$-topos. They are the analogues of Giraud's conditions
which characterize topoi in the case of ordinary categories.

\begin{fact}
The $\infty$-category $\SSet$ is presentable.
\end{fact}

In particular, $\SSet$ admits all limits and colimits, so that
fiber products exist in $\SSet$. However, fiber products have
special properties in the $\infty$-category $\SSet$:

\begin{fact}\label{co-fib}
In $\SSet$, the formation of colimits commutes with pullback.
\end{fact}

It follows from Fact \ref{co-fib} that the initial object
$\emptyset$ of $\SSet$ really behaves as if it were ``empty''.
Since $\emptyset$ is the colimit of the empty diagram, we see that
for any map $X \rightarrow \emptyset$, $X$ is the colimit of the
pullback of the empty diagram (which is also empty). Thus $X
\simeq \emptyset$.

\begin{fact}
Sums in $\SSet$ are disjoint. That is, given any two spaces $X$
and $Y$, the fiber product $X \times_{X \coprod Y} Y$ is an
initial object.
\end{fact}

Finally, we recall Proposition \ref{steep}, which implies in
particular:

\begin{fact}\label{triv}
Any groupoid object in $\SSet$ is effective.
\end{fact}

We are now in a position to define the notion of an
$\infty$-topos.

\begin{definition}\label{toposs}
An $\infty$-topos is an $\infty$-category $\calX$ with the
following properties:

\begin{itemize}
\item $\calX$ is presentable.

\item The formation of colimits in $\calX$ commutes with pullback.

\item Sums in $\calX$ are disjoint.

\item Every $\calX$-groupoid is effective.
\end{itemize}

\end{definition}

Then it follows from the discussion up to this point that:

\begin{proposition}
The $\infty$-category $\SSet$ is an $\infty$-topos.
\end{proposition}

We also note that any slice of an $\infty$-topos is an
$\infty$-topos.

\begin{proposition}\label{tivv}
Let $\calX$ be an $\infty$-topos, and let $E \in \calX$ be an
object. Then the slice $\infty$-category $\calX_{/E}$ is an
$\infty$-topos.
\end{proposition}

\begin{proof}
Proposition \ref{slicepresentable} implies that $\calX_{/E}$ is
presentable. Since colimits and fiber products in $\calX_{/E}$
agree with colimits and fiber products in $\calX$, the other
axioms follow immediately.
\end{proof}

\begin{remark}\label{mark}
Let $\calX$ be a presentable $\infty$-category in which pullbacks
commute with colimits. For any morphism $p: X \rightarrow E$ in
$\calX$, we define a presheaf $\calF$ on $\calX$ in the following
way: $\calF(Y) = \Hom_{E}( Y \times E, X)$. Since forming the
product with $E$ preserves colimits, we deduce from Proposition
\ref{representable} that $\calF$ is representable. In other words,
there exists an object of $\calX$ which represents the functor
``sections of $p$''.
\end{remark}

\subsection{Geometric Morphisms}

The proper notion of morphism between $\infty$-topoi is that of a
{\it geometric morphism}, which we now introduce. Recall that a
functor $F: \calC \rightarrow \calC'$ between $\infty$-categories
is said to be {\it left exact} if $F$ preserves finite limits.

\begin{definition}
Let $\calX$ and $\calY$ be $\infty$-topoi. A {\it geometric
morphism} $f$ from $\calX$ to $\calY$ is a left exact functor
$f^{\ast}: \calY \rightarrow \calX$ which preserves all colimits.
\end{definition}

It follows from Corollary \ref{adjointfunctor} that $f^{\ast}$ has
a right adjoint, which we shall denote by $f_{\ast}$.

It is clear that the class of geometric morphisms is stable under
composition. We may summarize the situation by saying that there
is an $(\infty,2)$-category $\Top^{\infty}$ of $\infty$-topoi and
geometric morphisms. If $\calX$ and $\calY$ are $\infty$-topoi, we
will write $\Top^{\infty}(\calX, \calY)$ for the $\infty$-category
of geometric morphisms from $\calX$ to $\calY$.

\begin{example}\label{cant}
Let $\calX$ be an $\infty$-topos, and let $E \in \calX$ be an
object. Then there is a geometric morphism $f_E: \calX_{/E}
\rightarrow \calX$, defined by $f_E^{\ast}(E') = E' \times E$. The
geometric morphisms which arise in this way are important enough
to deserve a name: a geometric morphism $f: \calY \rightarrow
\calX$ is said to be {\it \etale}  if it factors as a composite
$f_E \circ f'$, where $f': \calY \rightarrow \calX_{/E}$ is an
equivalence.

Consider an \etale geometric morphism $f: \calX_{/E} \rightarrow
\calX$. Then the forgetful inclusion of $\calX_{/E}$ into $\calX$
is a {\em left} adjoint to $f^{\ast}$, which we shall denote by
$f_{!}$. It follows that $f^{\ast}$ commutes with {\em all}
limits. The right adjoint of $f^{\ast}$ is slightly more difficult
to describe: it is given by $f_{\ast} X = Y$, where $Y$ is the
object of $X$ representing sections of $X \rightarrow E$ (see
Remark \ref{mark}).
\end{example}

\begin{proposition}\label{cantt}
The $(\infty,2)$-category $\Top^{\infty}$ admits arbitrary
$(\infty,2)$-categorical colimits. These are formed by taking the
limits of the underlying $\infty$-categories.
\end{proposition}

\begin{proof}
We have seen that any limit of $\infty$-pretopoi is an
$\infty$-pretopos. In the case where the $\infty$-categories
involved are $\infty$-topoi and the functors between them are
left-exact, one shows that the limit is actually an $\infty$-topos
by checking the axioms directly.
\end{proof}

\begin{corollary}
Let $\calC_0$ be a small $\infty$-category. Then the category of
presheaves $\calP = \SSet^{\calC_0^{op}}$ is an $\infty$-topos.
\end{corollary}

\begin{proof}
The $\infty$-category $\calP$ is a limit of copies of $\SSet$,
indexed by $\calC_0^{op}$. Hence it is a colimit of copies of
$\SSet$ in $\Top^{\infty}$.
\end{proof}

In fact, $\Top^{\infty}$ has all $(\infty,2)$-categorical limits
as well, but this is more difficult to prove and we will postpone
a discussion until the sequel to this paper.

To conclude this section, we shall show that the $\infty$-category
of geometric morphisms between two $\infty$-topoi is reasonably
well-behaved. First, we need to introduce a technical notion and a
lemma (which will also be needed later).

\begin{definition}
Let $\kappa$ be a regular uncountable cardinal and $\calX$ an
$\infty$-topos. We shall say that $\calX$ is {\it
$\kappa$-coherent} if the following conditions are satisfied:

\begin{itemize}
\item $\calX$ is $\kappa$-accessible.
\item $\calX_{\kappa}$ is stable under the formation of finite limits (in $\calX$).
\item The formation of finite limits in $\calX$ commutes with
$\kappa$-filtered colimits.
\end{itemize}
\end{definition}

\begin{remark}
If $\calX$ is $\kappa$-coherent, then $\calX$ is
$\kappa'$-coherent for any regular $\kappa' \geq \kappa$.
\end{remark}

\begin{remark}
This does not seem to be the right notion when $\kappa=\omega$. In
$\SSet$, the compact objects (retracts of finite complexes) are
not stable under finite limits. For example, the fiber product
$\ast \times_{S^1} \ast$ is the infinite discrete set $\pi_1(S^1)
= \Z$, which is not compact. Compare with Remark \ref{brak}.
\end{remark}

To conclude this section, we will show that the $\infty$-category
of geometric morphisms between two $\infty$-topoi is always an
accessible category. We begin with the following Lemma, which will
also be needed in the proof of Theorem \ref{giraud}.

\begin{lemma}\label{limitfiber}
Let $\calX$ be an $\infty$-topos. Then there exists a regular
cardinal $\kappa > \omega$ such that $\calX$ is $\kappa$-coherent.
\end{lemma}

\begin{proof}
Suppose that there exists a regular cardinal $\kappa_0$ such that
the formation of finite limits in $\calX$ commutes with
$\kappa_0$-filtered colimits. Enlarging $\calX$ if necessary, we
may assume that $\calX$ is $\kappa_0$-accessible. It then suffices
to choose any uncountable $\kappa \geq \kappa_0$ such that all
finite limits of diagrams in $\calX_{\kappa_0}$ lie in
$\calX_{\kappa}$.

It remains to prove the existence of $\kappa_0$. If $\calX =
\SSet$, then we may take $\kappa_0 = \omega$. Working
componentwise, we see that the $\kappa_0 = \omega$ works whenever
$\calX$ is an $\infty$-category of presheaves.

By Theorem \ref{giraud}, we may assume that $\calX$ is the
essential image of some left exact localization functor $L: \calP
\rightarrow \calP$, where $\calP$ is an $\infty$-category of
presheaves. Then $L$ commutes with fiber products, and with
$\kappa_0$-filtered colimits for $\kappa_0$ sufficiently large.
Thus, the assertion for $\calX$ follows from the assertion for
$\calP$.
\end{proof}

The next result will not be used in the rest of the paper, and may
be safely omitted by the reader. The proof involves routine, but
tedious, cardinal estimates.

\begin{proposition}
Let $\calX$ and $\calY$ be $\infty$-topoi. Then the
$\infty$-category $\Top(\calX, \calY)$ geometric morphisms from
$\calX$ to $\calY$ is accessible.
\end{proposition}

\begin{proof}
The $\infty$-category $\Top(\calX,\calY)$ is a full subcategory of
$\PTop(\calY, \calX)$, which is an $\infty$-pretopos in which
colimits may be computed pointwise. Let $\kappa$ be a regular
cardinal which satisfies the conclusion of Lemma \ref{limitfiber}
for both $\calX$ and $\calY$. Without loss of generality, we may
enlarge $\kappa$ so that $\calX$, $\calY$, and $\PTop(\calY,
\calX)$ are $\kappa$-accessible.

We next note that the condition that some functor $F \in
\PTop(\calY, \calX)$ be left-exact can be checked by examining a
bounded amount of data. Namely, we claim that $F$ is left exact if
and only if for any objects $Y' \rightarrow Y \leftarrow Y''$ in
$\calY_{\kappa}$ and any $X \in \calX_{\kappa}$, the natural map
$\pi^F_{X,Y}: \Hom_{\calX}(X, F(Y' \times_{Y} Y'')) \rightarrow
\Hom_{\calX}(X,FY') \times_{ \Hom_{\calX}(X,FY) }
\Hom_{\calX}(X,FY'')$ is an equivalence. This condition is
obviously necessary. Assume that the condition holds. {\em Any}
diagram $Y' \rightarrow Y \leftarrow Y''$ in $\calY$ may be
obtained as a $\kappa$-filtered colimit of diagrams $Y'_{\alpha}
\rightarrow Y_{\alpha} \leftarrow Y''_{\alpha}$ in
$\calY_{\kappa}$. Since $\calY$ is $\kappa$-coherent, we deduce
that $Y' \times_Y Y''$ is the colimit of the diagram $\{
Y'_{\alpha} \times_{Y_{\alpha}}  Y''_{\alpha} \}$. Using the
continuity of $F$, we see that $\pi^F_{X,Y}$ is an equivalence for
every $\kappa$-compact $X$ (with no compactness assumptions on the
triple $(Y,Y',Y'')$. Since every object of $\calX$ can be written
as a colimit of $\kappa$-compact objects, we deduce that
$\pi^F_{X,Y}$ is an equivalence in general.

Let $\kappa'$ be a regular cardinal $\geq \kappa$, such that
$\calX_{\kappa}$ and $\calY_{\kappa}$ have cardinality $<
\kappa'$, and the space $\Hom_{\calX}(A,B)$ is $\kappa'$-presented
for all $A,B \in \calX_{\kappa}$. To complete the proof, we will
show that if $F \in \PTop(\calY, \calX)$ is left exact, then $F$
can be written as a $\kappa'$-filtered colimit of left-exact
objects in $\PTop(\calY, \calX)_{\kappa'}$. Since $\PTop(\calY,
\calX)$ is $\kappa'$-accessible, we may write $F$ as a
$\kappa'$-filtered colimit of some diagram $\pi: \calI \rightarrow
\PTop(\calY, \calX)_{\kappa'}$. To simplify the discussion, let us
assume that $\calI$ is a $\kappa'$-directed partially ordered set.

Let $P$ denote the set of all directed, $\kappa'$-small subsets of
$\calI$. For each $S \in P$, let $F_S$ denote the colimit of the
system $\{ \pi(s) \}_{s \in S}$. Then $F$ may also be written as a
colimit of the functors $\{ F_S \}_{S \in P}$. To complete the
proof, it will suffice that $P_0 = \{ S \in P: F_S \in \Top(\calX,
\calY) \}$ is cofinal in $P$.

Choose $S \in P$. We must show that there exists $S' \in P_0$
which contains $S$. We will construct a transfinite sequence of
subsets $S_{\alpha} \in P$, with $S_0 = S$ and $S_{\lambda} =
\bigcup_{\alpha < \lambda} S_{\alpha}$ for limit ordinals
$\lambda$. For successor ordinals, our construction must be a
little more complicated. Suppose that $S_{\alpha}$ has been
defined, and let $F_{\alpha} = F_{S_{\alpha}}$. Consider all
instances of the following data:

\begin{itemize}
\item $X \in \calX_{\kappa}$
\item $Y' \rightarrow Y \leftarrow Y'' \in \calY_{\kappa}$
\item $p$ is a point of $\Hom_{\calX}(X, F_{\alpha}Y'
\times_{F_{\alpha}Y} F_{\alpha}Y'')$
\item $\phi$ is a map from $S^n$ to the fiber of
$\pi^{F_{\alpha}}_{X,Y}$, for some $n \geq 0$.
\end{itemize}

Let $p'$ denote the image of $p$ in $\Hom_{\calX}(X, FY'
\times_{FY} FY'')$. Since $F$ is left-exact, we deduce that the
induced map $\phi'$ from $S^n$ to the fiber of $\pi^F_{X,Y}$ over
$p'$ is homotopic to a constant. Since the formation of this fiber
is $\kappa'$-continuous, the same holds if we replace $F$ by
$F_{S_{\alpha} \cup T(x)}$, where $T(x) \subseteq \calI$ is
$\kappa'$-small, and $x$ denotes the complicated data described
above.

The hypotheses guarantee that there are fewer than $\kappa'$
possibilities for the data $x$. Consequently, we may choose a
filtered subset $S_{\alpha+1} \subseteq \calI$ containing
$S_{\alpha}$ and each $T(x)$. This completes the construction. We
note that $S_{\alpha}$ is $\kappa'$-small for each $\alpha <
\kappa'$. Set $S' = S_{\kappa}$.

To complete the proof, it suffices to show that $F_{\kappa}$ is
left exact. In other words, we must show that for any $X \in
\calX_{\kappa}$, $Y' \rightarrow Y \leftarrow Y''$ in
$\calY_{\kappa}$, the map $\pi^{F_{\kappa}}_{X,Y}$ is an
equivalence. By Whitehead's theorem, it suffices to show that for
any data $x$ of the sort described above for $\alpha = \kappa$,
$\phi$ is already homotopic to a constant map. This follows from
the construction, once we note that the formation of all relevant
objects commutes with $\kappa$-filtered colimits.
\end{proof}

\subsection{Giraud's Theorem}\label{leftexactloc}

We now come to the first non-trivial result of this paper, which
asserts that any $\infty$-topos can be obtained as a left-exact
localization of an $\infty$-category of prestacks. This should be
considered a version of Giraud's theorem for ordinary topoi, which
shows that any abstract category satisfying a set of axioms
actually arises as a category of sheaves on some site. Our result
is not quite as specific because we do not have an explicit
understanding of the localization functor as a ``sheafification
process".

\begin{theorem}\label{giraud}
Let $\calX$ be an $\infty$-category. The following are equivalent:
\begin{itemize}
\item The $\infty$-category $\calX$ is an $\infty$-topos.
\item There exists a cardinal $\kappa$ such that $\calX$ is
$\kappa$-accessible and the Yoneda embedding $\calX \rightarrow
\SSet^{\calX_{\kappa}^{op}}$ has a left exact left adjoint.
\item There exists an $\infty$-topos $\calP$ such that $\calX$ is
a left exact localization of $\calP$.
\end{itemize}
\end{theorem}

In order to prove Theorem \ref{giraud}, we will need a rather
technical lemma regarding the structure of ``free groupoids'',
which we will now formulate. Let $\calC$ be an $\infty$-category.
If $U$ is an object of $\calC$, then we shall call a pair of
morphisms $\pi_0, \pi_1: R \rightarrow U$ a {\it coequalizer
diagram over $U$}. For fixed $U$, the coequalizer diagrams over
$Y$ form an $\infty$-category $\coeq_U$. A {\it groupoid over $U$}
is a $\calC$-groupoid $U_{\bigdot}$ together with an
identification $U_0 \simeq U$. These groupoids also form an
$\infty$-category which we shall denote by $\group_U$.

Any groupoid $U_{\bigdot}$ over $U$ determines a coequalizer
diagram over $U$, using the two natural maps $\pi_0, \pi_1: U_1
\rightarrow U$. This determines a functor $G^{\calC}: \group_U
\rightarrow \coeq_U$. We are interested in constructing a
left-adjoint to $G^{\calC}$. In other words, we would like to
construct the ``free groupoid over $Y$ generated by a coequalizer
diagram over $Y$''. Moreover, we will need to know that the
construction of this free groupoid is compatible with various
functors.

\begin{lemma}\label{technical}
Let $\calX$ be an $\infty$-topos, and let $U$ be an object of
$\calX$. Then the following assertions hold:

\begin{itemize}
\item The functor $G^{\calX}: \group_U \rightarrow \coeq_U$ has a
left adjoint $F^{\calX}: \coeq_U \rightarrow \group_U$.

\item If $T: \calX \rightarrow \calX'$ is a functor between
$\infty$-topoi, then the identification $G^{\calC'} \circ T \simeq
T \circ G^{\calC}$ induces a natural transformation $F^{\calC'}
\circ T \rightarrow T \circ F^{\calC}$. This natural
transformation is an equivalence if $T$ preserves colimits and
fiber products over $U$.
\end{itemize}
\end{lemma}

\begin{proof}
We sketch the proof, which is based on Proposition \ref{free} from
the appendix, and other ideas introduced there. We begin by noting
that the $\infty$-category $\coeq_U$ has a monoidal structure. The
identity object is given by the diagram $1 = \{ \id_U, \id_U: U
\rightarrow U$. Given two coequalizer diagrams $A = \{ \pi_0,
\pi_1: R \rightarrow U \}$ and $A' = \{ \pi'_0, \pi'_1: R'
\rightarrow U \}$, we may define $A \otimes A' = \{ \pi_0, \pi'_1:
R \otimes_U R' \rightarrow U \}$, where the fiber product is
formed using $\pi_1$ and $\pi'_0$.

Since $\calX$ is an $\infty$-topos, we deduce that $\coeq_U$ is
presentable and that $\otimes$ is colimit-preserving (since
colimits commute with pullback in $\calX$). We may therefore apply
Proposition \ref{free} to produce monoid objects in $\coeq_U$. We
note that a monoid object in $\coeq_U$ is simply a {\it category
object} $U_{\bigdot}$ of $\calX$, together with an identification
$U \simeq U_0$.

To prove the first part of the lemma, we need to show that for any
object $N \in \coeq_U$, the ``free groupoid generated by $N$''
exists. To see this, we begin with the object $1 \in \coeq_U$ and
note that it comes equipped with a natural monoid structure (this
corresponds to the constant simplicial object $U_{\bigdot}$ having
value $U$). Applying Proposition \ref{free} to the morphism $1
\rightarrow 1 \coprod N$ in $\coeq_U$, we deduce the existence of
a monoid object $M^0 \in \coeq_U$. We may identify this monoid
object with a simplicial object $M^0_{\bigdot}$ in $\calX$ which
is a category object equipped with an identification $M^0_0 \simeq
U$. However, we are not yet done, because $M^0_{\bigdot}$ is not
necessarily a groupoid.

Our next goal is to promote $M^0$ to a groupoid. To do this, we
will define a transfinite sequence of monoid objects $M^{\alpha}
\in \coeq_U$, equipped with a coherent system of maps $M^{\beta}
\rightarrow M^{\alpha}$ for $\beta < \alpha$. We have already
defined $M^{\alpha}$ for $\alpha=0$, and when $\alpha$ is a limit
ordinal we will simply define $M^{\alpha}$ to be the appropriate
colimit. We are thereby reduced to handling the successor stages.

Assume that $M^{\alpha}$ has been defined. We will construct
$M^{\alpha+1}$ by freely adjoining inverses for all of the
morphisms in $M^{\alpha}$. This can be achieved using four
applications of Proposition \ref{free}. The idea is
straightforward but notationally difficult to describe, so we just
sketch the idea in a simplified example. Suppose that $\calC$ is
an ordinary category, and that $f: C \rightarrow C'$ is a morphism
in $\calC$ which we would like to formally adjoin an inverse for.
We can obtain this inverse in three steps:

\begin{itemize}
\item Adjoin a new morphism $g: C' \rightarrow C$.

\item Impose the relation $f \circ g = \id_{C'}$.

\item Impose the relation $g \circ f = \id_{C}$.
\end{itemize}

If $\calC$ is an $\infty$-category, then we need to be a bit more
careful, since we have now identified $g$ with the inverse of $f$
in two different ways. In this case, we need four steps:

\begin{itemize}
\item Adjoin a new morphism $g: C' \rightarrow C$.

\item Adjoin a path from $f \circ g$ to $\id_{C'}$.

\item Adjoin a path from $g \circ f$ to $\id_{C}$.

\item Adjoin a homotopy between the two paths that we now have
from $g \circ f \circ g$ to $g$.
\end{itemize}

The construction described above can be carried out in a relative
situation, using Proposition \ref{free}. The result is that we
obtain a new monoid $M^{\alpha+1} \in \coeq_U$. This monoid is
equipped with a monoid map $M^{\alpha} \rightarrow M^{\alpha+1}$.
If we think of $M^{\alpha+1}$ as a category object of $\calX$,
then it is obtained from $M^{\alpha}$ by forcing all of the
morphisms in $M^{\alpha}$ to become invertible.

By induction, one can easily check that for any groupoid object
$V_{\bigdot}$ of $\calC$, the map $\Hom( M^{\alpha}_{\bigdot},
V_{\bigdot} ) \rightarrow \Hom( M^{\beta}_{\bigdot}, V_{\bigdot})$
is an equivalence. Moreover, using the fact that $\calX$ is
presentable and some cardinality estimates, one can show that
$M^{\alpha}$ is itself a groupoid for sufficiently large $\alpha$
(it seems plausible that we can even take $\alpha=1$, but this is
not important). Thus for large $\alpha$, $M^{\alpha}$ is the
desired free groupoid.

The second assertion of the Lemma follows because the free
groupoid was constructed using colimits and iterated application
of Proposition \ref{free}, and the second assertion of Proposition
\ref{free} guarantees that the free constructions that we use are
compatible with colimit preserving, monoidal functors.
\end{proof}

If $\calX$ is an $\infty$-topos, there is an easier way to see
that $T$ has a left adjoint. Indeed, all groupoids in $\calC$ are
effective, so a groupoid $Y_{\bigdot}$ over $Y$ is uniquely
determined by the surjective map $Y \rightarrow |Y_{\bigdot}|$,
and conversely. It is then easy to see that the left adjoint to
$T$ should assign to any coequalizer diagram $X
\stackrel{\rightarrow}{\rightarrow} Y$ the groupoid $Y_{\bigdot}$
associated to the map $Y \rightarrow Z$, where $Z$ is the
coequalizer of $\pi_0, \pi_1: X \rightarrow Y$. In particular, we
have $Y_1 = Y \times_{Z} Y$. In particular, this proves the
following:

\begin{lemma}\label{tech2}
Let $\pi_0, \pi_1: X \rightarrow Y$ be a coequalizer diagram in an
$\infty$-topos $\calX$, with coequalizer $Z$. Suppose $F: \calX
\rightarrow \calX'$ is a functor between $\infty$-topoi which
commutes with all colimits, and with fiber products over $Y$. Then
the natural map
$$F(Y \times_{Z} Y) \rightarrow FY \times_{FZ} FY$$ is an
equivalence.
\end{lemma}

\begin{proof}
This follows immediately from the preceding discussion and Lemma
\ref{technical}.
\end{proof}

We now proceed to the main point.

\begin{lemma}
Let $\calX$ be an $\infty$-topos which is $\kappa$-coherent, let
$\calP = \SSet^{\calX_{\kappa}^{op}}$ be the $\infty$-category of
prestacks on $\calX_{\kappa}$, and let $G: \calX \rightarrow
\calP$ be the Yoneda embedding. Then the left adjoint $F$ of $G$
is left exact.
\end{lemma}

\begin{proof}
We first note that since $\calX$ is $\kappa$-accessible, the left
adjoint $F$ of $G$ exists by the proof of Proposition
\ref{simpsongiraud}. We must show that $F$ commutes with all
finite limits. It will suffice to show that $F$ commutes with
pullbacks and preserves the final object. The second point is
easy: since the final object $1 \in \calX$ is $\kappa$-compact, it
follows that the final object in $\calP$ is the presheaf
represented by $1 \in \calX_{\kappa}$, and $F$ carries this
presheaf into $1$.

Let us now show that $F$ commutes with pullbacks. In other words,
we must show that for any diagram $X \rightarrow Y \leftarrow Z$
in $\calP$, the natural morphism
$$\eta: F( X \times_Y Z) \rightarrow FX \times_{FY} FZ$$
is an equivalence. We proceed as follows: let $S$ denote the
collection of all objects $Y \in \calP$ such that $\eta$ is an
isomorphism for {\em any} pair of objects $X,Z$ over $Y$. We need
to show that every object of $\calP$ lies in $S$. Since every
prestack on $\calX_{\kappa}$ is a colimit of representable
prestacks, it will suffice to show that every representable
prestack lies in $S$ and that $S$ is closed under the formation of
colimits.

First of all, we note that since both $\calX$ and $\calP$ are
$\infty$-topoi, colimits commute with pullbacks. Since $F$
commutes with all colimits, we see that both $F(X \times_Y Z)$ and
$FX \times_{FY} FZ$ are compatible with arbitrary colimits in $X$
and $Z$. Since every prestack (over $Y$) is a colimit of
representable prestacks (over $Y$), in order to show that $Y \in
S$ it suffices to show that $\eta$ is an isomorphism whenever $X$
and $Z$ are representable. Suppose that $X$, $Y$, and $Z$ are all
representable by objects $x,y,z \in \calX_{\kappa}$, and let $w =
x \times_y z$ so that $w$ represents the prestack $W=X \times_Y
Z$. Then $L( X \times_Y Z) \simeq L(W) \simeq w \simeq x \times_y
z \simeq LX \times_{LY} LZ$, so that $Y \in S$. Therefore $S$
contains every representable prestack.

To complete the proof, we must show that $S$ is stable under the
formation of colimits. It will suffice to show that $S$ is stable
under the formation of sums and coequalizers (see Appendix
\ref{appendixdiagram}). We consider these two cases separately.

Suppose $Y \in \calP$ is a sum of some family of objects $\{
Y_{\alpha} \}_{\alpha \in A} \subseteq S$. As above, we may assume
that $X$ and $Z$ are represented by objects $x,z \in
\calX_{\kappa}$. The maps $X,Z \rightarrow Y$ correspond to points
$p_x \in Y_{\alpha}(x)$, $p_z \in Y_{\beta}(z)$ for some indices
$\alpha$ and $\beta$. Then $L( X \times_{Y} Z) = L( \emptyset) =
\emptyset$ if $\alpha \neq \beta$, and $$L( X \times_{Y} Z) = L(X
\times_{Y_{\alpha}} Z) = LX \times_{LY_{\alpha}} LZ$$ if
$\alpha=\beta$, where the last equality uses the fact that
$Y_{\alpha} \in S$. On the other hand, $$LX \times_{LY} LZ = x
\times_{ \coprod_{\gamma} LY_{\gamma} } z = x \times_{LY_{\alpha}}
(LY_{\alpha} \times_{ \coprod_{\gamma} LY_{\gamma}} LY_{\beta})
\times_{LY_{\beta}} z$$. Since sums are disjoint in $\calX$, this
fiber product is empty if $\alpha \neq \beta$ and is equal to $x
\times_{LY_{\alpha}} z$ otherwise.

We now come to the core of the argument, which is showing that $S$
is stable under the formation of coequalizers. Fortunately, the
hard work is already done. Suppose that $Y$ is the coequalizer of
a diagram $Y_1 \stackrel{\rightarrow}{\rightarrow} Y_0$ in
$\calP$, where $Y_0, Y_1 \in S$. As above, we may assume that $X$
and $Z$ are representable by objects $x,z \in \calC_{\kappa}$.
Then any point of $\Hom_{\calP}(X,Y) = Y(x)$ may be lifted to a
point of $Y_0(x)$. Thus we may assume that the maps $X, Z
\rightarrow Y$ both factor through $Y_0$. Since $X \times_Y Z= X
\times_{Y_0} (Y_0 \times_{Y} Y_0) \times_{Y_0} Z$ and we already
know that $L$ commutes with fiber products over $Y_0$, we can
reduce to the case where $X = Z = Y_0$ (no longer assuming $X$ and
$Z$ to be representable). Now we simply apply Lemma \ref{tech2},
noting that $L$ commutes with all colimits and with fiber products
over $Y_0$.
\end{proof}

We can now give the proof of Theorem \ref{giraud}:

\begin{proof}
We have just seen that $(1) \Rightarrow (2)$. Since prestack
$\infty$-categories are $\infty$-topoi, it is obvious that $(2)
\Rightarrow (3)$. Let us prove that $(3) \Rightarrow (1)$. By
Proposition \ref{simpsongiraud}, $\calX$ is presentable. It
remains to check the other axioms. Let $L: \calP \rightarrow
\calX$ denote the left exact localization functor.

Sums in $\calX_0$ are disjoint: given two objects $E$ and $E'$ of
$\calX_0$, their sum in $\calX_0$ is $L(E \coprod E')$, so that
the fiber product $E \times_{L(E \coprod E')} E' \simeq LE
\times_{L(E \coprod E')} LE' \simeq L(E \times_{E \coprod E'} E')
\simeq L(\emptyset)$, which is the initial object of $\calX_0$.

Finally, suppose that $X_{\bigdot}$ is a groupoid in $\calX_0$.
Then for any $Y \in \calX$, $\Hom(Y,X_\bigdot) = \Hom(LY,
X_{\bigdot})$, so that $X_{\bigdot}$ is also a groupoid in
$\calX$. Thus $X_{\bigdot}$ is effective in $\calX$. Since $L$ is
left exact, it follows that $X_{\bigdot}$ is effective in
$\calX_0$.
\end{proof}

Theorem \ref{giraud} is extremely useful. It allows us to reduce
the proofs of many statements about arbitrary $\infty$-topoi to
the case of $\infty$-categories of prestacks. We can then often
reduce to the case of the $\infty$-topos $\SSet$ by working
componentwise. This places the full apparatus of classical
homotopy theory at our disposal.

\subsection{Truncated Objects}\label{toposmade}

In this section, we will show that for any $\infty$-topos $\calX$,
the full subcategory of $\calX$ consisting of ``discrete objects''
is an (ordinary) topos. This observation gives us a functor from
the $(\infty,2)$-category of $\infty$-topoi to the $2$-category of
topoi.

We will begin with some generalities concerning truncated objects
in an $\infty$-category.

\begin{definition}
Let $\calC$ be an $\infty$-category and $k$ an integer $\geq -1$.
An object $C \in \calC$ is called {\it $k$-truncated} if for any
morphism $\eta: D \rightarrow C$, one has $\pi_{n}(
\Hom_{\calC}(D,C), \eta) = \ast$ for $n > k$. By convention we
shall say that an object is $(-2)$-truncated if it is a final
object of $\calC$. An object $C$ is called {\it discrete} if it is
$0$-truncated.
\end{definition}

In other words, an object $C \in \calC$ is discrete if
$\Hom_{\calC}(D,C)$ is a discrete set for any $D \in \calC$. We
see immediately that the discrete objects of $\calC$ form an
ordinary category.

\begin{remark}
The $k$-truncated objects in an $\infty$-category $\calC$ are
stable under the formation of all limits which exist in $\calC$.
\end{remark}

\begin{example}
Suppose $\calC = \SSet$. Then the $(-2)$-truncated objects are
precisely the contractible spaces. The $(-1)$-truncated objects
are spaces which are either empty or contractible. The
$0$-truncated objects are the spaces which are discrete (up to
homotopy). More generally, the $n$-truncated objects are the
spaces $X$ all of whose homotopy groups $\pi_k(X,x)$ vanish for $k
> n$.
\end{example}

Let us say that a morphism $f: C \rightarrow C'$ in an
$\infty$-category $\calC$ is {\it $k$-truncated} if it exhibits
$C$ as a $k$-truncated object in the slice category $\calC_{/C'}$
of objects over $C'$. For example, the $(-2)$-truncated morphisms
are precisely the equivalences in $\calC$. The $(-1)$-truncated
morphisms are called {\it monomorphisms}. A morphism $f:C
\rightarrow C'$ is a monomorphism if the induced map
$\Hom_{\calC}(X,C)\rightarrow \Hom_{\calC}(X,C')$ is equivalent to
the inclusion of a summand, for any $X \in \calC$.

The following easy lemma, whose proof is left to the reader, gives
a recursive characterization of $n$-truncated morphisms.

\begin{lemma}\label{trunc}
Let $\calC$ be an $\infty$-category with finite limits and let $k
\geq -1$. A morphism $f: C \rightarrow C'$ is $k$-truncated if and
only if the diagonal $C \rightarrow C \times_{C'} C$ is
$(k-1)$-truncated.
\end{lemma}

This immediately implies the following:

\begin{proposition}\label{eaa}
Let $F: \calC \rightarrow \calC'$ be a left-exact functor between
$\infty$-categories which admit finite limits. Then $F$ carries
$k$-truncated objects into $k$-truncated objects and $k$-truncated
morphisms into $k$-truncated morphisms.
\end{proposition}

\begin{proof}
An object $C$ is $k$-truncated if and only if the morphism $C
\rightarrow 1$ to the final object is $k$-truncated. Since $F$
preserves the final object, it suffices to prove the assertion
concerning morphisms. Since $F$ commutes with fiber products,
Lemma \ref{trunc} allows us to use induction on $k$, thereby
reducing to the case where $k=-2$. But the $(-2)$-truncated
morphisms are precisely the equivalences, and these are preserved
by any functor.
\end{proof}

Proposition \ref{eaa} applies in particular whenever $F =
f_{\ast}$ or $F=f^{\ast}$, if $f$ is a geometric morphism between
$\infty$-topoi.

\begin{proposition}\label{maketrunc}
Let $\calC$ be a presentable $\infty$-category, $k \geq -2$. Then
there exists an accessible functor $\tau_{k}: \calC \rightarrow
\calC$, together with a natural transformation $\id_{\calC}
\rightarrow \tau_{k}$, with the following property: for any $C \in
\calC$, the object $\tau_k C$ is $k$-truncated, and the natural
map $\Hom_{\calC}(\tau_k C, D) \rightarrow \Hom_{\calC}(C,D)$ is
an equivalence for any $k$-truncated object $D$.
\end{proposition}

\begin{remark}
If the $\infty$-category $\calC$ is potentially unclear in
context, then we will write $\tau^{\calC}_{k}$ for the truncation
functor in $\calC$.
\end{remark}

\begin{proof}
Suppose first that $\calC = \SSet$. In this case, the functor
$\tau_k$ is the classical ``Postnikov truncation" functor. It is
characterized by the fact that there exists a surjective map $X
\rightarrow \tau_k X$ which induces an isomorphism on $\pi_{n}$
for $n \leq k$, and $\pi_n(\tau_k X) =0$ for $n > k$, for any
choice of basepoint. If we represent spaces by fibrant simplicial
sets, then the functor $\tau_k$ can be implemented by replacing a
simplicial set $X$ by its $k$-coskeleton. In particular, we note
that $\tau_k$ is $\omega$-continuous.

Now suppose that $\calC$ is an $\infty$-category of prestacks. The
existence of $\tau_k$ in this case follows easily from the
existence of $\tau_k$ when $\calC=\SSet$: we simply work
componentwise. Once again, the functor $\tau_k$ is actually
$\omega$-continuous.

We now handle the case of a general presentable $\infty$-category
$\calC$. By Proposition \ref{simpsongiraud}, we may view $\calC$
as the image of a localization functor $L: \calC' \rightarrow
\calC'$, where $\calC'$ is an $\infty$-category of prestacks. We
will suppose that a truncation functor $\tau^{\calC'}_{k}$ has
already been constructed for $\calC'$.

We now construct a transfinite sequence of functors
$\tau^{\alpha}_{k}: \calC \rightarrow \calC$, together with a
coherent collection of natural transformations $\tau^{\alpha_0}_k
\rightarrow \ldots \rightarrow \tau^{\alpha_n}_k$ for $\alpha_0
\leq \ldots \leq \alpha_n$. Let $\tau^0_k$ be the identity
functor, and for limit ordinals $\lambda$ we let
$\tau^{\lambda}_k$ be the (filtered) colimit of $\{\tau^{\alpha}_k
\}_{\alpha < \lambda}$. Finally, let $\tau^{\alpha+1}_k = L
\tau^{\calC'}_k \tau^{\alpha}_k$.

One verifies easily by induction that for any $k$-truncated object
$D \in \calC$, the natural morphism $$\Hom_{\calC}(
\tau^{\alpha}_k C, D) \rightarrow \Hom_{\calC}( C,D)$$ is an
equivalence. Note also that each of the functors $\tau^{\alpha}_k$
is accessible. To complete the proof, we note that $L$ is
$\kappa$-continuous for some $\kappa$. Thus $\kappa$-filtered
colimits in $\calC$ agree with $\kappa$-filtered colimits in
$\calC'$, so that $\tau^{\kappa}_k C$ is a $\kappa$-filtered
colimit of $k$-truncated objects in $\calC'$, hence $k$-truncated
in $\calC'$ and thus also in $\calC$. We may then take $\tau_k =
\tau^{\kappa}_k$.
\end{proof}

\begin{remark}
The preceding proof can be simplified if $\calC$ is an
$\infty$-topos. In this case, the localization functor $L$ can be
chosen to be left exact, and therefore to preserve the property of
being $k$-truncated, so that one may simply define $\tau_k = L
\tau^{\calC'}_k$.
\end{remark}

\begin{remark}
It is immediate from the definition of $\tau_k$ that it is a
localization functor. We deduce immediately that if $\calC$ is a
presentable $\infty$-category, then the subcategory $\tau_k \calC$
of $k$-truncated objects of $\calC$ is also presentable. In
particular, the discrete objects of $\calC$ form a presentable
category.
\end{remark}

\begin{proposition}\label{compattrunc}
Let $F: \calC \rightarrow \calC'$ be a functor between presentable
$\infty$-categories, which is left-exact and colimit preserving.
Then there are natural equivalences $F(\tau_k X) \simeq \tau_k
F(X)$.
\end{proposition}

\begin{proof}
Since $F$ is left-exact, it preserves $k$-truncated objects. Thus
the natural map $FX \rightarrow F(\tau_k X)$ factors uniquely
through some map $\phi: \tau_k FX \rightarrow F(\tau_k X)$. To
show that $\phi$ is an equivalence, we must show that
$\Hom_{\calC'}(FX, W) \simeq \Hom_{\calC}( F (\tau_k X), W)$ is an
equivalence for any $k$-truncated $W$ in $\calC'$.

Since $F$ preserves colimits, it has a right adjoint $G$. Thus we
can rewrite the desired conclusion as $\Hom_{\calC}(X, GW) \simeq
\Hom_{\calC}(\tau_k X, GW)$. Since $G$ preserves all limits, $GW$
is $k$-truncated and we are done.
\end{proof}

Proposition \ref{compattrunc} applies in particular whenever $F=
f^{\ast}$, where $f$ is a geometric morphism between
$\infty$-topoi.

We will now need a few basic facts about truncated objects in an
$\infty$-topos.

\begin{proposition}\label{pushpull}
Let $\calX$ be an $\infty$-topos containing a diagram $X
\rightarrow Y \leftarrow Z$, where $Y$ and $Z$ are $k$-truncated.
Then the natural map $f: \tau_k (X \times_Y Z) \rightarrow (\tau_k
X) \times_Y Z$ is an equivalence.
\end{proposition}

\begin{proof}
By Theorem \ref{giraud}, we may view $\calX$ as a left-exact
localization of an $\infty$-category $\calP$ of prestacks. Let $L$
denote the localization functor.  Since $\tau_k^{\calX} = L
\tau_k^{\calP}$, it will suffice to prove the proposition after
$\calX$ has been replaced by $\calP$. Working componentwise, we
can reduce to the case where $\calX = \SSet$.

In order to show that $f$ is a homotopy equivalence, it will
suffice to show that $f$ induces isomorphisms on homotopy sets for
any choice of basepoint $p \in \tau_k(X \times_Y Z)$. We may lift
$p$ to a point $\widetilde{p} \in X \times_Y Z$.  Since both the
source and target of $f$ are $k$-truncated, it will suffice to
prove that $\pi_n(X \times_Y Z, \widetilde{p}) \simeq \pi_n(\tau_k
X \times_Y Z, f(p))$ is bijective for $n \leq k$. If $n=0$, this
follows by inspection, since $X$ and $\tau_k X$ have the same
connected component structure when $k \geq 0$.

If $n > 0$, we use the fiber sequences $F \rightarrow X \times_Y Z
\rightarrow X$ and $F \rightarrow (\tau_k X) \times_Y Z
\rightarrow \tau_k X$ with the same fiber (containing the point
$\widetilde{p}$). Using the long exact sequence associated to
these fibrations, we deduce the bijectivity from the (slightly
nonabelian) five lemma, using the fact that $\pi_j(X,x)
\rightarrow \pi_j(\tau_k X,x)$ is bijective for $j \leq k$ and
surjective for $j = k+1$.
\end{proof}

We can now show how every $\infty$-topos $\calX$ determines an
ordinary topos. The proof requires Corollary \ref{equivv}, which
is proven in the next section.

\begin{theorem}
Let $\calX$ be an $\infty$-topos. Then the full subcategory of
discrete objects of $\calX$ is a topos.
\end{theorem}

\begin{proof}
We have seen that these objects form a presentable category
$\tau_0 \calX$. To complete the proof, it will suffice to verify
the remainder of Giraud's axioms (see Proposition
\ref{toposdefined}), namely:

\begin{itemize}
\item Sums in $\tau_0 \calX$ are disjoint. It will
suffice to show that fiber products in $\tau_0 \calX$ agree with
the (homotopy) fiber products in $\calX$ and that sums in $\tau_0
\calX$ agree with sums in $\calX$. The first assertion follows
from the definition of a discrete object (both fiber products
enjoy the same universal property). To prove the second, it
suffices to show that a sum of discrete objects in $\calX$ is
discrete, which follows easily from the disjointness of sums in
$\calX$.

\item Colimits in $\tau_0 \calX$ commute with pullback. Since $\tau_0 \calX$
is a localization of $\calX$, colimits in $\tau_0 \calX$ may be
computed by first forming the appropriate colimit in $\calX$ and
then applying $\tau_0$. The required commutativity follows
immediately from Proposition \ref{pushpull}.

\item All equivalence relations in $\calX_1$ are effective. If $R \subseteq U \times U$
is an equivalence relation in $\tau_0 \calX$, then we obtain a
groupoid $X_{\bigdot}$ in $\tau_0 \calX$, with $X_n = R \times_U
\ldots \times_U R$. Since $\calX$ is an $\infty$-topos, this
groupoid is effective. Let $X=|X_{\bigdot}|$. To complete the
proof, it will suffice to show that $X$ is discrete. This follows
from Corollary \ref{equivv}, since $U \times_X U \rightarrow U
\times U$ is a monomorphism (by the definition of an equivalence
relation).
\end{itemize}
\end{proof}

The construction $\calX \rightsquigarrow \tau_0 \calX$ is
functorial on the $(\infty,2)$-category of $\infty$-topoi. Indeed,
if $f: \calX \rightarrow \calX'$ is a geometric morphism of
$\infty$-topoi, then $f^{\ast}$ and $f_{\ast}$ are both left exact
and therefore carry discrete objects into discrete objects. They
therefore induce (left exact) adjoint functors between $\tau_0
\calX$ and $\tau_0 \calX'$, which is precisely the definition of a
geometric morphism of topoi.

One should think of $\tau_0 \calX$ as a sort of ``Postnikov
truncation'' of $\calX$. The classical $1$-truncation of a
homotopy type $X$ remembers only the fundamental groupoid of $X$.
It therefore knows all about local systems of sets on $X$, but
nothing about fibrations over $X$ with non-discrete fibers. The
relationship between $\calX$ and $\tau_0 \calX$ is analogous:
$\tau_0 \calX$ knows about the ``sheaves of sets'' on $\calX$, but
has forgotten about sheaves with nondiscrete spaces of sections.

\begin{remark}
In view of the above, the notation $\tau_0 \calX$ is unfortunate
because the analogous notation for the $1$-truncation of a
homotopy type $X$ is $\tau_1 X$. We caution the reader not to
regard $\tau_0 \calX$ not as the result of applying an operation
$\tau_0$ to $\calX$; it instead denotes the essential image of the
functor  $\tau_0: \calX \rightarrow \calX$.
\end{remark}

We would expect, by analogy with homotopy theory, that there in
some sense a geometric morphism $\calX \rightarrow \tau_0 \calX$.
In order to properly formulate this, we need to show how to build
an $\infty$-topos from an ordinary topos. We will return to this
idea after a brief digression.

\subsection{Surjections}

Let $\calX$ be an $\infty$-topos. The assumption that every
groupoid in $\calX$ is effective leads to a good theory of
surjections (also called ``effective epimorphisms") in $\calX$, as
we shall explain in this section.

Let $U_0$ be an object in an $\infty$-topos $\calX$. We define a
groupoid object $U_{\bigdot}$ of $\calX$ by letting $U_n$ denote
the $(n+1)$-fold product of $U_0$ with itself. Let $U =
|U_{\bigdot}|$. Then we have the following:

\begin{proposition}
The morphism $U_0 \rightarrow U$ is equivalent to the natural
morphism $U_0 \rightarrow \tau_{-1} U_0$.
\end{proposition}

\begin{proof}
We first show that $U$ is $(-1)$-truncated. It suffices to show
that the diagonal map $U \rightarrow U \times U$ is an
equivalence. Note that $U = U \times_U U$. Since colimits commute
with fiber products and $U = |U_{n}|$, it suffices to show that
$p_{n}: U_{n} \times_U U_{n} \rightarrow U_{n} \times U_{n}$ is an
equivalence for all $n \geq 0$. Note that $p_{n} = U_{n}
\times_{U_0} p_0 \times_{U_0} U_{n}$; thus it suffices to consider
the case where $n=0$. In this case, the fiber product $U_0 \times
U_0 = U_1$ by definition, and $U_0 \times_{U} U_0 = U_1$ since
$U_{\bigdot}$ is effective.

To complete the proof, it suffices to show that $\Hom_{\calX}(U,E)
\rightarrow \Hom_{\calX}(U_0, E)$ is an equivalence whenever $E$
is $(-1)$-truncated. By definition, the condition on $E$ means
that both spaces are either empty or contractible; we must show
that if the target $\Hom_{\calX}(U_0, E)$ is nonempty, then so is
$\Hom_{\calX}(U,E)$. But this is clear from the description of $U$
as a colimit of objects, each of which admits a map to $E$
(automatically unique).
\end{proof}

If $U \simeq 1$, then we shall say that $U_0 \rightarrow 1$ is a
{\it surjection}. More generally, we shall say that a morphism
$U_0 \rightarrow E$ is a {\it surjection} if it is a surjection in
the $\infty$-topos $\calX_{/E}$. In other words, $U_0 \rightarrow
E$ is a surjection if $E \simeq |U_{\bigdot}|$, where $U_n$
denotes the $(n+1)$-fold fiber power of $U_0$ over $E$.

\begin{remark}

\begin{enumerate}
\item Any equivalence is surjective.
\item Any morphism which is homotopic to a surjection is
surjective.
\item If $f: \calX \rightarrow \calY$ is a geometric morphism and
$Y \rightarrow Y'$ is a surjection in $\calY$, then $f^{\ast} Y
\rightarrow f^{\ast} Y'$ is a surjection in $\calX$. In
particular, we see that $Y \times_{Y'} Y'' \rightarrow Y''$ is
surjective for any $Y'' \rightarrow Y'$.
\end{enumerate}
\end{remark}

To show that surjections are closed under composition, it is
convenient to recast the definition of a surjection in the
following way. Given an object $E \in \calX$, let $\Sub(E)$ denote
the collection of all ``subobjects" of $E$: namely, all
equivalence classes $(-1)$-truncated morphisms $E_0 \rightarrow
E$. We may partially order this collection so that $E_0 \leq E_1$
if there is a factorization (automatically unique) $E_0
\rightarrow E_1 \rightarrow E$. It is not difficult to verify that
$\Sub(E)$ is actually a set (it is the set of subobjects of $1$ in
the topos $\tau_0 \calX_{/E}$), but we shall not need this. Note
that given any morphism $f: E' \rightarrow E$, we obtain an
order-preserving pullback map $f^{\ast}: \Sub(E) \rightarrow
\Sub(E')$, given by $f^{\ast}(E'_0) = E'_0 \times_{E'} E$.

\begin{proposition}
A morphism $f: E' \rightarrow E$ in $\calX$ is surjective if and
only if $f^{\ast}: \Sub(E) \rightarrow \Sub(E')$ is injective.
\end{proposition}

\begin{proof}
Suppose first that $f^{\ast}$ is injective. Let $U_{n}$ denote the
$(n+1)$-fold fiber power of $E'$ over $E$. Then $|U_{\bigdot}|
\rightarrow E$ is $(-1)$-truncated. Thus $|U_{\bigdot}|$ and $E$
give two elements of $\Sub(E)$. It is easy to see that $f^{\ast}
|U_{\bigdot}| = f^{\ast} E \in \Sub(E')$. The injectivity of
$f^{\ast}$ then implies that $|U_{\bigdot}| \rightarrow E$ is an
equivalence, so that $f$ is surjective.

Conversely, suppose that $f$ is surjective, so that $|U_{\bigdot}|
\rightarrow E$ is an equivalence. Consider any pair of injections
$E_0 \rightarrow E$, $E_{1} \rightarrow E$. Suppose that $f^{\ast}
E_0 = f^{\ast} E_1$ in $\Sub(E')$; we must show that $E_0 = E_1$
in $\Sub(E)$. It will suffice to show that $E_0 \times_E E_1$ maps
isomorphically to both $E_0$ and $E_1$. By symmetry it will
suffice to prove this for $E_1$. Replacing $E_0$ by $E_0 \times_E
E_1$, we may suppose that there is a factorization $E_0
\stackrel{g}{\rightarrow} E_1 \rightarrow E$. We wish to show that
$g$ is an equivalence.

Since pullbacks commute with colimits, the morphism $g$ is
obtained as a colimit of morphisms $g_n: E_0 \times_{E} U_n
\rightarrow E_1 \times_E U_n$. Thus it suffices to show that each
$g_n$ is an equivalence. Since $g_n = g_0 \times_{U_0} U_{n}$, it
suffices to show that $g_0$ is an equivalence. But this is
precisely what the assumption that $f^{\ast} E_0 = f^{\ast} E_1$
tells us.
\end{proof}

From this we immediately deduce some corollaries.

\begin{corollary}\label{composite}
A composite of surjective morphisms is surjective.
\end{corollary}

\begin{corollary}\label{sums}
Let $f_{\alpha}: E'_{\alpha} \rightarrow E_{\alpha}$ be a
collection of morphisms in an $\infty$-topos $\calX$, and let $f:
E' \rightarrow E$ be their sum. Then $f$ is surjective if and only
if each $f_{\alpha}$ is surjective.
\end{corollary}

\begin{proof}
Note that $\calP(E) = \Pi_{\alpha} \calP(E_{\alpha})$ and
$\calP(E') = \Pi_{\alpha} \calP(E'_{\alpha})$. If each
$f_{\alpha}^{\ast}$ is injective, then $f^{\ast} = \Pi_{\alpha}
f_{\alpha}^{\ast}$ is also injective. The converse also holds,
since each of the factors $\calP(E_{\alpha})$ is nonempty.
\end{proof}

\begin{corollary}\label{cover}
Suppose $X \stackrel{f}{\rightarrow} Y \stackrel{g}{\rightarrow}
Z$ is a diagram with $g \circ f$ surjective. Then $g$ is
surjective.
\end{corollary}

The notion of a surjective morphism gives a good mechanism for
proving theorems by descent. The following gives a sampler:

\begin{proposition}
Let $\calX$ be an $\infty$-topos, let $f: E \rightarrow S$ be a
morphism in $\calX$, let $g: S' \rightarrow S$ be another
morphism, and let $f': E' \rightarrow S'$ be the base change of
$f$ by $S' \rightarrow S$. Then:

\begin{itemize}
\item If the morphism $f$ is $n$-truncated, then so is $f'$. The
converse holds if $S' \rightarrow S$ is surjective.
\item If the morphism $f$ is surjective, then so is $f'$. The
converse holds if $S' \rightarrow S$ is surjective.
\end{itemize}
\end{proposition}

\begin{proof}
We begin with the first claim. If $f$ is $n$-truncated, then $f'$
is clearly $n$-truncated. For the converse, we note that $f$ is
$n$-truncated if and only it induces an equivalence $E \rightarrow
\tau_n^{\calX_{/ S}} E$. Since base change by $S'$ commutes with
truncation, it will suffice to prove that a morphism which becomes
an equivalence after base change to $S'$ is an equivalence to
begin with (in other words, it suffices to consider the case where
$n=-2$).

Let $S'_{\bigdot}$ denote the groupoid obtained by taking fiber
powers of $S'$ over $S$, and let $E'_{\bigdot}$ be the groupoid
obtained by taking iterated fiber powers of $E'$ over $E$. Then if
$f'$ is an equivalence, we get an induced equivalence
$E'_{\bigdot} \rightarrow S'_{\bigdot}$, hence an equivalence
$|E'_{\bigdot} | \rightarrow | S'_{\bigdot} |$. Since $g$ is
surjective, the target is equivalent to $S$. The map $E'
\rightarrow E$ is a base change of $g$, hence surjective, so the
same argument shows that the source $|E'_{\bigdot} |$ is
equivalent to $E$, which completes the proof.

Now let us consider the second claim. We have already seen (and
used) that the surjectivity of $f$ implies the surjectivity of
$f'$. For reverse direction, use Corollary \ref{cover}.
\end{proof}

The following consequence was needed in the last section:

\begin{corollary}\label{equivv}
Let $\calC$ be an $\infty$-topos, $k \geq -1$, and let $f: U
\rightarrow X$ be a surjective morphism in $\calC$. Suppose that
$U \times_X U \rightarrow U \times U$ is $(k-1)$-truncated. Then
$X$ is $k$-truncated.
\end{corollary}

\begin{proof}
The morphism $U \times_X U \rightarrow U \times U$ is obtained
from the diagonal morphism $X \rightarrow X \times X$ by a
surjective base change. It follows that the diagonal of $X$ is
$(k-1)$-truncated, so that $X$ is $k$-truncated by Lemma
\ref{trunc}.
\end{proof}

We saw in the last section that if $\calX$ is an $\infty$-topos,
then the full subcategory $\calX_1$ consisting of discrete objects
is an ordinary topos. By Giraud's theorem \ref{giraud}, $\calX_1$
actually arises as the category of sheaves on some site equipped
with a Grothendieck topology. We will now construct such a site.

Suppose that $\calX$ is $\kappa$-coherent. Via the Yoneda
embedding, every $1$-truncated object determines a presheaf of
{\em sets} on $\calX_{\kappa}$. Note that a presheaf on an
$\infty$-category $\calC$ is the same thing as a presheaf on the
homotopy category $h \calC$. We will show that $\calX_{1}$ is
precisely the collection of sheaves of sets for a Grothendieck
topology on $h \calX_{\kappa}$ which we now describe. We will say
that a collection of morphisms $\{ E_{\alpha} \rightarrow E \}$ in
$h \calX_{\kappa}$ is a covering of $E$ if the induced map (which
is well-defined up to homotopy) $\coprod_{\alpha} E_{\alpha}
\rightarrow E$ is surjective.

\begin{proposition}
The covering families defined above determine a Grothendieck
topology on $h \calX_{\kappa}$. The presheaf represented by any
discrete object of $\calX$ is a sheaf of sets on $h
\calX_{\kappa}$. The Yoneda embedding induces an equivalence of
categories between $\calX_{1}$ and the category of sheaves of sets
on $h\calX_{\kappa}$.
\end{proposition}

\begin{proof}
Let us first show that the covering families determine a
Grothendieck topology on $h \calX_{\kappa}$.

\begin{itemize}
\item We must show that the one element family $\{E \rightarrow E\}$ is
covering for any $E$, and that any covering family remains a
covering family when it is enlarged. Both of these properties are
clear from the definitions.

\item Suppose that some family of morphisms $\{ E_{\alpha} \rightarrow E
\}$ give a covering of $E$, and the for each $\alpha$ we have a
covering $\{ E_{\alpha \beta} \rightarrow E_{\alpha} \}$. We must
show that the family of composites $\{ E_{\alpha\beta} \rightarrow
E \}$ is also a covering family. The induced morphism
$$\coproduct_{\alpha,\beta} E_{\alpha\beta} \rightarrow \coproduct_{\alpha}
E_{\alpha}$$ is surjective by Corollary \ref{sums}, hence the
morphism $$\coproduct_{\alpha,\beta} E_{\alpha\beta} \rightarrow
E$$ is surjective by Corollary \ref{composite}.

\item
Let $\{ E_{\alpha} \rightarrow E \}$ be a covering and $f: E'
\rightarrow E$ any morphism. We must show that there exists a
covering of $\{ E'_{\beta} \rightarrow E'\}$ such that each
composite map $E'_{\beta} \rightarrow E$ factors through some
$E_{\alpha}$. In fact, we may take the family of morphisms $\{
E_{\alpha} \times_{E} E' \rightarrow E' \}$. (Note that
$E_{\alpha} \times_{E} E'$ is {\em not necessarily} a fiber
product in $h \calC_{\kappa}$: these do not necessarily exist, and
are not needed).
\end{itemize}

Now suppose that $X \in \calX_1$ is a discrete object of $\calX$.
We claim that $\calF(\bigdot) = \Hom(\bigdot, X)$ is a sheaf on $h
\calC_{\kappa}$. Suppose given a covering $\{E_\alpha \rightarrow
E\}$ of an object $E$ and elements $\eta_{\alpha} \in
\calF(E_{\alpha})$. We must show that if the $\eta_{\alpha}$ agree
on overlaps, then they glue to a unique section of $\calF(E)$. Set
$E' = \coprod E_{\alpha}$ and $U_n$ denote the $(n+1)$-fold fiber
power of $E'$ over $E$. Then $U_{\bigdot}$ is a groupoid with
$|U_{\bigdot}| = E$. Consequently, we have $\Hom_{\calX}(E,X) =
|\Hom_{\calX}(U_{\bigdot},X)|$. Since $X$ is discrete, the right
hand side is simply the kernel of the pair of maps $p,q: \calF(E')
\rightarrow \calF(E' \times_E E')$. The condition that the
$\eta_{\alpha}$ agree on overlap is precisely the condition that
they give rise to an element of this kernel.

Let $\calC$ denote the topos of sheaves on $h \calX_{\kappa}$.
Since $\calX$ was assumed to be $\kappa$-accessible, the Yoneda
embedding $\pi: \calX_{1} \rightarrow \calC$ is fully faithful. To
complete the proof, it will suffice to show that $\pi$ is
essentially surjective. In other words, we must prove that any
sheaf of sets $\calF$ on $h \calC_{\kappa}$ is representable by an
object of $\calX$ (automatically discrete). Regard $\calF$ as a
presheaf on $\calC_{\kappa}$; we need to show that it is
representable. If $\calF$ carries all $\kappa$-small colimits into
limits, then it has a unique $\kappa$-continuous extension to
$\calX$ which is representable by Proposition \ref{representable}.
In order to show that $\calF$ commutes with all $\kappa$-small
limits, it will suffice to show that $\calF$ commutes with
$\kappa$-small sums and with coequalizers (see Appendix
\ref{appendixdiagram}).

If $E$ is a sum of $\kappa$-compact objects $E_{\alpha}$, then $\{
E_{\alpha} \rightarrow E \}$ is a covering family so that
$\calF(E) \rightarrow \Pi_{\alpha} \calF(E_{\alpha})$ is
injective. The surjectivity follows from the fact that sums are
disjoint in $\calX$, so that no compatibility conditions are
required to glue elements of $\calF(E_{\alpha})$ to an element of
$\calF(E)$.

Finally, suppose $E$ is presented as a homotopy coequalizer of a
diagram $\pi_0, \pi_1: F_1 \rightarrow F_0$ of $\kappa$-compact
objects. We must show that the diagram $\calF(E) \rightarrow
\calF(F_0) \stackrel{\rightarrow}{\rightarrow} \calF(F_1)$ is a
(homotopy) equalizer. Since $\calF$ is set-valued, this means only
that $\calF(E)$ is the kernel of the pair of maps $\calF(F_0)
\stackrel{\rightarrow}{\rightarrow} \calF(F_1)$. The injectivity
of $\calF(E) \rightarrow \calF(F_0)$ follows from the fact that
$\{ F_0 \rightarrow E \}$ is a covering (and that $\calF$ is a
separated presheaf).

To complete the proof, we need to show that given any $\eta$ in
the kernel of $\pi_0, \pi_1: \calF(F_0) \rightarrow \calF(F_1)$
arises from an object in $\calF(E)$. Since $F_0 \rightarrow E$ is
surjective, it suffices to show that the two restrictions of
$\eta$ to $\calF( F_0 \times_{E} F_0)$ coincide. Let $\pi'_0,
\pi'_1: F_0 \times_{E} F_0 \rightarrow F_0$ denote the two
projections.

To complete the proof, we need to recall the notation of Lemma
\ref{technical}: the $\infty$-category $\coeq_{F_0}$ has a
coherently associative multiplication $\otimes$. Let $\calC$
denote the full subcategory of $\coeq_{F_0}$ consisting of those
triples $(A,p,q)$ such that $p^{\ast} \eta = q^{\ast} \eta \in
\calF(A)$. It is clear from the definition that $\calC$ is stable
under $\otimes$ and under direct sum. Since $\calF$ is a separated
presheaf, any surjective morphism $\psi: A \rightarrow A'$ induces
an injection $\psi^{\ast}: \calF(A') \rightarrow \calF(A)$. It
follows that if $(A,p,q) \in \calC$ and there is a morphism
$(A,p,q) \rightarrow (A',p',q')$ in $\coeq_{F_0}$ which induces a
surjection $A \rightarrow A'$, then $(A',p',q') \in \calC$.

Let $R_{ij}$ denote the object $(F_1, \pi_i,\pi_j) \in \coeq_Y$.
The proof of Lemma \ref{technical} implies that we can construct
$F_0 \times_{E} F_0$ as a colimit of $\otimes$-products of the
objects $R_{ij}$. Consequently, there exists a morphism $(A,p,q)
\rightarrow (F_0 \times_{E} F_0, \pi'_0, \pi'_1)$ which induces a
surjection $A \rightarrow F_0 \times_{E} F_0$, where $A$ is a
disjoint union of $\otimes$-products of the objects $R_{ij}$. We
need to show that $(F_0 \times_{E} F_0, \pi'_0, \pi'_1) \in
\calC$. By the preceding argument, it suffices to show that each
$R_{ij}$ lies in $\calC$, which is clear.
\end{proof}

Now we come to a crucial fact which will be needed subsequently: a
map is surjective if and only if it induces a surjective map on
connected components. Note that a map between discrete objects of
$\calX$ is surjective if and only if it is a surjective map in the
ordinary topos $\tau_0 \calX$, in the usual sense.

\begin{proposition}\label{pi0detects}
Let $\calX$ be an $\infty$-topos, and let $\phi: E \rightarrow F$
be a morphism in $\calX$. Then $\phi$ is surjective if and only if
$\tau_0 \phi: \tau_0 E \rightarrow \tau_0 F$ is surjective.
\end{proposition}

\begin{proof}
Suppose that $\phi$ is surjective, so that $F \simeq
|E_{\bigdot}|$ where $E_{n}$ denotes the $(n+1)$-fold fiber power
of $E$ over $F$. Since $\tau_0: \calX \rightarrow \tau_0 \calX$
commutes with colimits (it is a left adjoint), we deduce that
$\tau_{0} F$ is the geometric realization of $\tau_0 E_{\bigdot}$
in the ordinary category $\calX_{1}$. In other words, it is the
coequalizer of a pair of morphisms $\tau_0 E_1 \rightarrow \tau_0
E_0$. This proves that $\tau_0 E \rightarrow \tau_0 F$ is
surjective.

For the converse, let us suppose that $\tau_0 E \rightarrow \tau_0
F$ is surjective. Choose $\kappa$ such that $\calX$ is
$\kappa$-coherent and $E,F \in \calX_{\kappa}$. Then $E$
represents a prestack on $\calX_{\kappa}$; let $hE$ denote the
corresponding presheaf on $h \calX_{\kappa}$, defined by $hE(X) =
\pi_0 E(X) = \pi_0 \Hom_{\calX}(X,E)$. Similarly let $hF$ denote
the presheaf corresponding to $F$. Then $\tau_0 E$ and $\tau_0 F$
are the sheafifications of these presheaves. The identity map $F
\rightarrow F$ gives an element of $\eta \in hF(F)$. Since $\tau_0
E$ surjects onto $\tau_0 F$, there exists a covering $\{
F_{\alpha} \rightarrow F \}$ and liftings $\eta_{\alpha} \in
hE(F_{\alpha})$. Enlarging $\kappa$ if necessary, we may assume
that the sum $F'$ of the $F_{\alpha}$ is $\kappa$-compact, so that
$\eta$ lifts to $\eta' \in \pi_0 \Hom_{\calX}(F',E)$ with $F'
\rightarrow F$ surjective. By Corollary \ref{cover}, $\phi$ is
also surjective.
\end{proof}

\subsection{Stacks on a Topos}

In this section, we introduce a construction which passes from
ordinary topoi to $\infty$-topoi, and show that it is right
adjoint to the functor $\calX \rightsquigarrow \tau_0 \calX$
constructed in \S \ref{toposmade}. The construction is in some
sense dual to what we did in the last section: instead of starting
with an $\infty$-topos and producing a site, we will begin with a
site and produce an $\infty$-topos of ``sheaves of spaces" on that
site.

\begin{definition}
Let $\toposX$ be a topos, and $\calF$ a prestack on $\toposX$. We
will say that $\calF$ is a {\it stack} on $\toposX$ if it
satisfies the following descent conditions:
\begin{itemize}
\item Given any collection $\{C_i\}$ of objects of $\toposX$ with
(disjoint) sum $C$, the natural map $\calF(C) \rightarrow \Pi_i
\calF(C_i)$ is an equivalence.

\item Given any surjection $C \rightarrow D$ in $\toposX$, let
$C_{\bigdot}$ denote the simplicial object of $\toposX$ with
$C_{n}$ the $(n+1)$-fold fiber power of $C$ over $D$. Then the
natural map $\calF(D) \rightarrow |\calF(C_{\bigdot})|$ is an
equivalence. (Here $|\calF(C_{\bigdot})|$ denotes the geometric
realization of the cosimplicial space).

\end{itemize}
\end{definition}

\begin{lemma}\label{eqq}
Let $\toposX$ be a topos, and let $\toposX_0$ be a full
subcategory of $\toposX$ which is closed under finite limits and
generates $\toposX$. Let $f: \calF \rightarrow \calG$ be a map of
sheaves on $\toposX$. If $f$ induces an equivalence $\calF(C)
\rightarrow \calG(C)$ for all $C \in X_0$, then $f$ is an
equivalence. \end{lemma}

\begin{proof}
Without loss of generality, we may enlarge $\toposX_0$ so that
$\toposX_0$ is closed under sums: $\toposX_0$ remains stable under
the formation of fiber products, and the value of a sheaf on a sum
is determined by its values on the summands.

We must show that for all $C \in \toposX$, $\calF(C) \rightarrow
\calG(C)$ is an equivalence. By hypothesis, this is true when $C
\in \toposX_0$. Next, suppose that there is a monomorphism $C
\rightarrow C'$, with $C' \in \toposX_0$. Choose a surjection $U
\rightarrow C$, where $U \in \toposX_0$. Then we have $U
\times_{C} U \times_{C} \ldots \times_{C} U = U \times_{C'} U
\times_{C'} \ldots \times_{C'} U$, so that all fiber powers of $U$
over $C$ lie in $\toposX_0$. Let $U_{n}$ be the $(n+1)$-fold fiber
power of $U$ over $C$. Then $\calF(C) \simeq | \calF(U_{\bigdot})
|$, and similarly for $\calG$. Since $f$ induces homotopy
equivalences $\calF(U_n) \rightarrow \calG(U_n)$, it induces a
homotopy equivalence between the geometric realizations $\calF(C)
\rightarrow \calG(C)$.

Now consider the general case. As before, we choose a surjection
$U \rightarrow C$, where $U \in \toposX_0$. In this case, each
fiber power $U \times_C U \times_C \ldots \times_C U$ is a
subobject of $U \times U \times \ldots U$ which belongs to
$\toposX_0$. The above argument shows that $\calF(U_n) \rightarrow
\calF'(U_n)$ is an equivalence for each $n$, and once again we
obtain an equivalence $\calF(C) \rightarrow \calF'(C)$.
\end{proof}

\begin{proposition}\label{plus}
Let $\toposX$ be a topos. Then the $\infty$-category of stacks on
$\toposX$ is an $\infty$-topos.
\end{proposition}

\begin{proof}
Let $\calC$ be a small category with finite limits and equipped
with a Grothendieck topology such that $\toposX$ is equivalent to
the category of sheaves of sets on $\calC$. We will obtain the
$\infty$-category of sheaves on $\toposX$ as a left-exact
localization of the $\infty$-category $\calP = \SSet^{\calC^{op}}$
of prestacks on $\calC$.

Let $\calF$ be a prestack on $\calC$. We define a new presheaf
$\calF^{+}$ as follows: $\calF^{+}(U) = \colim_{\calS} \lim_{S \in
\calS} \calF(S)$, where $\calS$ ranges over the (directed)
collection of sieves in $\calC$ which cover $U$. Note that there
is a morphism $\calF \rightarrow \calF^{+}$ which is natural in
$\calF$. Consequently we may obtain a transfinite sequence of
functors $s_{\alpha}$ by iterating the construction $\calF \mapsto
\calF^{+}$. Let $s_{0}(\calF) = \calF$, $s_{\alpha+1}(\calF) =
s_{\alpha}(\calF)^{+}$, and $s_{\lambda}(\calF) = \colim
\{s_{\alpha}(\calF) \}_{\alpha < \lambda}$ when $\lambda$ is a
limit ordinal. One verifies easily by induction that each
$s_{\alpha}$ is left-exact and accessible.

Let us call a prestack $\calF$ on $\calC$ a {\it stack} if, for
any sieve $\calS$ in $\calC$ covering $U$, the natural map
$$\calF(U) \rightarrow \lim_{S \in \calS} \calF(S)$$ is an
equivalence. It follows easily by induction that for any prestack
$\calF$ and any stack $\calF'$, the natural map
$$\Hom(s_{\alpha} \calF, \calF') \rightarrow \Hom( \calF,
\calF')$$ is an equivalence for every ordinal $\alpha$. We also
note that for $\kappa$ sufficiently large, $s_{\kappa}(\calF)$ is
a stack for {\em any} prestack $\calF$. It follows immediately
that $s_{\kappa}$ is a localization functor, whose essential image
consists of the sheaves on $\calC$. Since $s_{\kappa}$ is
left-exact, we deduce that the $\infty$-category $\calX$ of stacks
on $\calC$ is an $\infty$-topos.

By construction, the discrete objects of $\calX$ are precisely the
sheaves of {\em sets} on $\calC$. Thus $\tau_0 \calX$ is
equivalent to the topos $\toposX$. Let $\calX'$ denote the
$\infty$-category of sheaves on $\toposX$. The restricted Yoneda
embedding induces a functor $F: \calX \rightarrow \calX'$. To
complete the proof, it will suffice to show that $F$ is an
equivalence.

We next construct a left adjoint $G$ to $F$. Given a sheaf $\calF$
Given a sheaf $\calF$ on $X$ and an object $U \in \calC$, let
$G\calF(U)$ denote the value of $\calF$ on the sheafification of
$U$. One verifies readily that $G: \calX' \rightarrow \calX$ is
left adjoint to $F$. It now suffices to show that the adjunction
morphisms $\id_{\calX'} \rightarrow FG$ and $GF \rightarrow
\id_{\calX}$ are equivalences.

It is clear from the definitions that $GF \rightarrow \id_{\calX}$
is an equivalence. For the other adjunction, we must show that if
$\calF$ is a sheaf on $X$ and $\calF'$ denotes the induced sheaf
$G\calF$ on $\calC$, then $\calF(W) \rightarrow \Hom_{\calX}(W,
\calF')$ is an equivalence. The right hand side here may be
written as a homotopy limit
$$\lim_{U \in \calC} \Hom(W(U), \calF'(U)) \simeq \lim_{(U,\eta)}
\calF'(U),$$ where the latter homotopy limit is taken over all
pairs $U \in \calC$, $\eta \in W(U)$. If $W$ is the sheafification
of an object $U$ of $\calC$, then this category has a final object
and hence the homotopy limit is given by $\calF'(U) = \calF(W)$,
as desired.

Now we have a natural morphism $\calF \rightarrow FG(\calF)$,
which induces equivalences when evaluated at any $W$ which is the
sheafification of a representable presheaf. Since every object of
$\toposX$ is a colimit of sheafifications of representable
presheaves, we may conclude using Lemma \ref{eqq}.
\end{proof}

If $\toposX$ is a topos, we will denote the $\infty$-topos of
stacks on $\toposX$ by $\stacks \toposX$. We now formulate a
universal property enjoyed by $\stacks \toposX$:

\begin{proposition}\label{universality}
Let $\toposX$ be a topos and $\calY$ an $\infty$-topos. The
category of geometric morphisms $\Top( \tau_0 \calY, \toposX)$ is
naturally equivalent to the $\infty$-category of geometric
morphisms $\Top^{\infty}(\calY, \stacks \toposX)$.
\end{proposition}

\begin{proof}
Let $\calC$ be a small category with finite limits, equipped with
a Grothendieck topology such that $\toposX$ is equivalent to the
category of sheaves of sets on $\calC$. Then the proof of
Proposition \ref{plus} shows that we may identify $\stacks
\toposX$ with the $\infty$-category of sheaves on $\calC$.

We first assume that the topology on $\calC$ is discrete, in the
sense that any sieve which covers an object $C \in \calC$ actually
contains the identity map $C \rightarrow C$. In this case,
geometric morphisms from $\calY$ into $\stacks \toposX$ are just
given by left-exact functors $f^{\ast}: \calC \rightarrow \calY$.
A similar argument shows that $\Top(\tau_0 \calY,\toposX)$ is
equivalent to the category of left exact functors $\calC
\rightarrow \tau_0 \calY$. Note that every object of $\calC$ is
discrete, and left-exact functors preserve discrete objects by
Lemma \ref{eaa}. Therefore every left exact functor $\calC
\rightarrow \calY$ factors uniquely through $\tau_0 \calY$, which
completes the proof in this case.

Now we consider the general case. Let $\calP$ denote the
$\infty$-category of presheaves on $\calC$, so that $\stacks
\toposX$ is a left-exact localization of $\calP$ and $\toposX$ is
a left-exact localization of $\tau_0 \calP$. We will denote both
localization functors by $L$. The case that we have just treated
shows that $\Top(\tau_0 \calY, \tau_0 \calP) \simeq \Top(\calY,
\calP)$. On the other hand, $\Top(\calY, \stacks \toposX)$ can be
identified with the full subcategory of $\Top(\calY, \calP)$
consisting of left-exact, colimit preserving functors $f^{\ast}:
\calP \rightarrow \calY$ which factor through $L$, in the sense
that the natural morphism $f^{\ast}(E) \rightarrow f^{\ast}(LE)$
is an equivalence for all $E \in \calP$. Similarly, $\Top(\tau_0
\calY, \toposX)$ can be identified with the collection of all
left-exact, colimit preserving functors $f^{\ast}: \tau_0 \calP
\rightarrow \tau_0 \calY$ which factor through $L$. To complete
the proof, it suffices to show that a geometric morphism
$f^{\ast}: \calP \rightarrow \calY$ factors through $L$ if and
only if the restriction $f^{\ast}|\tau_0 \calP$ factors through
$L$.

The ``only if" part is trivial, so let us focus on the ``if". Let
$S$ denote the collection of all morphisms of the form $U
\rightarrow LU$, where $U \in \tau_0 \calP$. Since every morphism
of $S$ becomes invertible after applying $f^{\ast}$, there exists
a factorization $\calP \rightarrow S^{-1} \calP \rightarrow
\calY$. It now suffices to verify that $S^{-1} \calP \rightarrow
\calY$ is an equivalence. In other words, we must show that every
$S$-local object of $\calP$ is actual $L$-local. This follows
immediately from the definition of $L$.
\end{proof}

\begin{remark}
It is possible to give a definition of $n$-topos for any $0 \leq n
\leq \infty$: this is a particular kind of $n$-category, which
``looks like'' the $n$-category of $(n-1)$-truncated stacks on a
topological space $X$. For $n = \infty$, one recovers Definition
\ref{toposs}. For $n = 1$, one recovers the classical definition
of a topos. For $n = 0$, one recovers the notion of a locale.

For each $n \geq m$, one has a ``truncation'' functor from
$n$-topoi to $m$-topoi, which replaces an $n$-topos by its full
subcategory of $m$-truncated objects. Each of these truncation
functors has a right adjoint. We have discussed the situation here
when $n = \infty$, $m=1$. When $n=1$, $m=0$, the ``truncation''
replaces a topos by the locale consisting of its ``open subsets''
(that is, subobjects of the final object).
\end{remark}

\subsection{Homotopy Groups}\label{homotopysheaves}

In this section, we will discuss the homotopy groups of objects
and morphisms in an $\infty$-topos. These will be needed in \S
\ref{hyperstack} on hyperdescent and in \S \ref{dimension} on
dimension theory, but are not needed for the main result of \S
\ref{paracompactness}.

Let $\calX$ be an $\infty$-topos, and let $X$ be an object of
$\calX$. Fixing a base point $\ast \in S^n$, we get an
``evaluation'' map $X^{S^n} \rightarrow X$. Let $\pi_n(X) =
\tau^{\calX_{/X}}_0(X^{S^n})$, where the notation indicates that
we consider $X^{S^n}$ as an object of the $\infty$-topos
$\calX_{/X}$. This is a sheaf of pointed sets on $X$: the base
point is induced by the retraction of $S^n$ onto $\ast$. As in the
classical case, the co-$H$-space structure on the sphere gives
rise to a group structure on the sheaf $\pi_n(X)$ if $n \geq 1$,
which is abelian if $n \geq 2$.

In order to work effectively with homotopy sets, it is convenient
to define the homotopy sets $\pi_n(f)$ of a morphism $f: X
\rightarrow Y$ to be the homotopy sets of $X$, considered via $f$
as an object of $\calX_{/Y}$. This construction gives again
sheaves of pointed sets on $X$ (groups if $n \geq 1$, abelian
groups if $n \geq 2$). The intuition is that the stalk of these
sheaves at a ``point" $p \in X$ is the $n$th homotopy group of the
mapping fiber of $f$, with base point $p$.

It will be useful to have the following recursive definition of
homotopy groups. Regarding $X$ as an object of the topos
$\calX_{/Y}$, we may take its $0$th truncation
$\tau_0^{\calX_{/Y}} X$. This is a sheaf of sets on $\calX_{/Y}$;
by definition we set $\pi_0(f) = f^{\ast} \tau_0^{\calX_{/Y}}(X) =
X \times_{Y} \tau_0^{\calX_{/Y}}(X)$. The natural map $X
\rightarrow \tau_0^{\calX_{/Y}}(X)$ gives a global section of
$\pi_0(f)$. Note that in this case, $\pi_0(f)$ is the pullback of
a sheaf on $Y$: this is because the definition of $\pi_0$ does not
require a base point. If $n \geq 0$, then $\pi_{n}(f) =
\pi_{n-1}(f')$, where $f': X \rightarrow X \times_{Y} X$ is the
diagonal map. Finally, in the case where $f= X \rightarrow 1$, we
let $\pi_{n}(X) = \pi_{n}(f)$.

\begin{remark}\label{emmy}
Let $f: \calX \rightarrow \calY$ be a geometric morphism of
$\infty$-topoi, and let $g: Y \rightarrow Y'$ be a morphism in
$\calY$. Then $f^{\ast}( \pi_n(g) ) \simeq \pi_n(f^{\ast}(g))$.
This follows immediately from the definition and Proposition
\ref{compattrunc}.
\end{remark}

\begin{remark}\label{sequence}
Given a pair of composable morphisms $X \stackrel{f}{\rightarrow}
Y \stackrel{g}{\rightarrow} Z$, we have a sequence of pointed
sheaves
$$\ldots \rightarrow f^{\ast} \pi_{i+1}(g) \rightarrow \pi_i(f) \rightarrow \pi_i(g \circ f) \rightarrow f^{\ast}
\pi_i(g) \rightarrow \pi_{i-1}(f) \rightarrow \ldots$$ with the
usual exactness properties. To see this, one first mimics the
usual construction to produce the boundary maps in the above
sequence. To prove exactness, one can embed the $\infty$-topos
$\calX$ in an $\infty$-topos $\calP$ of presheaves. Using Remark
\ref{emmy}, we reduce to the problem of showing that the
corresponding sequence is exact in $\calP$. This can be checked
componentwise. We are thereby reduced to the case $\calX =\SSet$,
which is classical.
\end{remark}

\begin{remark}
If $\calX = \SSet$, and $\eta: 1 \rightarrow X$ is a pointed
space, then $\eta^{\ast} \pi_{n}(X)$ is the $n$th homotopy group
of $X$ with base point $\eta$.
\end{remark}

We now study the implications of the vanishing of homotopy groups.

\begin{proposition}\label{ditz}
Let $f: X \rightarrow Y$ be an $n$-truncated morphism. Then
$\pi_{k}(f) = \ast$ for all $k > n$. If $n \geq 0$ and $\pi_{n}(f)
= \ast$, then $f$ is $(n-1)$-truncated.
\end{proposition}

\begin{proof}
The proof goes by induction on $n$. If $n = -2$, then $f$ is an
equivalence and there is nothing to prove. Otherwise, we know that
$f': X \rightarrow X \times_Y X$ is $(n-1)$-truncated. The
inductive hypothesis then allows us to infer that $\pi_{k}(f) =
\pi_{k-1}(f') = \ast$ whenever $k > n$ and $k > 0$. Similarly, if
$n \geq 1$ and $\pi_{n}(f) = \pi_{n-1}(f') = \ast$, then $f'$ is
$(n-2)$-truncated by the inductive hypothesis, so that $f$ is
$(n-1)$-truncated.

The case of small $k$ and $n$ requires special attention: we must
show that if $f$ is $0$-truncated, then $f$ is $(-1)$-truncated if
and only if $\pi_0(f) = \ast$. The fact that $f$ is $0$-truncated
implies that $\tau_0^{\calX_{Y}} X = X$, so that $\pi_0(f) = X
\times_Y X$. To say $\pi_0(f) = \ast$ is to assert that the map
$f': X \rightarrow \pi_0(f)$ is an equivalence, which is to say
that $X$ is $(-1)$-truncated.
\end{proof}

\begin{remark}
The Proposition \ref{ditz} implies that if $f$ is $n$-truncated
for $n \gg 0$, then we can test whether or not $f$ is
$n$-truncated for any particular value of $n$ by computing the
homotopy groups of $f$. In contrast to the classical situation, it
is not always possible to drop the assumption that $f$ is
$n$-truncated for $n \gg 0$.
\end{remark}

\begin{lemma}\label{truncatepin}
Let $X$ be an object in an $\infty$-topos $\calX$. Then the
natural map $p: X \rightarrow \tau_{n} X$ induces isomorphisms
$\pi_{k}(X) \rightarrow p^{\ast} \pi_{k}(\tau_n X)$ for all $k
\leq n$.
\end{lemma}

\begin{proof}
Let $\phi: \calX \rightarrow \calY$ be a geometric morphism such
that $\phi_{\ast}$ is fully faithful. By Proposition
\ref{compattrunc} and Remark \ref{emmy}, it will suffice to prove
the lemma in the case where $\calX = \calY$. By Proposition
\ref{giraud}, we may assume that $\calY$ is an $\infty$-category
of presheaves. In this case, homotopy groups and truncations are
computed pointwise. Thus we may reduce to the case $\calX =
\SSet$, where the conclusion is evident.
\end{proof}

\begin{definition}
Let $f: X \rightarrow Y$ be a morphism in an $\infty$-topos
$\calX$. We shall say that $f$ is {\it $(-1)$-connected} if $f$ is
surjective. If $0 \leq n \leq \infty$, we shall say that $f$ is
{\it $n$-connected} if it is surjective and $\pi_k(f) = \ast$ for
$k \leq n$. We shall say that the object $X$ is {\it
$n$-connected} if $f: X \rightarrow 1$ is $n$-connected.
\end{definition}

\begin{proposition}\label{goober}
Let $X$ be an object in an $\infty$-topos $\calX$. Then $X$ is
$n$-connected if and only if $\tau_n X \rightarrow 1$ is an
equivalence.
\end{proposition}

\begin{proof}
The proof goes by induction on $n \geq -1$. If $n=-1$, then the
conclusion holds by definition. Suppose $n \geq 0$. Let $p: X
\rightarrow \tau_n X$ denote the natural map. If $\tau_n X \simeq
1$, then $\pi_{k} X = p^{\ast} \pi_k(\tau_n X) = p^{\ast}
(\pi_k(1)) = \ast$ for $k \leq n$ by Lemma \ref{truncatepin}.
Conversely, suppose that $X$ is $n$-connected. Then $p^{\ast}
\pi_n(\tau_n X) = \ast$. Since $p$ is surjective, we deduce that
$\pi_n(\tau_n X) = \ast$. This implies that $\tau_n X$ is
$(n-1)$-truncated, so that $\tau_n X \simeq \tau_{n-1} X$.
Repeating this argument, we reduce to the case where $n=-1$ which
was handled above.
\end{proof}

Since $\Hom_{\calX}(X,Y) \simeq \Hom_{\calX}(\tau_n X, Y)$ when
$Y$ is $n$-truncated, we deduce that $X$ is $n$-connected if and
only if $\Hom_{\calX}(1,Y) \rightarrow \Hom_{\calX}(X,Y)$ is an
equivalence for all $n$-truncated $Y$. From this, we can
immediately deduce the following relative version of Proposition
\ref{goober}:

\begin{corollary}\label{goober2}
Let $f: X \rightarrow X'$ be a morphism in an $\infty$-topos
$\calX$. Then $f$ is $n$-connected if and only if $\Hom_{X'}(X',Y)
\rightarrow \Hom_{X'}(X,Y)$ is an equivalence, for any
$n$-truncated $Y \rightarrow X'$.
\end{corollary}

\begin{corollary}
The class of $n$-connected morphisms in an $\infty$-topos is
stable under cobase extension.
\end{corollary}

\begin{proof}
This follows from the preceding corollary and the fact that
$n$-truncated morphisms are stable under base extension.
\end{proof}

A morphism $f: X \rightarrow X'$ is $\infty$-connected if and only
if it is $n$-connected for all $n$. By the preceding corollary,
this holds if and only if $\Hom_{X'}(X',Y) \simeq \Hom_{X'}(X,Y)$
for any morphism $Y \rightarrow X'$ which is $n$-truncated for
some $n$. In $\SSet$, this implies that $f$ is an equivalence, but
this is not true in general. We will investigate this issue
further in the next section.

We conclude by noting the following stability properties of the
class of $n$-connected morphisms:

\begin{proposition}
Let $\calX$ be an $\infty$-topos.
\begin{enumerate}
\item Any $n$-connected morphism of $\calX$ is $m$-connected for
any $m \leq n$.
\item Any equivalence is $\infty$-connected.
\item The class of $n$-connected morphisms is closed under
composition.
\item Let $f: X \rightarrow Y$ be an $n$-connected morphism and
let $Y' \rightarrow Y$ be any map. Then the induced map $f': X' =
X \times_Y Y' \rightarrow Y'$ is $n$-connected. The converse holds
if $Y' \rightarrow Y$ is surjective.
\end{enumerate}
\end{proposition}

\begin{proof}
The first two claims are obvious. The third follows from the long
exact sequence of Remark \ref{sequence}. The first assertion of
$(4)$ follows from the fact that homotopy groups are compatible
with base change (a special case of Remark \ref{emmy}). The second
assertion of $(4)$ follows from Remark \ref{emmy} and descent for
equivalences ($\pi_k(f') \simeq Y'$ if and only if $\pi_k(f)
\simeq Y$).
\end{proof}

\subsection{Hyperstacks}\label{hyperstack}

If $\calC$ is a small category equipped with a Grothendieck
topology, then Jardine (see \cite{jardine}) shows that the
category of simplicial presheaves on $\calC$ is equipped with the
structure of a simplicial model category, which gives rise to an
$\infty$-category. On the other hand, we have associated to
$\calC$ an $\infty$-category $\calX$ of ``stacks on $\calC$".
These two constructions are not equivalent: the weak equivalences
in Jardine's model structure are the ``local homotopy
equivalences'', which correspond to the $\infty$-connected
morphisms in our setting. However, it turns out the Joyal-Jardine
theory can be obtained from ours by inverting the
$\infty$-connected morphisms. The resulting localization functor
is left exact, so that the $\infty$-category underlying Jardine's
model structure is also an $\infty$-topos.

Let us begin by studying the $\infty$-connected morphisms:

\begin{proposition}\label{goober3}
Let $\calX$ be an $\infty$-topos, and let $S$ denote the
collection of $\infty$-connected morphisms of $\calX$. Then $S$ is
saturated. \end{proposition}

\begin{proof}
This follows immediately from the description of
$\infty$-connected morphisms given in Corollary \ref{goober2}.
\end{proof}

In order to construct a localization of $\calX$ which inverts the
$\infty$-connected morphisms, we will need to show that the class
of $\infty$-connected morphisms is setwise generated. This follows
immediately from the following:

\begin{proposition}\label{card}
Let $\calX$ be an $\infty$-topos. Then there exists a cardinal
$\kappa$ such that any $\infty$-connected morphism $f: X
\rightarrow Y$ may be written as a $\kappa$-filtered colimit of
$\infty$-connected morphisms $f_{\alpha}: X_{\alpha} \rightarrow
Y_{\alpha}$ between $\kappa$-compact objects.
\end{proposition}

\begin{proof}
The proof is routine cardinality estimation, and may be safely
omitted by the reader. We choose $\kappa_0
> \omega$ large enough that $\calX$ is $\kappa_0$-coherent, the
functor $\tau_0$ is $\kappa$-continuous and preserves
$\kappa_0$-compact objects. Since homotopy sets are constructed
using finite limits and $\tau_0$, it follows that the formation of
homotopy sets in $\calX$ preserves $\kappa_0$-compactness and
commutes with $\kappa_0$-filtered colimits. Now choose $\kappa \gg
\kappa_0$.

Since $\calX$ is $\kappa$-coherent, we may write $f$ as a
$\kappa$-filtered colimit of morphisms $f_{\alpha}: X_{\alpha}
\rightarrow Y_{\alpha}$ where each $X_{\alpha}$ and each
$Y_{\alpha}$ is $\kappa$-compact. For simplicity, let us assume
that the colimit is indexed by a partially ordered set $\calI$.
Let $\calP$ denote the collection of all $\kappa_0$-filtered,
$\kappa$-small subsets of $\calI$ (ordered by inclusion). For each
$S \in \calP$, we let $f_S: X_S \rightarrow Y_S$ denote the
colimit indexed by $S$. Let $\calP_0 \subseteq \calP$ denote the
collection of all subsets $S$ such that $f_S$ is
$\infty$-connected. It is easy to see that $\calP$ is
$\kappa$-filtered, and the colimit of the $f_S$ is $f$. To
complete the proof, it suffices to show that $\calP_0$ is cofinal
in $\calP$.

Given any $\alpha \in \calI$ and any $k \geq 0$, we know that
$\pi_k(f_{\alpha}) \rightarrow \pi_k(f)$ factors through $X$
(since $\pi_k(f) \simeq X$). Since the source is $\kappa$-compact,
we deduce that there exists $\beta > \alpha$ such that
$\pi_i(f_{\alpha}) \rightarrow \pi_i(f_{\beta})$ factors through
$X_{\beta}$.

Given $S \in \calP$, we define a transfinite sequence of
enlargements of $S$ as follows. Let $S_0 = S$, let $S_{\lambda} =
\bigcup_{\gamma < \lambda} S_{\gamma}$ when $\lambda$ is a limit
ordinal. Finally let $S_{\lambda+1}$ be some $\kappa_0$-filtered,
$\kappa$-small subset of $\calI$ containing $S_i$ with the
property that for any $\alpha \in S_i$ and any $k \geq 0$, there
exists $\beta \in S_{i+1}$ such that $\beta \geq \alpha$ and
$\pi_k( f_{\alpha} ) \rightarrow \pi_k( f_{\beta})$ factors
through the base point. It is clear that $S_{\lambda}$ is
well-defined for $\lambda < \kappa$, and that $S_{\kappa_0}$ is
$\kappa$-filtered. Since the formation of homotopy sets commutes
with $\kappa$-filtered colimits, we deduce that $f_{S_{\kappa_0}}$
is $\infty$-connected, so that $S_{\kappa_0} \in \calP_0$ as
desired.
\end{proof}

Let $\calX$ be an $\infty$-topos, and let $S$ be the class of
$\infty$-connected morphisms of $\calX$. Propositions
\ref{goober3} and \ref{card} imply that $S$ is saturated and
setwise generated, so that we may construct a localization functor
$L: \calX \rightarrow \calX$ which inverts precisely the morphisms
in $S$.

\begin{proposition}\label{carder}
Let $\calC$ be a presentable $\infty$-category, and let $S$ denote
a setwise generated, saturated collection of morphisms of $\calC$.
Let $L: \calC \rightarrow \calC$ denote the corresponding
localization functor. The following are equivalent:

\begin{enumerate}
\item The class of morphisms $S$ is stable under base change.
\item The localization functor $L$ is left exact.
\end{enumerate}
\end{proposition}

\begin{proof}
Since $S$ is precisely the collection of morphisms $f$ such that
$Lf$ is an equivalence, it is immediate that $(2)$ implies $(1)$.
Let us therefore assume $(1)$. Since the final object $1 \in
\calX$ is obviously $S$-local, we have $L1 \simeq 1$. Thus it will
suffice to show that $L$ commutes with pullbacks. In other words,
we must show that $f: X \times_Y Z \rightarrow LX \times_{LY} LZ$
is an $S$-localization for any pair of objects $X,Z$ over $Y$ in
$\calX$. Since $LX \times_{LY} LZ$ is clearly $S$-local, it
suffices to prove that $f$ belongs to $S$. Now $f$ is a composite
of maps $$X \times_Y Z \rightarrow X \times_{LY} Z \rightarrow LX
\times_{LY} Z \rightarrow LX \times_{LY} LZ.$$ The last two maps
are obtained from $X \rightarrow LX$ and $Z \rightarrow LZ$ by
base change. It follows that they belong to $S$. Thus, it will
suffice to show that $f': X \times_Y Z \rightarrow X \times_{LY}
Z$ belongs to $S$. Since this map is obtained from $Y \times_Y Y
\rightarrow Y \times_{LY} Y$ by a base-change, it suffices to
prove that $f'': Y \rightarrow Y \times_{LY} Y$ belongs to $S$.
Projection to the first factor gives a section $s: Y \times_{LY} Y
\rightarrow Y$ of $f''$, so it suffices to prove that $s$ belongs
to $S$. But $s$ is a base change of the morphism $Y \rightarrow
LY$.
\end{proof}

Since the formation of homotopy sheaves is compatible with base
change, we see that the class $S$ of $\infty$-connected morphisms
satisfies $(1)$ of Proposition \ref{carder}. It follows from
Theorem \ref{giraud} that the collection of $S$-local objects of
$\calX$ forms an $\infty$-topos, which we shall denote by
$\calX^{\red}$. In a moment we shall characterize $\calX^{\red}$
by a universal property. First we need a simple lemma.

\begin{lemma}
Let $f: X \rightarrow Y$ be a morphism in $\calX^{\red}$. We may
form $\pi_n(f)$ in $\calX$. Then $\pi_n(f)$ lies in
$\calX^{\red}$, and is the $n$th homotopy sheaf of $f$ in
$\calX^{\red}$.
\end{lemma}

\begin{proof}
Homotopy sheaves are constructed using finite limits and $\tau_0$.
Since $\calX^{\red}$ is stable under finite limits, it suffices to
show that for any $g: E \rightarrow E'$ in $\calX^{\red}$, we have
$\tau_0^{\calX_{E'}} E \in \calX^{\red}$. This follows from the
fact that $\tau_0^{\calX_{E'}}$ is discrete over $E'$ and $E'$ is
$S$-local.
\end{proof}

Following \cite{toen}, we shall say that an $\infty$-topos $\calX$
{\it $t$-complete} if every $\infty$-connected morphism of $\calX$
is an equivalence.

\begin{lemma}
The $\infty$-topos $\calX^{\red}$ is $t$-complete.
\end{lemma}

\begin{proof}
Since the formation of homotopy sheaves in $\calX^{\red}$ is
compatible with the formation of homotopy sheaves in $\calX$, any
morphism $f: X \rightarrow Y$ which is $\infty$-connected in
$\calX^{\red}$ is also $\infty$-connected in $\calX$, and
therefore an equivalence in $\calX^{\red}$. Thus $\calX^{\red}$ is
$t$-complete.
\end{proof}

We can give the following universal characterization of
$\calX^{\red}$.

\begin{proposition}
Let $\calX$ be an $\infty$-topos and $\calY$ an $\infty$-topos
which is $t$-complete. Then composition with the natural geometric
morphism $\calX^{\red} \rightarrow \calX$ induces an equivalence
of $\infty$-categories
$$\Top^{\infty}( \calY, \calX^{\red}) \rightarrow \Top^{\infty}(\calY, \calX).$$
\end{proposition}

\begin{proof}
In other words, we must show that if $f: \calY \rightarrow \calX$
is a geometric morphism, then $f^{\ast}$ factors through $L$: that
is, it carries $\infty$-connected morphisms into equivalences.
Since $\calY$ is $t$-complete, it suffices to show that $f^{\ast}$
preserves $\infty$-connected morphisms. This follows immediately
from the fact that $f^{\ast}$ commutes with the formation of
homotopy sheaves.
\end{proof}

The objects of $\calX^{\red}$ are the $S$-local objects of
$\calX$. It is of interest to describe these objects explicitly.
This description was given in \cite{hollander} in the case where
$\calX$ is a left exact localization of a category of prestacks on
an ordinary category, and in \cite{toen} in general. We will
summarize their results here. What follows will not be needed
later in this paper and may be safely omitted by the reader.

We will need to employ a bit more terminology concerning
simplicial objects in an $\infty$-category. Let $\Delta$ denote
the category of combinatorial simplices, and let $\Delta_{\leq n}$
denote the full subcategory consisting of combinatorial simplices
of dimension $\leq n$. If $\calC$ is any $\infty$-category, then a
simplicial object of $\calC$ is a contravariant functor
$E_{\bigdot}: \Delta \rightarrow \calC$. This induces a functor
$\Delta_{\leq n} \rightarrow \calC$ by restriction, called the
{\it skeleton} of $E_{\bigdot}$. If $\calC$ admits finite limits,
then the skeleton functor $\calC^{\Delta^{op}} \rightarrow
\calC^{\Delta_{\leq n}^{op}}$ has a right adjoint, called the
$n$-coskeleton. If $E_{\bigdot}$ is a simplicial object of
$\calC$, we will let $\cosk^{n}$ denote the $n$-coskeleton of the
$n$-skeleton of $E$; this is a new simplicial object equipped with
a map $E_{\bigdot} \rightarrow \cosk^{n}_{\bigdot}(E_{\bigdot})$
which is an equivalence on the $n$-skeleton (and is universal with
respect to this property). We say that $E_{\bigdot}$ is {\it
$n$-coskeletal} if $E_{\bigdot} \rightarrow
\cosk^n_{\bigdot}(E_{\bigdot})$ is an equivalence.

\begin{definition}
Let $\calX$ be an $\infty$-topos, and let $E_{\bigdot}$ denote a
simplicial object of $\calX$. We shall say that $E_{\bigdot}$ is a
{\it hypercovering} if, for each $n \geq 0$, the natural map
$E_{n} \rightarrow \cosk^{n-1}_{n}(E_{\bigdot})$ is surjective.

More generally, if $E_{\bigdot} \rightarrow E$ is an augmented
simplicial object, we shall say that $E_{\bigdot}$ is a {\it
hypercovering of $E$} if $E_{\bigdot}$ is a hypercovering in the
$\infty$-topos $\calX_{/E}$.
\end{definition}

\begin{theorem}{\cite{toen}, \cite{hollander}}\label{hyer}
Let $\calX$ be an $\infty$-topos, let $S$ denote the class of
$\infty$-connected morphisms of $\calX$, and let $S'$ denote the
class of morphisms having the form $| E_{\bigdot} | \rightarrow
E$, where $E_{\bigdot}$ is a hypercovering of $E$ in $\calX$. An
object $E' \in \calX$ is $S$-local if and only if it is
$S'$-local.
\end{theorem}

\begin{proof}
Suppose first that $E'$ is $S'$-local. Let $f: X \rightarrow Y$ be
an $\infty$-connected morphism. Regard $X$ as a constant
simplicial object over $Y$. If we can show that $X$ is a
hypercovering, then $X \simeq |X| \rightarrow Y$ induces an
equivalence $\Hom_{\calX}(Y,E') \rightarrow \Hom_{\calX}(X,E')$
and we are done. In other words, we need to show that for each $n
\geq 0$, the natural map $X \rightarrow \cosk^{n-1}_{n}(X)$ is
surjective (here the coskeleton is formed in $\calX_{/Y}$). This
is precisely the condition that $\pi_n(f)$ vanishes.

To prove the converse, we show that if $E_{\bigdot} \rightarrow E$
is a hypercovering then $|E_{\bigdot}| \rightarrow E$ is
$\infty$-connected. Without loss of generality, we may replace
$\calX$ by $\calX_{E}$ and thereby assume that $E=1$. We need to
prove that $|E_{\bigdot}|$ is $n$-connected for every $n$. By
Corollary \ref{goober2}, it suffices to show that $\Hom_{\calX}(
1, Y) \rightarrow \Hom_{\calX}( |E_{\bigdot}|, Y)$ is an
equivalence whenever $Y$ is $n$-truncated. In this case,
$\Hom_{\calX}( |E_{\bigdot}|, Y)$ does not depend on $E_{k}$ for
$k > n+1$. Thus we may replace $E_{\bigdot}$ by its
$(n+1)$-coskeleton without loss of generality. One then shows, by
induction on $k$, that $| \cosk^n_{\bigdot}(E_{\bigdot}) |
\rightarrow E$ is an equivalence, using the fact that
$E_{\bigdot}$ is a hypercovering and the fact that all groupoids
are effective in $\calX$. Taking $k = n+1$, we obtain the desired
result.
\end{proof}

In other words, the $S$-local objects of $\calX$ are precisely the
objects which {\it satisfy hyperdescent}.

\subsection{Stacks Versus Hyperstacks}\label{versus}

If $\toposX$ is a Grothendieck topos, then we have seen how to
construct two potentially different $\infty$-topoi associated to
$\toposX$. First, we have the $\infty$-topos $\calX= \stacks
\toposX$ consisting of stacks on $X$. Second, we can form the
$\infty$-topos $\calX^{\red}$ consisting of stacks on $X$
satisfying hyperdescent, which is a localization of $\calX$. It
seems that the majority of the literature is concerned with
$\calX^{\red}$, while $\calX$ itself has received less attention
(although there is some discussion in \cite{hollander}). We would
like to make the case that $\calX$ is the more natural choice.
Proposition \ref{universality} provides a formal justification for
this claim: the functor which constructs $\calX$ from $\toposX$ is
the adjoint of the forgetful functor from $\infty$-topoi to topoi.
In this section, we would like to summarize some less formal
reasons why it may be nicer to work with $\calX$.

\begin{itemize}

\item Suppose that $\pi: X \rightarrow S$ and $\psi: S' \rightarrow S$
denote continuous maps of locally compact topological spaces, and
let $X' = X \times_{S} S'$ with projections $\pi': X' \rightarrow
S'$ and $\psi': X' \rightarrow X$. In classical sheaf theory, one
has a natural transformation $\psi^{\ast} \pi_{\ast} \rightarrow
\pi'_{\ast} {\psi'}^{\ast}$ of functors between the derived
categories of left-bounded complexes of sheaves on $X$ and on
$S'$. The proper base change theorem asserts that this
transformation is an equivalence whenever the map $\pi$ is proper.

The functors $\psi^{\ast} \pi_{\ast}$ and $\pi'_{\ast}
{\psi'}^{\ast}$ may also be defined on $\infty$-categories of
stacks and $\infty$-categories of hyperstacks, and one again has a
base-change transformation as above. It is natural to ask if the
base change transformation is an equivalence when $\pi$ is proper.
It turns out that this is {\em true} for the $\infty$-categories
of stacks, but {\em false} for the $\infty$-categories of
hyperstacks:

\begin{counterexample}\label{thrust}
Let $Q$ denote the Hilbert cube $[0,1] \times [0,1] \times
\ldots$, and let $\pi: X \rightarrow [0,1]$ be the projection onto
the first factor. For each $i$, we let $Q_i$ denote ``all but the
first $i$'' factors of $Q$, so that $Q = [0,1]^i \times Q_i$.

We construct a stack $\calF$ on $X = Q \times [0,1]$ as follows.
Begin with the empty stack. Adjoin to it two sections, defined
over the open sets $[0,1) \times Q_1 \times [0,1)$ and $(0,1]
\times Q_1 \times [0,1)$. These sections both restrict to give
sections of $\calF$ over the open set $(0,1) \times Q_1 \times
[0,1)$. We next adjoin paths between these sections, defined over
the smaller open sets $(0,1) \times [0,1) \times Q_2 \times
[0,\frac{1}{2})$ and $(0,1) \times (0,1] \times Q_2 \times [0,
\frac{1}{2})$. These paths are both defined on the smaller open
set $(0,1) \times (0,1) \times Q_2 \times [0, \frac{1}{2})$, so we
next adjoin two homotopies between these paths over the open sets
$(0,1) \times (0,1) \times [0,1) \times Q_3 \times [0,
\frac{1}{3})$ and $(0,1) \times (0,1) \times (0,1] \times Q_3
\times [0, \frac{1}{3})$. Continuing in this way, we obtain a
stack $\calF$. On the closed subset $Q \times \{0\} \subset X$,
the stack $\calF$ is $\infty$-connected by construction, and
therefore the associate hyperstack admits a global section.
However, the hyperstack associated to $\calF$ does not admit a
global section in any neighborhood of $Q \times \{0\}$, since such
a neighborhood must contain $Q \times [ 0, \frac{1}{n})$ for $n
\gg 0$ and the higher homotopies required for the construction of
a section are eventually not globally defined.
\end{counterexample}

The same issue arises in classical sheaf theory if one wishes to
work with unbounded complexes. In \cite{spaltenstein},
Spaltenstein defines a derived category of unbounded complexes of
sheaves on $X$, where $X$ is a topological space. His definition
forces all quasi-isomorphisms to become invertible, which is
analogous to procedure of obtaining $\calX^{\red}$ from $\calX$ by
inverting the $\infty$-connected morphisms. Spaltenstein's work
shows that one can extend the {\it definitions} of all of the
basic objects and functors. However, it turns out that the {\it
theorems} do not all extend: in particular, one does not have the
proper base change theorem in Spaltenstein's setting
(Counterexample \ref{thrust} can easily be adapted to the setting
of complexes of abelian sheaves). The problem may be rectified by
imposing weaker descent conditions, which do not invert all
quasi-isomorphisms.

The proof of the proper base change theorem in the context of
$\infty$-topoi will be given in a sequel to this paper.

\item The $\infty$-topos $\calX$ has better finiteness properties
than $\calX^{\red}$. Suppose that $X$ is a coherent topological
space (that is, $X$ has a basis of compact open sets which is
stable under the formation of finite intersections), and let
$\calX$ denote the $\infty$-category of stacks on $X$ (in other
words, stacks on the ordinary topos of sheaves on $X$). We may
construct the $\infty$-category $\calX$ as a localization of the
$\infty$-category $\calP$ of prestacks on $\calC$, where $\calC$
denotes the partially ordered set of compact open subsets of $X$.
The localization functor $L$ is given by transfinitely iterating
the construction $\calF \mapsto \calF^{+}$ of Proposition
\ref{plus}. Since every covering in $X_{\omega}$ may be refined to
a finite covering, one can show that $\calF \mapsto \calF^{+}$
commutes with filtered colimits in $\calF$. It follows easily that
$L$ commutes with filtered colimits. From this, we may deduce that
$L$ carries compact objects of $\calP$ into compact objects of
$\calX$. Since $\calP$ is generated by compact objects, we deduce
that $\calX$ is generated by compact objects.

By contrast, one can give an example of a coherent topological
space for which the ``hypersheafification'' of prestacks does not
commute with filtered colimits. The construction is modeled on
Counterexample \ref{thrust}, replacing each copy of the interval
$[0,1]$ with a finite topological space having a generic point
(analogous to the interior $(0,1)$) and two special points
(analogous the the endpoints $0$ and $1$). The compatibility of
sheafification with filtered colimits is a crucial feature of
coherent topoi, which might be lost in the $\infty$-categorical
setting if $\calX^{\red}$ is used.

\begin{remark}
In particular, we note that $\calX \neq \calX^{\red}$ when
$\toposX$ is a coherent topos, so that $\calX$ does not
necessarily have enough points. Thus, one does not expect an
$\infty$-categorical version of the G\"{o}del-Deligne theorem to
hold.
\end{remark}

\item Although $\calX$ and $\calX^{\red}$ are not equivalent in
general, they do coincide whenever certain finite-dimensionality
conditions are satisfied (see Corollary \ref{fdfd}). These
conditions are satisfied in most of the situations in which
Jardine's model structure is commonly used, such as the Nisnevich
topology on a scheme of finite Krull dimension.

\item Let $X$ be a paracompact space, $\calX$ the $\infty$-topos
of stacks on $X$, $K$ a homotopy type, and $p: \calX \rightarrow
\SSet$ the natural geometric morphism. Then $\pi_0 p_{\ast}
p^{\ast} K$ is a natural definition of the ``sheaf cohomology'' of
$X$ with coefficients in $K$, and this agrees with the definition
$[X,K]$ given by classical homotopy theory (Theorem \ref{nice}).
As we mentioned in the introduction, this fails if we replace
$\calX$ by $\calX^{\red}$.

\item If the topos $\toposX$ has enough points, then the local
equivalences of $\calX$ are precisely the morphisms which induce
weak equivalences on each stalk. It follows immediately that
$\calX^{\red}$ has enough points, and that if $\calX$ has enough
points then $\calX \simeq \calX^{\red}$. The possible failure of
the Whitehead theorem in $\calX$ may be viewed either as a bug or
a feature. While this failure has the annoying consequence that
$\calX$ may not have ``enough points'' even though $\toposX$ does,
it also means that $\calX$ might detect certain global phenomena
which cannot be properly understood by restricting to points. Let
us consider an example from classical geometric topology. If $f: X
\rightarrow Y$ is a continuous map between compact ANRs, then $f$
is said to be {\it cell-like} if for any $y \in Y$, the inverse
image $f^{-1}(y)$ is contractible in arbitrarily small
neighborhoods of itself. In the finite dimensional case, this
holds if and only if each fiber $f^{-1}(y)$ has trivial shape, but
this fails in general.

Let us indicate how the notion of an $\infty$-topos can shed some
light on the situation. First, we note the following:

\begin{proposition}\label{cell}
Let $f: X \rightarrow Y$ be a map of compact ANRs, and let $f$
also denote the associated geometric morphism from $\calX$ to
$\calY$, where $\calX$ and $\calY$ denote the $\infty$-topoi of
stacks on $X$ and $Y$, respectively. The following conditions are
equivalent:

\begin{enumerate}
\item The map $f$ is cell-like.

\item For every $E \in \calY$, the adjunction morphism $E
\rightarrow f_{\ast} f^{\ast} E$ is an equivalence.

\item The functor $f^{\ast}$ is fully faithful.
\end{enumerate}

\end{proposition}

\begin{proof}
A proof will be given in a sequel to this paper.
\end{proof}

Since $f$ is proper, the stalk of $f_{\ast} f^{\ast} E$ at a point
$y$ is equal to the cohomology of $f^{-1}(y)$ with coefficients in
$E_y$ (by the proper base change theorem promised above), which is
equivalent to $E_y$ if $f^{-1}(y)$ has trivial shape (essentially
by the definition of the strong shape). It follows that if all of
the fibers of $f$ have trivial shape, then the adjunction map $E
\rightarrow f_{\ast} f^{\ast} E$ induces an equivalence on all
stalks, and is therefore $\infty$-connected, for any $E \in
\calY$. If $Y$ is finite dimensional, we can conclude that $f$ is
cell-like, but not in general.

It is crucial that we work with stacks, rather than hyperstacks,
in the above discussion. If we chose instead to work with
hyperstacks, then we could take Proposition \ref{cell} as the
definition of a concept analogous to, but distinct from, that of a
cell-like map. However, this analogue is of more limited geometric
interest.
\end{itemize}

\section{Paracompact Spaces}\label{paracompactness}

Let $X$ be a topological space, and let $p: X \rightarrow \ast$ be
the projection. Let $K$ be a homotopy type. Then we may regard $K$
as a stack on $\ast$, and form the stack $p_{\ast} p^{\ast} K$.
Alternatively, we can form an actual topological space $|K|$
(which we assume to be nice, say a metric ANR) and study the
collection of homotopy classes of maps $[X, |K|]$. The goal of
this section is to prove the following:

\begin{theorem}\label{nice}
If $X$ is paracompact, then there is a bijection
$$\phi: [X,|K|] \rightarrow \pi_0 p_{\ast} p^{\ast} K.$$
\end{theorem}

In fact, the map $\phi$ always exists and is natural in both $K$
and $X$; only the bijectivity of $\phi$ requires $X$ to be
paracompact.

Let us outline the proof of Theorem \ref{nice}. In order to do
calculations with stacks on $X$, we will need to introduce some
models for the stacks. It is most common to use simplicial
presheaves on $X$, but this is inconvenient for two reasons. The
first is that general open subsets of $X$ can behave badly; we
will therefore restrict our attention to presheaves defined with
respect to a nice basis for the topology of $X$. The second
problem is more subtle, but for technical reasons it will be
convenient for us to introduce a category $\calK$ of
$\Ind$-compact metric spaces for use as a model of homotopy
theory, rather than simplicial sets. We do not require much from
the category $\calK$: it does not need to be a Quillen model
category, for example. The main properties of $\calK$ which we
shall need are as follows:

\begin{itemize}
\item There is a good formal relationship between $\calK$-valued
presheaves on $X$ and topological spaces {\em over} $X$.
\item Any $\calK$-valued presheaf on $X$ gives rise to a prestack
on $X$, and all prestacks arise in this way.
\item The homotopy limit constructions needed to ``stackify'' a
prestack may be carried out easily on the level of $\calK$-valued
presheaves.
\end{itemize}

After establishing these properties, we will be in a position to
prove Theorem \ref{main}, which is closely related to Theorem
\ref{nice} and has a few other interesting implications. Finally,
in \S \ref{dooky} we will show how to deduce Theorem \ref{nice}
from Theorem \ref{main}.

\subsection{Some Point-Set Topology}

Let $X$ denote a paracompact topological space. In order to prove
Theorem \ref{main}, we will need to investigate the homotopy
theory of prestacks on $X$. We then encounter the following
technical obstacle: an open subset of a paracompact space need not
be paracompact. Because we wish to deal only with paracompact
spaces, it will be convenient to restrict our attention to
presheaves which are defined only with respect to a particular
basis $\calB$ for $X$ consisting of paracompact open sets. More
precisely, we need the following:

\begin{lemma}\label{goofy}
Let $X$ be a paracompact topological space. There exists a
collection $\calB$ of open subsets of $X$ with the following
properties:

\begin{itemize}
\item The elements of $\calB$ form a basis for the topology of
$X$.
\item Each element of $\calB$ is paracompact.
\item The elements of $\calB$ are stable under finite
intersections.
\end{itemize}
\end{lemma}

\begin{proof}
If $X$ is a metric space, then we may take $\calB$ to consist of
the collection of {\em all} open subsets of $X$. In general, this
does not work because the property of paracompactness is not
inherited by open subsets. A more general solution is to take
$\calU$ to consist of all open $F_{\sigma}$ subsets of $X$,
together with $X$ itself. Recall that a subset of $X$ is called an
$F_{\sigma}$ if it is a countable union of closed subsets of $X$.
It is then clear that $\calB$ is stable under finite
intersections. The paracompactness of each element of $\calB$
follows since each is a countable union of paracompact spaces (see
for example \cite{munkres}). To see that $\calB$ forms a basis for
$X$, we argue as follows. Given any $x \in V \subseteq X$, we may
construct a sequence of closed subsets $K_i \subseteq X$ with $K_0
= \{x\}$ and such that $K_{i+1} \subseteq V$ contains some
neighborhood of $K_i$ (by normality). Then $U = \bigcup K_i$ is an
open $F_{\sigma}$ containing $x$ and contained in $V$.
\end{proof}

In the sequel, we will frequently suppose some collection of open
sets $\calB$ has been chosen so as to have the above properties.
Then $\calB$ may be viewed as a category with finite limits, and
is equipped with a natural Grothendieck topology. The topos of
sheaves of sets on $X$ is equivalent to the topos of sheaves of
sets on $\calB$. Thus, the proof of Proposition \ref{plus} shows
that the $\infty$-topos of $X^{\infty}$ is equivalent to the
$\infty$-topos of stacks on $\calB$ (in the sense of Proposition
\ref{plus}).

\subsection{Presheaves of Spaces}

Throughout this section, we will assume that $X$ is a topological
space and that $\calB$ is a collection of open subsets of $X$
which is stable under finite intersections. We will abuse notation
by identifying $X$ with the topos of sheaves on $X$.

In order to prove Theorem \ref{nice}, we will need to do some
calculations in the $\infty$-topos $\stacks X$. For this, we will
need to use some concrete models for objects in $\stacks X$, and
to exhibit these we need to choose a model for homotopy theory.
The homotopy theory of ``sheaves of spaces'' is usually formalized
in terms of simplicial presheaves, as in \cite{jardine}. However,
for our purposes it will be more convenient to use a somewhat
obscure model for homotopy types.

Let $\calK_0$ denote the (ordinary) category of compact metrizable
spaces (and continuous maps). We let $\calK$ denote the category
$\Ind(\calK_0)$. Passing to filtered colimits gives a functor $F$
from $\calK$ to topological spaces. Inside of $\calK_0$ one may
find the ordinary category of finite CW-complexes, and inside of
$\calK$ one may find the ordinary category of {\em all}
CW-complexes. On these subcategories the functor $F$ is fully
faithful. In general, the functor $F$ is neither full nor
faithful, and its essential image contains many ``bad'' spaces.
However, for our purposes it will be convenient to use the
category $\calK$ as a model for homotopy theory.

\begin{remark}
Since the category $\calK_0$ is essentially small and has finite
colimits, the category $\calK$ is presentable, and may therefore
be identified with the category of colimit-preserving presheaves
on $\calK$, which are in turn determined by their restriction to
$\calK_0$. In other words, we may view an object $Z \in \calK$ as
being given by a contravariant functor $\Hom( \bigdot, Z)$ on the
category of compact metrizable spaces. This functor is required to
be compatible with finite colimits in $\calK_0$, and every such
functor is representable by an object of $\calK$.
\end{remark}

We let $\calK_{\calB}$ denote the (ordinary) category of
presheaves on $\calB$ with values in $\calK$. Since $\calK_0$
contains the category of finite simplices, every object of $\calK$
determines a (fibrant) simplicial set and we obtain a functor
$\Psi: \calK \rightarrow \SSet$, and therefore a functor
$\Psi_{\calB}: \calK_{\calB} \rightarrow \SSet^{\calB^{op}}$.
Proposition \ref{buildmodel}, which we shall prove later, implies
that $\Psi_{\calB}$ is essentially surjective. We wish to use the
source category $\calK_{\calB}$ as ``models'' for the target
$\infty$-category $\SSet^{\calB^{op}}$. We use this terminology
only to suggest intuition; we neither need nor use the language of
model categories.

We will let $\Sp_X$ denote the (ordinary) category of {\it
Hausdorff topological spaces over $X$}; that is, an object of
$\Sp_X$ is a Hausdorff space $Y$ equipped with a continuous map
$f: Y \rightarrow X$. The proof of Theorem \ref{nice} will make
use of a pair of adjoint functors
$$S: \Sp_X \rightarrow \calK_{\calB}$$
$$|\,|: \Sp_X \leftarrow \calK_{\calB}$$
which we now describe.

We will begin with the functor $S$, whose definition is more
intuitive. Let $f: Y \rightarrow X$ be an object of $\Sp_X$. We
define an object $S(Y) \in \calK_{\calB}$ by the following
formula:
$$ \Hom_{\calK}( Z, S(Y)(U) ) = \Hom_{X}( Z \times U, Y)$$
where $Z$ is any compact metrizable space, $U$ is any element of
$\calB$, and the right hand side denotes the space of continuous
maps {\em over} the space $X$. Using the fact that $Y$ is
Hausdorff, one readily checks that this formula is compatible with
colimits in $Z$, and therefore uniquely defines $S(Y)(U) \in
\calK$ (which is automatically functorial in $U$). We note that
$S$ is a functor from $\Sp_Y$ to $\calK_{\calB}$.

The functor $S$ has a left adjoint, which we will denote by
$$ \calF \mapsto | \calF |.$$ We shall call $|\calF|$ the {\it
geometric realization} of $\calF$. The construction of $|\calF|$
is straightforward: one begins by gluing together the spaces
$\calF(U) \times U$ in the obvious way, and then takes a maximal
Hausdorff quotient. We leave the details to the reader.

\subsection{Homotopy Limits}

Let $X$ be a topological space, $\calB$ a basis for the topology
of $X$ which is closed under finite intersections.

Suppose that $\calF \in \calK_{\calB}$ is a $\calK$-valued
presheaf on the basis $\calB$. Then $\Psi_{\calB} \calF$ is a
prestack on $\calB$. In the process of sheafifying $\Psi_{\calB}
\calF$, we encounter the limit
$$\lim_{V \in \calS} \Psi_{\calB}( \calF )(V) $$
where $\calS$ is a sieve on an open set $U \in \calB$. This limit
is canonically determined in the $\infty$-category $\SSet$ of
spaces. However, for concrete computations, we would like to
represent this limit by an element of $\calK$ which depends
functorially (in the usual, $1$-categorical sense) on $\calF$. For
this, we shall fix an open cover $\calU = \{ U_{\alpha} \}_{\alpha
\in A}$ of $U$ which generates the sieve $\calS$.

Now define an object $\holim_{\calU} \calF \in \calK$ by the
following universal property:

\begin{itemize}
\item For any compact metrizable space $K$,
$\Hom(K, \holim_{\calU} \calF )$ is equal to the set of all
collections of compatible maps $\phi_v: K \times \Delta^n
\rightarrow \calF_{\bigdot}( U_{v(0)} \cap U_{v(1)} \cap \ldots
\cap U_{v(n)} )$, where $v$ ranges over all functions $[0,n]
\rightarrow A$. The compatibility requirement is that $\phi_v
\circ (\id_K \times \Delta^{f}) = \phi_{v \circ f}$, for any
factorization $[0,m] \stackrel{f}{\rightarrow} [0,n]
\stackrel{v}{\rightarrow} \calU$.
\end{itemize}

\begin{remark}
The existence of $\holim_{\calU} \calF$ depends on the
representability of the functor described above. It suffices to
show that the functor carries finite (colimits) in $K$ to finite
limits, which is obvious.
\end{remark}

If $\calF \in \calK_{\calB}$, then $\holim_{\calU} \calF$ actually
represents the limit
$$\lim_{V \in \calS} \Psi \calF(V) \in \SSet.$$ In
particular, it is independent of the choice of $\calU$, up to
canonical homotopy equivalence. This is straightforward and left
to the reader.

If $\calU$ is an open cover of $U$ generating a sieve $\calS$,
then the functor
$$ \calF \mapsto \holim_{\calU} \calF$$ is
co-represented by an object $N(\calU) \in \calK_{\calB}$. More
precisely, let $N(\calU)(V)$ denote the object of $\calK_{\calB}$
whose $n$-simplices are indexed by $(n+1)$-element subsets of $\{
\alpha \in A: V \subseteq U_{\alpha} \}$. Then we have
$$\Hom_{\calK}(K, \holim_{\calU} \calF) = \Hom_{\calK_{\calB}}(K
\times N(\calU), \calF)$$ where on the right side, we regard $K$
as a constant presheaf in $\calK_{\calB}$.

Now we have the following easy result, which explains the
usefulness of a basis of {\em paracompact} open sets.

\begin{lemma}\label{partit}
Suppose that $U \in \calB$ is paracompact, and let $\calS$ be a
sieve on $U$, generated by an open cover $\calU = \{ U_{\alpha}
\}_{\alpha \in A}$. The natural map $\pi: | N(\calU) | \rightarrow
U$ is a fiberwise homotopy equivalence of spaces over $X$. In
other words, there exists a map $s: U \rightarrow | N(\calU) |$
such that $\pi \circ s$ is the identity and $s \circ \pi$ is
fiberwise homotopic to the identity.
\end{lemma}

\begin{proof}
Any partition of unity subordinate to the cover $\calU$ gives rise
to a section $s$. (Indeed, a section $s$ is almost the same thing
as a partition of unity subordinate to $\calU$; there is only a
slight difference in the local finiteness conditions.) To show
that $s \circ \pi$ is fiberwise homotopic to the identity, use a
``straight line'' homotopy.
\end{proof}

\begin{proposition}\label{aese}
Suppose that $\calB$ is a basis for a topological space $X$ which
is closed under finite intersections and consists of paracompact
sets. Then for any space $Y \in \Sp_X$, the prestack $\Psi_{\calB}
S(Y)$ on $\calB$ is a stack.
\end{proposition}

\begin{proof}
Let $\calF = S(Y)$. Choose any $U \in \calB$ and any sieve $\calS$
which covers $U$; we must show that
$$f: \calF(U) \rightarrow \lim_{\calU} \calF$$
is a weak homotopy equivalence (that is, it induces an equivalence
after applying $\Psi$). It will suffice to prove the following
statement for any inclusion of finite simplicial complexes $K_0
\subseteq K$:

\begin{itemize}
\item Given any map $g: K_0 \rightarrow \calF(U)$, and any
extension $h$ of $f \circ g$ to $K$, there exists an extension

$\widetilde{g}$ of $g$ to $K_0$ such that $f \circ \widetilde{g}$
is homotopic to $h$ by a homotopy fixed on $K_0$.
\end{itemize}

Since $K_0$ and $K$ are compact metric spaces, the above assertion
is equivalent to the following:

\begin{itemize}
\item Given any map $g: K_0 \times U \rightarrow Y$ which is a
section of $Y \rightarrow X$ and any extension $h$ of the induced
map $K_0 \times N(\calU) \rightarrow Y$ to $K \times N(\calU)
\rightarrow Y$, there exists an extension $\widetilde{g}: K \times
U \rightarrow Y$ such that the induced map $K \times N(\calU)
\rightarrow Y$ is homotopic to $h$ by a homotopy fixed on $K_0
\times N(\calU)$.
\end{itemize}

This follows immediately from Lemma \ref{partit} and the fact that
$K_0 \rightarrow K$ is a cofibration.
\end{proof}

\begin{remark}
The prestack $\Psi_{\calB} S(Y)$ is {\em not necessarily a
hyperstack}. This observation was one of the main motivations for
this paper.
\end{remark}

\subsection{Simplicial Complexes}

In this section, we discuss the use of simplicial complexes as a
model for homotopy theory. Recall that a {\it combinatorial
simplicial complex} consists of the following data:

\begin{itemize}
\item A set $V$ of vertices.
\item A collection $J(V)$ of finite subsets of $V$, containing the
empty set, and with the property that $S \in J(V)$ and $S'
\subseteq S$ implies $S' \in J(V)$
\end{itemize}

A morphism between combinatorial simplicial complexes $(V,J(V))$
and $(V', J'(V'))$ is a map $\alpha: V \rightarrow V'$ such that
$\alpha(S) \in J'(V')$ for any $S \in J(V)$. Let $\Simp$ denote
the category of combinatorial simplicial complexes.

Given a combinatorial simplicial complex $K=(V,J(V))$, we can
associate an object $rK \in \calK$, the {\it geometric
realization} of $K$. If $V$ is finite, then we take $rK$ to be the
full subcomplex of the simplex $\Delta^V$ spanned by the vertices
in $V$, consisting of all faces $F \subseteq \Delta^V$ such that
$\{ v \in V : v \in F\} \in J(V)$. In the general case, we define
$rK$ to be a formal filtered colimit of the spaces $rK_{\alpha}$,
where $K_{\alpha} = (V_{\alpha}, \{ S \in J(V): S \subseteq
V_{\alpha} \})$ and $V_{\alpha}$ ranges over the finite subsets of
$V$.

In particular, after composition with $\Psi$, we obtain a
geometric realization functor $\Simp \rightarrow \SSet$. It is
well-known that this functor is essentially surjective (in other
words, every topological space has the weak homotopy type of a
simplicial complex). In this section, we sketch a proof of a
relative version of this statement.

Let $\calB$ be a partially ordered set (the particular case we
have in mind is that $\calB$ is a basis of paracompact open sets
for a paracompact space $X$, but that will play no role in this
section). We let $\Simp_{\calB}$ denote the category of
$\Simp$-valued presheaves on $\calB$ for which the induced
presheaf of vertex sets is constant. In other words, an object of
$\Simp_{\calB}$ consists of a vertex set $V$ together with a
presheaf of families of finite subsets of $V$, which satisfy the
definition of a combinatorial simplicial complex over each $U \in
\calB$.

Applying the realization functor $R$ pointwise, we obtain a
functor $\Simp_{\calB} \rightarrow \calK_{\calB}$, which we shall
denote also by $r$. The main result we will need is the following:

\begin{proposition}\label{buildmodel}
The functor $\Psi r: \Simp_{\calB} \rightarrow \SSet^{\calB^{op}}$
is essentially surjective.
\end{proposition}

\begin{proof}
Let $\calF$ be a stack on $\calB$. If we model objects of $\SSet$
by simplicial sets, then we may view $\calF$ as a homotopy
coherent diagram of simplicial sets, indexed by $\calB$. Using
standard techniques, we may replace this homotopy coherent diagram
by a strictly commutative diagram, which we shall also denote by
$\calF$.

Not every simplicial set $S_{\bigdot}$ corresponds to a
combinatorial simplicial complex: roughly speaking, this requires
that simplices of $S_{\bigdot}$ are determined by their vertices.
However, this condition is always satisfied if we replace
$S_{\bigdot}$ by its first barycentric subdivision, after which we
may functorially extract a combinatorial simplicial complex whose
geometric realization is homeomorphic to the geometric realization
of $S_{\bigdot}$. Applying this procedure to the functor $\calF$,
we may replace $\calF$ with a strictly commutative diagram in the
category $\Simp$, which we shall denote by $\calF'$.

Let $\calF'(U) = (V_U, J_U(V_U))$. The problem now is that the
sets $V_U$ may vary; we need to replace $\calF'$ by a diagram in
$\Simp$ which has a constant vertex set. This may be achieved as
follows: let $V = \coprod_{U \in \calB} V_U$. For each $S \in
J_U(V_U)$, we let
$$S' = \{ v \in V_{U'}: (U \subseteq U') \wedge (v|U \in S) \},$$
and we let $J'_U(V)$ denote the collection of all finite subsets
of $V$ which are contained in $S'$ for some $S \in J_U(V_U)$. Then
the collection $\{ (V, J'_U(V))\}_{U \in \calB}$ gives an object
of $\Simp_{\calB}$, regarded as a $\calB$-indexed diagram $\calF'$
in $\Simp$. By construction, there is a natural transformation
$\calF \rightarrow \calF'$. One can easily check that it induces
homotopy equivalences at each level $U \in \calB$.
\end{proof}

\subsection{The Main Result}

The main goal of this section is to prove:

\begin{lemma}\label{hardd}
Let $X$ be a paracompact topological space, and $\calB$ a basis
for the topology of $X$ which is stable under finite intersections
and consists of paracompact sets. Let $K \in \Simp_{\calB}$. Then
the induced map
$$ (\Psi_{\calB} rK)^{+} \rightarrow (\Psi_{\calB} S|rK|)^{+}$$
is an equivalence of presheaves of spaces on $\calB$. Here $r:
\Simp_{\calB} \rightarrow \calK_{\calB}$ denotes the functor
obtained by applying $r: \Simp \rightarrow \calK$ componentwise,
and $\calF \rightsquigarrow \calF^{+}$ denotes the partial
sheafification functor of Proposition \ref{plus}.
\end{lemma}

\begin{proof}
We must show that the induced map $p: (\Psi rK)^{+}(U) \rightarrow
(\Psi S|rK|)^{+}(U)$ is an equivalence for each $U \in \calB$.
Replacing $X$ by $U$ and all other objects by their restrictions
to $U$, we may reduce to the case where $U = X$.

We prove the following: given any map $f: S^n \rightarrow (\Psi
rK)^{+}(X)$, and any extension $g$ of $p \circ f$ to $D^{n+1}$,
there exists an extension $\widetilde{f}$ of $f$ to $D^{n+1}$ such
that $p \circ \widetilde{f}$ is homotopic to $g$ by a homotopy
fixed on $S^{n}$. Since $(\Psi rK)^{+}(X)$ and $(\Psi S|rK|)^{+}$
are defined as filtered colimits, and $S^n$ and $D^{n+1}$ are
finite complexes, we may assume given factorizations $S^n
\rightarrow \lim_{V \in \calS} (\Psi rK)(V)$ and $D^{n+1}
\rightarrow \lim_{V \in \calS} (\Psi S|rK|)(V)$. Now choose an
open cover $\calU$ of $X$ which generates the sieve $\calS$; we
may then represent $f$ and $g$ by maps $S^n \times N(\calU)
\rightarrow RK$ and $D^{n+1} \times N(\calU) \rightarrow S|rK|$,
which we shall again denote by $f$ and $g$.

We will produce a map $\widetilde{f}: D^{n+1} \times N(\calU)
\rightarrow rK$, possibly after passing to some refinement of the
original covering $\calU$, which extends $f$ and induces a map
$D^{n+1} \times N(\calU) \rightarrow S|rK|$ which is homotopic to
$g$ by a homotopy fixed on $S^n \times N(\calU)$. Let $\calU =
\{U_{\alpha} \}_{\alpha in A}$. We first choose a covering $\calU'
= \{U'_{\alpha}\}_{\alpha \in A}$ such that $\calU'$ is locally
finite and $\overline{U'_{\alpha}} \subseteq U_{\alpha}$ for each
$\alpha \in A$.

Let $K = (V, J(V))$. Then we may regard $g$ as determining a
continuous map $D^{n+1} \times |N(\calU)| \rightarrow |rK|
\subseteq \Delta^V \times X$, where $\Delta^V$ denotes the
(possibly infinite) simplex spanned by the vertices in $V$. For
each $v \in V$, let $\Delta^V_v$ denote the open subset of
$\Delta^V$ consisting of points with nonzero $v$-coordinate, so
that the open sets $\Delta^V_v$ cover $\Delta^V$.

Choose any $\alpha \in A$ and any $x \in U'_{\alpha}$. Let $A_x =
\{ \alpha' \in A: x \in \overline{U'_{\alpha'}} \}$, so that $g$
determines a map $D^{n+1} \times \Delta^{A_x} \times W_{x,\alpha}
\rightarrow \Delta^V$ for some small neighborhood $W_{x,\alpha}
\subseteq U'_{\alpha}$ of $x$. Note that $A_x$ is a finite set.
Using the compactness of $D^{n+1} \times \Delta^{A_x}$ and a
tube-lemma argument, we can find (after possibly shrinking $W_x$)
a finite set $V_x \subseteq V$ and open subsets $Z_v \subseteq
D^{n+1} \times \Delta^{A_x}$ for $v \in V_x$ such that $g( Z_v
\times W_{x,\alpha} ) \subseteq \Delta^V_v$. Choose a partition of
unity dominated by the cover $\{ Z_v: v \in V_x \}$, and use this
partition of unity to define a map $g_{x,\alpha}: D^{n+1} \times
\Delta^{A_x} \times W_x \rightarrow \Delta^{V_x}$.

Now let $B = \{ (x,\alpha): (\alpha \in A) \wedge x \in
U'_{\alpha}$, and consider the open cover $\calW = \{ W_{x,\alpha}
\}_{(x,\alpha) \in B}$. This cover refines $\calU$ via the
projection $\pi: B \rightarrow A$. We define a map
$\widetilde{f}': D^{n+1} \times N(\calW) \rightarrow RK$ as
follows. Given $\beta_0, \ldots, \beta_m \in B$, we let
$\widetilde{f}'( D^{n+1} \times \Delta^{ \{\beta_0, \ldots,
\beta_n\} } \rightarrow RK( W_{\beta_0} \cap \ldots \cap
W_{\beta_m} )$ by the formula $$\widetilde{f}(z, \sum_{0 \leq i
\leq m} c_i \beta_i) = \sum_{0 \leq i \leq m} c_i g_{\beta_i}(x,
\sum_{0 \leq j \leq m} c_j \pi(\beta_j)).$$ One readily checks
that $\widetilde{f}|S^{n} \times N(\calW)$ is homotopic to $f$
(restricted to $N(\calW)$) and that $p \circ \widetilde{f}'$ is
homotopic to $g$ (restricted to $N(\calW)$) via straight-line
homotopies (which are obviously compatible with one another). A
standard argument, using the fact that $S^n \rightarrow D^{n+1}$
is a cofibration, allows us to replace $\widetilde{f}'$ by a
function $\widetilde{f}$ with the desired properties.
\end{proof}

Now, the hard work is done and we can merely collect up the
consequences.

\begin{theorem}\label{main}
Let $X$ be a paracompact topological space, and $\calB$ a basis
for the topology of $X$ which consists of paracompact open sets
and is stable under intersections. Then:
\begin{itemize}
\item For any prestack $\calF$ on $\calB$, $\calF^{+}$
is a stack.
\item If $\calF$ is a prestack on $\calB$ which is
represented by an object $K \in \Simp_{\calB}$, then $S|rK|$
represents the sheafification of $\calF$.
\end{itemize}
\end{theorem}

\begin{proof}
The claim follows immediately from Lemmas \ref{aese} and
\ref{hardd}.
\end{proof}

\begin{remark}
The first part of the theorem asserts that on paracompact
topological spaces, sheafification requires only one step. This
generalizes the well-known fact that \Cech resolutions can be used
to compute cohomology on paracompact spaces.
\end{remark}

\subsection{The Proof of Theorem \ref{nice}}\label{dooky}

Now that we have Theorem \ref{main} in hand, it is easy to give a
proof of Theorem \ref{nice}.

Let $K$ be a combinatorial simplicial complex, $X$ a paracompact
topological space, and $\calX$ the $\infty$-topos of stacks on
$X$. Choose a basis $\calB$ for $X$ consisting of paracompact open
sets and stable under finite intersections, and let $\widetilde{K}
\in \calK_{\calB}$ denote the corresponding presheaf of spaces on
$\calB$. Let $p: X \rightarrow \ast$ be the canonical map, which
induces a functor $p^{\ast}: \SSet \rightarrow \calX$. By Theorem
\ref{main}, the set $\pi_0 \Hom_{\calX}(1, p^{\ast} RK)$ may be
identified with the set of homotopy classes of sections of the
projection $f: | \widetilde{K} | \rightarrow X$. To complete the
proof, it will suffice to show the following:

\begin{proposition}
There is a natural bijection between the set of homotopy classes
of sections of $f$ and homotopy classes of maps $X \rightarrow
|K|$, where $|K|$ denotes the geometric realization of $K$.
\end{proposition}

\begin{proof}
There is a continuous bijection $b: | \widetilde{K} | \rightarrow
X \times |K|$, so that any section of $f$ gives a map $X
\rightarrow |K|$. The only problem is that $b$ is not necessarily
a homeomorphism, so it is not clear that every continuous map $X
\rightarrow |K|$ arises in this way. However, we will show that
this is always true ``up to homotopy''. In order to see that the
corresponding statement is also true for homotopies {\em between}
continuous maps $X \rightarrow |K|$, it is useful to formulate a
slightly more general assertion.

If $Z$ is a compact metrizable space, let us call a continuous map
$X \times Z \rightarrow |K|$ {\it completely continuous} if it
lifts to a continuous map $X \times Z \rightarrow |
\widetilde{K}|$ (over $X$). We will prove the following claim:
\begin{itemize}
\item If $Z_0 \subseteq Z$ is a cofibration of compact metrizable spaces, and
$g: Z \times X \rightarrow |K|$ is such that $g| Z_0 \times X$ is
completely continuous, then $g$ is homotopic to a completely
continuous map by a homotopy which is fixed on $Z_0$.
\end{itemize}

Assuming the claim for the moment, let us complete the proof.
Applying the claim in to the case where $Z$ is a point, we deduce
that every continuous map $X \rightarrow |K|$ is homotopic to one
induced by a section of $f$. The injectivity assertion is proved
by applying the claim in the case where $Z = [0,1]$, $Z_0 =
\{0,1\}$.

It remains to establish the claim. Since the inclusion of $Z_0$
into $Z$ is a cofibration, a standard argument shows that it will
suffice to produce {\em any} completely continuous map homotopic
to $g$; one can then correct the continuous map and the homotopy
so that it is fixed on $Z_0 \times X$. In other words, we may
assume that $Z_0$ is empty.

Let $K = (V, J(V))$. Then we may regard $|K|$ as a subset of the
infinite simplex spanned by the elements of $V$; a map $g: Z
\times X \rightarrow |K|$ is determined by component maps $g_v: Z
\times X \rightarrow [0,1]$ which satisfy the following
conditions:
\begin{enumerate}
\item Each $g_v$ is continuous.
\item For any point $(z,x) \in Z \times X$, the set $\{v \in V: g_v(z,x) \neq 0\}$
lies in $J(V)$. In particular, it is finite.
\item For each $(z,x) \in Z \times X$, the sum $\sum_{v} g_v(z,v)$ is equal to $1$.
\end{enumerate}

Every collection of functions $\{g_v\}$ satisfying the above
condition gives a map of {\em sets} $g: Z \times X \rightarrow
|K|$, but the function $g$ is not necessarily continuous: this is
equivalent to a complicated local finiteness condition which is
slightly stronger than condition $(2)$ above. The complete
continuity of $g$ is equivalent to an even stronger local
finiteness condition. However, both of these local finiteness
conditions are satisfied if the collection $\{g_v\}$ satisfies the
following condition:

\begin{itemize}
\item For any $(z,x) \in Z \times X$, there is an open
neighborhood $U \subseteq Z \times X$ of $(z,x)$ such that $g_v|U
= 0$ for almost every $v \in V$.
\end{itemize}

Given any continuous $g: Z \times X \rightarrow |K|$, we let $U_v
= \{ (z,x) \in Z \times X: g_v(z,x) \neq 0 \}$. Since $Z \times X$
is paracompact, we may find a locally finite refinement $\{ U'_v
\}_{v \in V}$ of the cover $\{ U_v \}_{v \in V}$. Let $\{ g'_v
\}_{v \in V}$ denote any partition of unity subordinate to the
cover $\{U'_v\}_{v \in V}$. The family $\{ g'_v \}$ satisfies all
of the conditions enumerated above, including the strongest local
finiteness assumption, so that it induces a completely continuous
map $g': Z \times X \rightarrow |K|$. To complete the proof, we
need only to show that $g'$ is homotopic to $g$. For this, one
uses a ``straight-line'' homotopy from $g$ to $g'$.
\end{proof}

\section{Dimension Theory}\label{dimension}

In this section, we will discuss the dimension theory of
topological spaces from the point of view of $\infty$-topoi. We
introduce the {\it homotopy dimension} of an $\infty$-topos, and
explain how it relates to various classical notions: covering
dimension for paracompact spaces, cohomological dimension, and
Krull dimension of Noetherian spaces. We also show that finiteness
of the homotopy dimension of an $\infty$-topos $\calX$ has nice
consequences: it implies that every object is the inverse limit of
its Postnikov tower, which proves that $\calX$ is $t$-complete. It
follows that for topological spaces of finite homotopy dimension,
our theory coincides with the Joyal-Jardine theory.

We will conclude by proving a generalization of Grothendieck's
vanishing theorem for the cohomology of abelian sheaves on
Noetherian topological spaces.

\subsection{Homotopy Dimension}

Let $\calX$ be an $\infty$-topos. We shall say that $\calX$ has
{\it homotopy dimension $\leq n$} if every $(n-1)$-connected
object $E \in \calX$ has a global section $1 \rightarrow E$. We
say that $\calX$ has {\it finite homotopy dimension} if there
exists $n \geq 0$ such that $\calX$ has homotopy dimension $\leq
n$.

We shall say that $\calX$ is {\it locally of homotopy dimension
$\leq n$} if every object of $\calX$ can be constructed as a
colimit of objects $E$ with the property that $\calX_{/E}$ has
homotopy dimension $\leq n$. We say that $\calX$ is {\it locally
of finite homotopy dimension} if every object of $\calX$ can be
constructed as a colimit of objects $E$ with the property that
$\calX_{/E}$ has finite homotopy dimension.

\begin{remark}
An $\infty$-topos $\calX$ is (locally or globally) of homotopy
dimension $\leq -1$ if and only if $\calX$ is equivalent to
$\ast$, the $\infty$-category with a single object and a
contractible space of endomorphisms (the $\infty$-category of
stacks on the empty space).
\end{remark}

\begin{remark}
If $\calX$ is a coproduct of $\infty$-topoi $\calX_{\alpha}$, then
$\calX$ is of homotopy dimension $\leq n$ (locally of homotopy
dimension $\leq n$, locally of finite homotopy dimension) if and
only if each $\calX_{\alpha}$ is of homotopy dimension $\leq n$
(locally of homotopy dimension $\leq n$, locally of finite
homotopy dimension).
\end{remark}

\begin{remark}
The $\infty$-topos $\SSet$ is of homotopy dimension $\leq 0$ and
locally of homotopy dimension $\leq 0$. For any object $E \in
\SSet$, the slice $\infty$-topos $\SSet_{/E}$ is of homotopy
dimension $\leq n$ if $E$ can be represented by a CW complex with
cells only in dimensions $\leq n$.
\end{remark}

\begin{lemma}\label{pie}
Let $\calX$ be an $\infty$-topos of homotopy dimension $\leq n$,
and let $E \in \calX$ be $k$-connected. Then $\Hom_{\calX}(1,E)$
is $(k-n)$-connected. In particular, if $E$ is $\infty$-connected,
then $\Hom_{\calX}(1,E)$ is contractible.
\end{lemma}

\begin{proof}
The proof goes by induction on $k$. If $k < n-1$ there is nothing
to prove, and if $k= n-1$ then the assertion follows immediately
from the definition of homotopy dimension. In the general case, we
note that $\Hom_{\calX}(1,E)$ is $(k-n)$-connected if and only if
for any pair of points $p,q \in \Hom_{\calX}(1,E)$, the space of
paths joining $p$ to $q$ is $(k-n-1)$-connected. Since this space
of paths is given by $\Hom_{\calX}(1, 1 \times_{E} 1)$, the result
follows from the inductive hypothesis since $1 \times_{E} 1$ is
$(k-n-1)$-connected.
\end{proof}

If $E$ is an object in an $\infty$-topos $X$, then we have seen
that there is a Postnikov tower $$\rightarrow \tau_n E \rightarrow
\tau_{n-1} E \rightarrow \ldots \rightarrow \tau_0 E \rightarrow
\tau_{-1} E.$$ Let $\tau_{\infty} E$ denote the inverse limit of
this tower. There is a natural map $E \rightarrow \tau_{\infty}
E$.

\begin{proposition}\label{eqvv}
Let $\calX$ be an $\infty$-topos which is locally of finite
homotopy dimension. Then the natural map $E \rightarrow
\tau_{\infty} E$ is an equivalence for any $E \in \calX$.
\end{proposition}

\begin{proof}
We must show that $\Hom_{\calX}(E', E) \rightarrow
\Hom_{\calX}(E', \tau_{\infty} E)$ is an equivalence for every
object $E' \in \calX$. Since $\calX$ is locally of finite homotopy
dimension, it will suffice to prove this in the case where
$\calX_{E'}$ has homotopy dimension $\leq k$.

Let $\eta: E' \rightarrow \tau_{\infty} E$ be an arbitrary map.
For each $n \geq 0$, let $\eta_n$ denote the induced map $E'
\rightarrow \tau_n E$, and let $E_n = E' \times_{\tau_n E} E$.
Since $E \rightarrow \tau_n E$ is $n$-connected, we deduce that
$E_n \rightarrow E'$ is $n$-connected. The space
$\Hom_{\tau_{\infty} E}(E', E)$ is equal to the homotopy inverse
limit of the sequence of spaces $\Hom_{\tau_n E}(E',E) =
\Hom_{E'}(E', E_n)$. Since $\Hom_{E'}(E',E_n)$ is
$(n-k)$-connected by Lemma \ref{pie}, we deduce that
$\Hom_{\tau_{\infty} E}(E',E)$ is contractible, as desired.
\end{proof}

\begin{corollary}\label{fdfd}
If $\calX$ is locally of homotopy dimension $\leq n$, then $\calX$
satisfies hyperdescent.
\end{corollary}

\begin{proof}
Let $S$ denote the collection of $\infty$-connected morphisms of
$\calX$. Every $n$-truncated object of $\calX$ is $S$-local by
Corollary \ref{goober2}. It follows that any limit of truncated
objects of $\calX$ is $S$-local. Since any object $E \in \calX$ is
the inverse limit of its Postnikov tower, we deduce that $E$ is
$S$-local. Thus $\calX \rightarrow S^{-1} \calX$ is an equivalence
of categories.
\end{proof}

\subsection{Cohomological Dimension}

Let $\calX$ be an $\infty$-topos. Given any sheaf $A$ of groups on
$\calX$ (that is, an abelian group object in the ordinary category
$\calX_{1}$), we may construct a groupoid $X_{\bigdot}$ with
$X_{n} = G^n$. Let $K(G,1) = |X_{\bigdot}|$; this is a
$1$-truncated object of $\calX$. A similar construction permits us
to construct Eilenberg-MacLane objects $K(G,n)$ for $n \geq 2$
when $G$ is abelian. We can then define the {\it cohomology} of
the $\infty$-topos $\calX$ with coefficients in $G$ by the formula
$$H^{n}(\calX, G) = \pi_0(\Hom_{\calX}(1, K(G,n))).$$

\begin{remark}
If $\calX$ is the $\infty$-topos of stacks on some ordinary topos
$X$, then $H^{n}(\calX, G) \simeq H^{n}(X,G)$, where the right
hand side denotes sheaf cohomology, defined using an injective
resolution of $G$ in the category of abelian sheaves on $X$. For
some discussion we refer the reader to \cite{jardine}.
\end{remark}

\begin{definition}
Let $\calX$ be an $\infty$-topos. We will say that $\calX$ has
{\it cohomological dimension $\leq n$} if, for any sheaf of
abelian groups $G$ on $\calX$, we have $H^k(\calX,G) = \ast$ for
$k > n$.
\end{definition}

\begin{remark}
For small values of $n$, some authors prefer to require a stronger
vanishing condition which applies also when $G$ is a non-abelian
coefficient system. The appropriate definition requires the
vanishing of cohomology for coefficient groups which are defined
only up to inner automorphisms, as in \cite{giraud}. With the
appropriate modifications, Theorem \ref{cohdim} below remains
valid for $n < 2$.
\end{remark}

In order to study the cohomological dimension, we will need to be
able to recognize Eilenberg-MacLane objects $K(G,n)$. By
construction, we note that $K(G,n)$ is $n$-truncated,
$(n-1)$-connected, and possesses a global section. We will next
prove that these properties characterize the Eilenberg-MacLane
objects. We will deduce this from the following more general
statement:

\begin{proposition}\label{obs}
Let $\calX$ be an $\infty$-topos, and let $E$ and $E'$ be objects
of $\calX$ with base points $\eta: 1 \rightarrow E$ and $\eta': 1
\rightarrow E'$. Suppose that $n \geq 1$, $E$ is $n$-connected,
and $E'$ is $n$-truncated. Then the space of pointed morphisms
from $E$ to $E'$ is equivalent to the $($discrete$)$ space of maps
of sheaves of groups $\eta^{\ast} \pi_n(E) \rightarrow
{\eta'}^{\ast} \pi_n E'$.
\end{proposition}

\begin{proof}
Suppose first that $\calX = \SSet$. In this case, the result
follows from classical obstruction theory. If $\calX$ is an
$\infty$-category of presheaves, then the result can be proved by
working componentwise. In the general case, we apply Theorem
\ref{giraud} to realize $\calX$ as a left-exact localization of
some $\infty$-category $\calP$ of presheaves. Let $f: \calX
\rightarrow \calP$ denote the natural geometric morphism.

The object $f_{\ast} E$ is not necessarily $n$-connected in
$\calP$. However, if we let $F$ denote the mapping fiber $1
\times_{ \tau_{n-1} f_{\ast} E} f_\ast E$, then $F$ is a pointed,
$n$-connected object of $\calP$ and the natural map $F \rightarrow
f_{\ast} E$ becomes an equivalence upon applying $f^{\ast}$.
Consequently, we have $\Hom_{\ast}(E,E') = \Hom_{\ast}(f^{\ast}F,
E') = \Hom_{\ast}(F, f_{\ast} E')$, where $\Hom_{\ast}$ denotes
the base-point compatible morphisms. Since $f_{\ast} E'$ is
$n$-truncated in $\calP$, we deduce that the latter space is
equivalent to the space of morphisms of groups from $\xi^{\ast}
\pi_n F$ to ${\eta'}^\ast \pi_n f_{\ast} E'$, where $\xi: 1
\rightarrow F$ denotes the natural base point.

Since $E$ is an $n$-truncated object of $\calX$, we see that the
space of pointed maps from an $n$-sphere into $E$ is already
discrete. Consequently, no truncation is required in the
definition of $\pi_n E$, and we see that the formation of $\pi_n
E$ is compatible with left-exact functors. It follows that $\pi_n
f_{\ast} E = f_{\ast} \pi_n E$, so that the set of group
homomorphisms from $\xi^{\ast} \pi_n F$ into $(\eta')^{\ast} \pi_n
f_\ast E$ is equivalent to the space of group homomorphisms from
$\xi^{\ast} \pi_n f_{\ast} E \simeq f_{\ast} \eta^{\ast} \pi_n E$
into $f_{\ast} {\eta'}^{\ast} \pi_n E'$. Since $f_{\ast}$ is fully
faithful and compatible with products, this is equivalent to the
set of group homomorphisms from $\eta^{\ast} \pi_n E$ to
${\eta'}^{\ast} \pi_n E'$, as desired.
\end{proof}

We can now give a characterization of the Eilenberg-MacLane
objects in an $\infty$-topos.

\begin{corollary}
Let $\calX$ be an $\infty$-topos, and let $n \geq 1$. Suppose that
$E$ is an object of $\calX$ which is $n$-truncated,
$(n-1)$-connected, and equipped with a section $\eta: 1
\rightarrow E$. Then there exists a canonical equivalence $E
\simeq K( \eta^{\ast} \pi_n E, n)$.
\end{corollary}

Before we can proceed further, we need a lemma.

\begin{lemma}\label{nicelemma}
Let $\calX$ be an $\infty$-topos, $n \geq 0$, and $E$ an
$n$-connected object of $\calX$, and $f: \calX_{/E} \rightarrow
\calX$ the canonical geometric morphism. Then $f^{\ast}$ induces
an equivalence of $\infty$-categories between the
$\infty$-category of $(n-2)$-truncated objects of $\calX$ and the
$(n-2)$-truncated objects of $\calX_{E}$.
\end{lemma}

\begin{proof}
We first prove that $f^{\ast}$ is fully faithful on the
$\infty$-category of $(n-1)$-truncated objects of $\calX$. Let $X,
Y \in \calX$ be objects, where $Y$ is $(n-1)$-truncated. Then
$\Hom_{\calX_{E}}(f^{\ast} X, f^{\ast} Y) = \Hom_{\calX}( E \times
X, Y) = \Hom_{\calX}( E, Y^X )$. Since $Y$ is $(n-1)$-truncated,
$Y^X$ is also $(n-1)$-truncated so that the $n$-connectedness of
$E$ implies that $\Hom_{\calX}( E, Y^X) \simeq \Hom_{\calX}(1,
Y^X) = \Hom_{\calX}(X,Y)$.

Now suppose that $X_E \in \calX_{/E}$ is $(n-2)$-truncated; we
must show that $X_E \simeq f^{\ast} X$ for some $X \in \calX$
(automatically $(n-2)$-connected by descent, so that $X$ is
canonically determined by the first part of the proof). The
natural candidate for $X$ is the object $f_{\ast} X_E$. We will
show that the adjunction $f^{\ast} f_{\ast} X_E \rightarrow X_E$
is an equivalence; then $f_{\ast} X_E$ will be $(n-2)$-connected
by descent and the proof will be complete.

Let $\pi_0, \pi_1: E \times E \rightarrow E$ denote the two
projections. We first claim that there exists an equivalence
$\pi_0^{\ast} X_E \simeq \pi_1^{\ast} X_E$. To see this, we note
that both sides become equivalent to $X_E$ after pulling back
along the diagonal $\delta: E \rightarrow E \times E$. Since $E$
is $n$-connected, the diagonal map $E \rightarrow E \times E$ is
$(n-1)$-connected, so that $\delta^{\ast}$ is fully faithful on
$(n-2)$-truncated objects (by the first part of the proof). Thus
$\pi_0^{\ast} X_E$ and $\pi_1^{\ast} X_E$ must be equivalent to
begin with (in fact canonically, but we shall not need this).

We wish to show that the adjunction morphism $p: f^{\ast} f_{\ast}
X_E$ is an equivalence; then we may take $X = f_{\ast} X_E$. It
suffices to show that $p$ is an equivalence after (surjective)
base change to $E$; in other words, we must show that
$$p_E: \pi_0^{\ast} (\pi_0)_{\ast} \pi_1^{\ast} X_E \rightarrow
\pi_1^{\ast} X_E$$ is an equivalence. Since $\pi_0^{\ast}$ is
fully faithful (on $(n-1)$-truncated objects), we see immediately
that the adjunction morphism $\pi_0^{\ast} (\pi_0)_{\ast} Y
\rightarrow Y$ is an equivalence whenever $Y$ lies in the
essential image of $\pi_0^{\ast}$. Since $\pi_1^{\ast} X_E \simeq
\pi_0^{\ast} X_E$, it lies in the essential image of
$\pi_0^{\ast}$ and the proof is complete.
\end{proof}

Suppose that $E$ is an $(n-1)$-connected, $n$-truncated object of
$\calX$ for $n \geq 2$. Then $\pi_{n}(E)$ is a sheaf of abelian
groups on $E$. By the preceding lemma, there is a unique sheaf of
abelian groups $G$ on $\calX$ such that $\pi_{n}(E) = E \times G$.
In this situation, we shall say that $E$ is an {\it
$(n,G)$-bundle}.

\begin{proposition}
Let $\calX$ be an $\infty$-topos, $n \geq 2$, and $G$ a sheaf of
abelian groups on $\calX$. Then the $(n,G)$-bundles on $X$ are
classified, up to equivalence, by the abelian group
$H^{n+1}(\calX,G)$. Under this equivalence, the identity element
of $H^{n+1}(\calX,G)$ corresponds to $K(G,n)$. An $(n,G)$ bundle
is equivalent to $K(G,n)$ if and only if it admits a global
section.
\end{proposition}

\begin{proof}
Given $\eta \in H^{n+1}(\calX,G)$, we may represent $\eta$ by a
morphism $\eta: 1 \rightarrow K(G,n+1)$, and then form the fiber
product $E_{\eta}= 1 \times_{K(G,n+1)} 1$ using the map $\eta$ and
the base point of $K(G,n+1)$. Since $K(G,n+1)$ is $0$-connected,
$\eta$ is locally equivalent to the base point and $E$ is locally
equivalent to $K(G,n)$. It follows by descent that $E$ is
$(n-1)$-connected and $n$-truncated. The long exact homotopy
sequence gives a canonical identification of $\pi_n E_{\eta}$ with
$G \times E_{\eta}$. Thus, $E_{\eta}$ is an $(n,G)$-bundle.

For each object $E \in \calX$, let $\calC_{E}$ denote the
$\infty$-category of $(n,G)$-bundles in $\calX_{/E}$ (where we
also regard $G$ as a sheaf of groups on $\calX_{/E}$, via the
pullback construction). Morphisms in $\calC_{E}$ are required
induce the identity on $G$. It is easy to see that all morphisms
in $\calC_{E}$ are equivalences, and with a little bit of work one
can show that $\calC_{E}$ is essentially small. Consequently we
may view $\calC_{E}$ as a space which varies contravariantly in
$E$ (via the pullback construction). A descent argument shows that
the functor $E \mapsto \calC_{E}$ carries colimits into limits, so
that by Theorem \ref{representable} this functor is representable
by some object $C \in \calX$. Then the construction of the first
part of the proof gives a map of spaces
$$\phi: K(G,n+1) \rightarrow C$$ which we shall show to be an
equivalence.

Since every $(n,G)$-bundle is locally trivial, we deduce that
$\phi$ is surjective. Hence, to show that $\phi$ is an equivalence
it suffices to prove that $K(G,n+1) \times_{C} K(G,n+1)$ is
equivalent to $K(G,n+1)$. In other words, we must show that given
two maps $\eta, \eta': E \rightarrow K(G,n+1)$, the space of paths
from $\eta$ to $\eta'$ in $\Hom_{\calX}(E, K(G,n+1))$ is
equivalent to the space of identifications between the associated
$(n,G)$-bundles on $E$. Both of these spaces form stacks on $E$,
so to prove that they are equivalent we are free to replace $E$ by
any object which surjects onto $E$. Thus, we may assume that
$\eta$ and $\eta'$ are trivial, in which case the associated
$(n,G)$ bundles are both equivalent to $E \times K(G,n)$. In other
words, we are reduced to proving that the natural map
$$\Hom_{\calX}(E, K(G,n)) \rightarrow \Hom_E(K(G,n) \times E,
K(G,n) \times E)$$ is a homotopy equivalence of the left hand side
onto the collection of components on the right hand side
consisting of maps which induce the identity on $G$.

Replacing $\calX$ by $\calX_{/E}$, we may assume that $E=1$. We
note that $K(G,n)$ has an infinite loop structure and in
particular we have a splitting $\Hom_{\calX}( K(G,n), K(G,n))
\simeq \Hom_{\ast}(K(G,n), K(G,n)) \times \Hom_{\calX}(1,
K(G,n))$, where the first factor denotes the space of pointed
morphisms. By Proposition \ref{obs}, this is a discrete space
consisting of all morphisms of sheaves of groups from $G$ to $G$.
In particular, there is only one point corresponding to the
identity map of $G$, which proves that $\phi$ is an equivalence
onto its image as desired.
\end{proof}

\begin{theorem}\label{cohdim}
Let $\calX$ be an $\infty$-topos and $n \geq 2$. Then $\calX$ has
cohomological dimension $\leq n$ if and only if it satisfies the
following condition: any $(n-1)$-connected, truncated object of
$\calX$ admits a global section.
\end{theorem}

\begin{proof}
Suppose that $\calX$ has the property that every
$(n-1)$-connected, truncated object of $\calX$ admits a global
section. It follows that for any truncated, $n$-connected object
$E$ of $\calX$, $\Hom_{\calX}(1,E)$ is connected. Let $k > n$, and
let $G$ be a sheaf of abelian groups on $\calX$. Then $K(G,k)$ is
$n$-connected, so that $H^k(\calX,G) = \ast$. Thus $\calX$ has
cohomological dimension $\leq n$.

For the converse, let us assume that $\calX$ has cohomological
dimension $\leq n$ and let $E$ denote an $(n-1)$-connected,
$k$-truncated object of $\calX$. We will show that $E$ admits a
global section by induction on $k$. If $k \leq n-1$, then $E = 1$
and there is nothing to prove. For the inductive step, we may
assume that $\tau_{k-1} E$ admits a global section $\eta$. We may
replace $E$ by $1 \times_{\tau_{k-1} E} E$, and thereby assume
that $E$ is $(k-1)$-connected. Let $G = \pi_k E$; then by Lemma
\ref{nicelemma}, $G = E \times G'$ for some sheaf of groups $G'$
on $\calX$. Since $H^{k+1}(\calX,G) = \ast$, we deduce that the $E
\simeq K(G',k)$ and therefore possesses a global section.
\end{proof}

\begin{corollary}\label{confusion}
Let $\calX$ be an $\infty$-topos. If $\calX$ has homotopy
dimension $\leq n$, then $\calX$ has cohomological dimension $\leq
n$. The converse holds provided that $\calX$ has finite homotopy
dimension and $n \geq 2$.
\end{corollary}

\begin{proof}
Only the last claim requires proof. Suppose that $\calX$ has
cohomological dimension $\leq n$ and homotopy dimension $\leq k$.
We must show that every $(n-1)$-connected object $E$ of $\calX$
has a global section. The object $\tau_{k-1} E$ is truncated and
$(n-1)$-connected, so it has a global section by Theorem
\ref{cohdim}. Replacing $E$ by $E \times_{ \tau_{k-1} E} 1$, we
can reduce to the case where $E$ is $(k-1)$-connected, which
follows from the definition of homotopy dimension.
\end{proof}

In the next two sections, we will examine classical conditions
which give bounds on the cohomological dimension and prove that
they also give bounds on the homotopy dimension. We do not know if
every $\infty$-topos of finite cohomological dimension also has
finite homotopy dimension (though this seems unlikely). In
particular, we do not know if the $\infty$-topos of stacks on
$B\hat{Z}$, the classifying topos of the profinite completion of
$\Z$, has finite homotopy dimension. This topos is known to have
cohomological dimension $2$; see for example \cite{serre}.

\subsection{Covering Dimension}

In this section, we will review the classical theory of covering
dimension for paracompact spaces, and then show that the covering
dimension of a paracompact space $X$ coincides with its homotopy
dimension.

Let $X$ be a paracompact space. Recall that $X$ has {\it covering
dimension $\leq n$} if the following condition is satisfied: for
any open covering $\{ U_{\alpha} \}$ of $X$, there exists an open
refinement $\{ V_{\alpha} \}$ of $X$ such that each intersection
$V_{\alpha_0} \cap \ldots \cap V_{\alpha_{n+1}} = \emptyset$
provided the $\alpha_i$ are pairwise distinct.

\begin{remark}
If $X$ is paracompact, the above definition is equivalent to the
(a priori weaker) requirement that such a refinement exist when
$\{U_{\alpha} \}$ is a finite covering of $X$. This weaker
condition gives a good notion whenever $X$ is a normal topological
space. Moreover, if $X$ is normal, then the covering dimension of
$X$ (by this new definition) is equal to the covering dimension of
the Stone-\Cech compactification of $X$. Thus, the dimension
theory of normal spaces is controlled by the dimension theory of
compact Hausdorff spaces.
\end{remark}

\begin{remark}
Suppose that $X$ is a compact Hausdorff space, which is written as
a filtered inverse limit of compact Hausdorff spaces $\{
X_{\alpha} \}$, each of which has dimension $\leq n$. Then $X$ has
dimension $\leq n$. Conversely, any compact Hausdorff space of
dimension $\leq n$ can be written as a filtered inverse limit of
finite simplicial complexes having dimension $\leq n$. Thus, the
dimension theory of compact Hausdorff spaces is controlled by the
(completely straightforward) dimension theory of finite simplicial
complexes.
\end{remark}

\begin{remark}
There are other approaches to classical dimension theory. For
example, a topological space $X$ is said to have {\it small
(large) inductive dimension $\leq n$} if every point of $X$ (every
closed subset of $X$) has arbitrarily small open neighborhoods $U$
such that $\bd U$ has small inductive dimension $\leq n-1$. These
notion are well-behaved for separable metric spaces, where they
coincides with the covering dimension (and with each other). In
general, the covering dimension has better formal properties.
\end{remark}

We now give an alternative characterization of the covering
dimension of a paracompact space. First, we need a technical
lemma.

\begin{lemma}\label{core}
Let $X$ be a paracompact space, $k \geq 0$, $\{U_{\alpha}
\}_{\alpha \in A}$ be a covering of $X$. Suppose that for any $J
\subseteq U_{\alpha}$ of size $k+1$, we are given a covering
$\{V_{J, \beta} \}_{\beta \in B_J}$ of the intersection $U_J =
\bigcap_{\alpha \in J} U_{\alpha}$. Then there exists a covering
$\{ W_{\alpha} \}_{\alpha \in \widetilde{A}}$ of $X$ and a map
$\pi: \widetilde{A} \rightarrow A$ with the following properties:
\begin{itemize}
\item For $\alpha \in \widetilde{A}$, $W_{\alpha} \subseteq
U_{\pi(\alpha)}$.
\item Suppose that $\alpha_0, \ldots, \alpha_k \in
\widetilde{A}$ have the property that $J = \{ \pi(\alpha_0),
\ldots, \pi(\alpha_k) \}$ has cardinality $k+1$. Then there exists
$\beta \in B_J$ such that $W_{\alpha_0} \cap \ldots \cap
W_{\alpha_k} \subseteq V_{J,\beta}$.
\end{itemize}
\end{lemma}

\begin{proof}
Since $X$ is paracompact, we may find a locally finite refinement
$\{ U'_{\alpha} \}_{\alpha \in A}$ which covers $X$, such that the
each closure $\overline{ U'_{\alpha} }$ is contained in
$U_{\alpha}$. Let $S$ denote the set of all subsets $J \subseteq
A$ having size $k+1$. For $J \in S$, let $K_J = \bigcap_{\alpha
\in J} \overline{U_{\alpha}}$. Now let $$\widetilde{A} = \{
(\alpha, J, \beta): \alpha \in A, J \in S, \alpha \in J, \beta \in
B_J \} \coprod A.$$ For $(\alpha, J, \beta) \in \widetilde{A}$, we
set $\pi(\alpha, J, \beta) = \alpha$ and $W_{\alpha, J, \beta} =
(U'_{\alpha} - \bigcup_{\alpha \in J'} K_{J'}) \cup (V_{J, \beta}
\cap U'_{\alpha})$. If $\alpha \in A$, we let $\pi(\alpha) =
\alpha$ and $W_{\alpha} = (U'_{\alpha} - \bigcup_{\alpha \in J'}
K_{J'})$. It is easy to check that this assignment has the desired
properties.
\end{proof}

\begin{theorem}\label{paradimension}
Let $X$ be a paracompact topological space of covering dimension
$\leq n$. Then $X^{\infty}$ has homotopy dimension $\leq n$.
\end{theorem}

\begin{proof}
Suppose that $\calF$ is an $(n-1)$-connected stack on $X$. We must
invoke the notation (but not the conclusions) of Section \S
\ref{paracompactness}. Let $\calB$ denote a basis for $X$
satisfying the conclusions of Lemma \ref{goofy}. Then $\calF$ may
be represented by some $\widetilde{\calF} \in \Simp_{\calB}$. We
wish to produce a global section of $\calF$. It will be sufficient
(and also necessary, by Theorem \ref{main}) to produce a map
$N(\calU) \rightarrow r \widetilde{\calF}$ for some open cover
$\calU$ of $X$ consisting of elements of $\calB$.

We will prove the following statement by induction on $i$, $-1
\leq i \leq n$:

\begin{itemize}
\item There exists an open cover $\{U_{\alpha} \}$ of $X$ (consisting of elements of the basis $\calB$)
and a map from the $i$-skeleton of the nerve of this cover into
$r\widetilde{\calF}$ (in the category $\calK_{\calB}$).
\end{itemize}

Assume that this statement holds for $i = n$. Passing to a
refinement, we may assume that the cover $\{U_{\alpha} \}$ has the
property that no more than $n+1$ of its members intersect (this is
the step where we shall use the assumption on the covering
dimension of $X$). It follows that the $n$-skeleton of the nerve
is the entire nerve, so that we obtain a global section of
$\calF$.

To begin the induction in the case $i = -1$, we use the cover
$\{X\}$; the $(-1)$-skeleton of the nerve of this cover is empty
so there is no data to provide.

Now suppose that we have exhibited the desired cover $\{U_{\alpha}
\}_{\alpha \in A}$ for some value $i < n$. Suppose $J \subseteq A$
has cardinality $i+2$. Over the open set $U_J = \bigcap_{\alpha
\in J} U_\alpha$, our data provides us with a map $f_J: \bd
\Delta^{i+1} \times U_J \rightarrow \calF$. By assumption this map
is locally trivial, so that we may cover $U_J$ by open sets
$\{V_{J,\beta} \}_{\beta \in B_J}$ over which the map $f_J$
extends to a map $f'_{J,\beta}: \Delta^{i+1} \times V_{J,\beta}
\rightarrow \calF$. We apply Lemma \ref{core} to this data, to
obtain an new open cover $\{ W_{\alpha} \}_{\alpha \in
\widetilde{A} }$ which refines $\{ U_{\alpha} \}_{\alpha \in A}$.
Refining the cover further if necessary, we may assume that each
of its members belongs to $\calB$. By functoriality, we obtain a
map $f$ from the $i$-skeleton of the nerve of $\{ W_{\alpha}
\}_{\alpha \in \widetilde{A}}$ to $\calF$. To complete the proof,
it will suffice to extend $f$ to the $(i+1)$-skeleton of the nerve
of $\{ W_{\alpha} \}_{\alpha \in \widetilde{A}}$. Let $\pi:
\widetilde{A} \rightarrow A$ denote the map of Lemma \ref{core},
and consider any $(i+2)$-element subject $J \subseteq
\widetilde{A}$. If $\pi(J)$ has size $< i+2$, then composition
with $\pi$ gives us a canonical extension of $f$ to $\Delta^{i+1}
\times \bigcap_{j \in J} W_j$ (since the corresponding simplex
becomes degenerate with respect to the covering $\{ U_{\alpha}
\}$). If $\pi(J)$ has size $i+2$, then Lemma \ref{core} assures us
that $\bigcap_{j \in J} W_j$ is contained in some $V_{\pi(J),
\beta}$, so that the desired extension exists.
\end{proof}

\begin{remark}
In fact, the inequality of Theorem \ref{paradimension} is an
equality: if the homotopy dimension of $X$ is $\leq n$, then the
cohomological dimension of $X$ is $\leq n$, so that the
paracompact space $X$ has covering dimension $\leq n$.
\end{remark}

\subsection{Heyting Dimension}

Let $X$ be a topological space. Recall that $X$ is {\it
Noetherian} if the collection of closed subsets of $X$ satisfies
the descending chain condition. In particular, we see that the
irreducible closed subsets of $X$ form a well-founded set.
Consequently there is a unique ordinal-valued rank function $r$,
defined on the irreducible closed subsets of $X$, having the
property that $r(C)$ is the smallest ordinal which is larger than
$r(C_0)$ for any proper subset $C_0 \subset C$. We call $r(C)$ the
{\it Krull dimension} of $C$, and we call $\sup_{C \subseteq X}
r(C)$ the {\it Krull dimension} of $X$.

Let $X$ be a topological space. We shall say that $X$ is a {\it
Heyting space} if satisfies the following conditions:

\begin{itemize}
\item The space $X$ has a basis consisting of compact open sets.
\item The compact open subsets of $X$ are stable under finite
intersections.
\item If $U$ and $V$ are compact open subsets of $X$, then the
interior of $U \cup (X-V)$ is compact.
\item Every irreducible closed subset of $X$ has a unique generic
point (in other words, $X$ is a {\it sober} topological space).
\end{itemize}

\begin{remark}
The last condition is not very important. Any topological space
which does not satisfy this condition can be replaced by a
topological space which does, without changing the lattice of open
sets.
\end{remark}

\begin{remark}
Recall that a {\it Heyting algebra} is a distributive lattice $L$
with the property that for any $x,y \in L$, there exists a maximal
element $z$ with the property that $x \wedge z \subseteq y$. It
follows immediately from our definition that the lattice of
compact open subsets of a Heyting space forms a Heyting algebra.
Conversely, given any Heyting algebra one may form its spectrum,
which is a Heyting space. This sets up a duality between the
category of Heyting spaces and the category of Heyting algebras,
which is a special case of a more general duality between coherent
topological spaces and distributive lattices. We refer the reader
to \cite{johnstone} for more details.
\end{remark}

\begin{remark}
Suppose that $X$ is a Noetherian topological space in which every
irreducible closed subset has a unique generic point. Then $X$ is
a Heyting space, since every open subset of $X$ is compact.
\end{remark}

\begin{remark}
If $X$ is a Heyting space and $U \subseteq X$ is a compact open
subset, then $X$ and $X-U$ are also Heyting spaces. In this case,
we say that $X-U$ is a {\it cocompact} closed subset of $X$.
\end{remark}

We next define the dimension of a Heyting space. The definition is
recursive. Let $\alpha$ be an ordinal. A Heyting space $X$ has
{\it Heyting dimension $\leq \alpha$} if and only if, for any
compact open subset $U \subseteq X$, the boundary of $U$ has
Heyting dimension $< \alpha$ (we note that the boundary of $U$ is
also a Heyting space); a Heyting space has dimension $< 0$ if and
only if it is empty.

\begin{remark}
A Heyting space has dimension $\leq 0$ if and only if it is
Hausdorff. The Heyting spaces of dimension $\leq 0$ are precisely
the compact, totally disconnected Hausdorff spaces. In particular,
they are also paracompact spaces and their Heyting dimension
coincides with their covering dimension.
\end{remark}

\begin{proposition}\label{closs}

\begin{enumerate}
\item Let $X$ be a Heyting space of dimension $\leq \alpha$. Then
for any compact open subset $U \subseteq X$, both $U$ and $X-U$
have Heyting dimension $\leq \alpha$.

\item Let $X$ be a Heyting space which is a union of finitely many
compact open subsets $U_{\alpha}$ of dimension $\leq \alpha$. Then
$X$ has dimension $\leq \alpha$.

\item Let $X$ be a Heyting space which is a union of finitely many
cocompact closed subsets $K_{\alpha}$ of Heyting dimension $\leq
\alpha$. Then $X$ has Heyting dimension $\leq \alpha$.
\end{enumerate}
\end{proposition}

\begin{proof}
All three assertions are proven by induction on $\alpha$. The
first two are easy, so we restrict our attention to $(3)$. Let $U$
be a compact open subset of $X$, having boundary $B$. Then $U \cap
K_{\alpha}$ is a compact open subset of $K_{\alpha}$, so that the
boundary $B_{\alpha}$ of $U \cap K_{\alpha}$ in $K_{\alpha}$ has
dimension $\leq \alpha$. We see immediately that $B_{\alpha}
\subseteq B \cap K_{\alpha}$, so that $\bigcup B_{\alpha}
\subseteq B$. Conversely, if $b \notin \bigcup B_{\alpha}$ then,
for every $\beta$ such that $b \in K_{\beta}$, there exists a
neighborhood $V_{\beta}$ containing $b$ such that $V_{\beta} \cap
K_{\beta} \cap U = \emptyset$. Let $V$ be the intersection of the
$V_{\beta}$, and let $W = V - \bigcup_{b \notin K_{\gamma}}
K_{\gamma}$. Then by construction, $b \in W$ and $W \cap U =
\emptyset$, so that $b \in B$. Consequently, $B = \bigcup
B_{\alpha}$. Each $B_{\alpha}$ is closed in $K_{\alpha}$, thus in
$X$ and also in $B$. The hypothesis implies that $B_{\alpha}$ has
dimension $< \alpha$. Thus the inductive hypothesis guarantees
that $B$ has dimension $< \alpha$, as desired.
\end{proof}

\begin{remark}\label{DVR}
It is not necessarily true that a Heyting space which is a union
of finitely many {\em locally closed} subsets of dimension $\leq
\alpha$ is also of dimension $\leq \alpha$. For example, a
topological space with $2$ points and a nondiscrete, nontrivial
topology has Heyting dimension $1$, but is a union of two locally
closed subsets of Heyting dimension $0$.
\end{remark}

\begin{proposition}
If $X$ is a sober Noetherian topological space, then the Krull
dimension of $X$ coincides with the Heyting dimension of $X$.
\end{proposition}

\begin{proof}
We first prove, by induction on $\alpha$, that if the Krull
dimension of a sober Noetherian space is $\leq \alpha$, then the
Heyting dimension of a sober Noetherian space is $\leq \alpha$.
Let $X$ be a Noetherian topological space of Krull dimension $\leq
\alpha$. Using Proposition \ref{closs}, we may assume that $X$ is
irreducible. Consider any open subset $U \subseteq X$, and let $Y$
be its boundary. We must show that $Y$ has Heyting dimension $\leq
\alpha$. Using Proposition \ref{closs} again, it suffices to prove
this for each irreducible component of $Y$. Now we simply apply
the inductive hypothesis and the definition of the Krull
dimension.

For the reverse inequality, we again use induction on $\alpha$.
Assume that $X$ has Heyting dimension $\leq \alpha$. To show that
$X$ has Krull dimension $\leq \alpha$, we must show that every
irreducible closed subset of $X$ has Krull dimension $\leq
\alpha$. Without loss of generality we may assume that $X$ is
irreducible. Now, to show that $X$ has Krull dimension $\leq
\alpha$, it will suffice to show that any {\em proper} closed
subset $K \subseteq X$ has Krull dimension $< \alpha$. By the
inductive hypothesis, it will suffice to show that $K$ has Heyting
dimension $< \alpha$. By the definition of the Heyting dimension,
it will suffice to show that $K$ is the boundary of $X - K$. In
other words, we must show that $X - K$ is dense in $X$. This
follows immediately from the irreducibility of $X$.
\end{proof}

We now prepare the way for our vanishing theorem. First, we
introduce a modified notion of connectivity:

\begin{definition}
Let $X$ be a Heyting space and $k$ any integer. A stack $\calF$ on
a compact open set $V \subseteq X$ is {\it strongly $k$-connected}
if the following condition is satisfied: for any $m \geq -1$, any
compact open $U \subseteq V$, and any map $\phi: S^{m} \rightarrow
\calF(U)$, there exists a cocompact closed subset $K \subseteq U$
such that $\overline{K} \subseteq X$ has dimension $< m-k$ and an
open cover $\{ V_{\alpha} \}$ of $V-K$ such that the restriction
of $\phi$ to $\calF(V_{\alpha})$ is nullhomotopic for each
$\alpha$ (if $m=-1$, this means that $\calF(V_{\alpha})$ is
nonempty).
\end{definition}

\begin{remark}
Since the definition involves taking the closure of $K$ in $X$,
rather than in $V$, we note that the strong connectivity of
$\calF$ may increase if we replace $X$ by $V$.
\end{remark}

\begin{remark}
Strong $k$-connectivity is an unstable analogue of the
connectivity conditions on complexes of sheaves, associated to a
the dual of the standard perversity (which is well-adapted to the
pushforward functor). For a discussion of perverse sheaves in the
abelian context, see for example \cite{deligne}.
\end{remark}

\begin{remark}
Suppose $X$ has Heyting dimension $\leq n$. If $\calF$ is
$(k+n)$-connected, then for $m \leq k+n$ we may take $K =
\emptyset$ and for $m > k+n$ we may take $K=V$. It follows that
$\calF$ is strongly $k$-connected. Conversely, it is clear from
the definition that strong $k$-connectivity implies
$k$-connectivity.
\end{remark}

The strong $k$-connectivity of $\calF$ is, by construction, a
local property. The key to our vanishing result is that this is
equivalent to a stronger {\it global} property.

\begin{theorem}\label{vanishing}
Let $X$ be a Heyting space of dimension $\leq n$, let $W \subseteq
X$ be a compact open set, and let $\calF$ be a stack on $W$. The
following conditions are equivalent:

\begin{enumerate}

\item For any compact open $U \subseteq V \subseteq W$, any $m
\geq -1$, map $\zeta: S^m \rightarrow \calF(V)$ and any
nullhomotopy $\eta$ of $\eta|U$, there exists an extension of
$\eta$ to $V-K$, where $K \subseteq V$ is a cocompact closed
subset and $\overline{K} \subseteq X$ has dimension $< m-k$.

\item For any compact open $V \subseteq W$, any $m \geq -1$, and
any map $\zeta: S^m \rightarrow \calF(V)$, there exists a
nullhomotopy of $\zeta$ on $V-K$, where $K \subseteq V$ is a
cocompact closed subset and $\overline{K} \subseteq X$ has
dimension $< m-k$.

\item The stack $\calF$ is strongly $k$-connected.
\end{enumerate}
\end{theorem}

\begin{proof}
It is clear that $(1)$ implies $(2)$ (take $U$ to be empty) and
that $(2)$ implies $(3)$ (by definition). We must show that $(3)$
implies $(1)$. Replacing $W$ by $V$ and $\calF$ by the stack
$\calF|V \times_{\calF|V \otimes S^m} 1$, we may reduce to the
case where $W=V$ and $m=-1$.

The proof goes by induction on $k$. For our base case, we take
$k=-n-2$, so that there is no connectivity assumption on the stack
$\calF$. We are then free to choose $K = X-U$ (it is clear that
$\overline{K}$ has dimension $\leq n$).

Now suppose that the theorem is known for strongly
$(k-1)$-connected stacks on any compact open subset of $X$; we
must show that for any strongly $k$-connected $\calF$ on $V$ and
any $\eta \in \calF(U)$, there exists an extension of $\eta$ to
$V-K$ where $\overline{K} \subseteq X$ has dimension $< -1 -k$.

Since $\calF$ is strongly $k$-connected, we deduce that there
exists an open cover $\{V_{\alpha} \}$ of some open subset
$V-K_0$, where $K_0$ has dimension $< -1-k$ in $X$, together with
elements $\psi_{\alpha} \in \calF( V_{\alpha})$. Replacing $V$ by
$V- K_0$ and $U$ by $U- K_0 \cap U$, we may suppose $K_0 =
\emptyset$.

Since $V$ is compact, we may assume that there exist only finitely
many indices $\alpha$. Proceeding by induction on the number of
indices, we may reduce to the case where $V = U \cup V_{\alpha}
\cup K_0$ for some $\alpha$. Let $\calF'$ denote the sheaf on $U
\cap V_{\alpha}$ of {\it paths from $\psi_{\alpha}$ to $\eta$} (in
other words, $\calF' = 1 \times_{\calF} 1$, where the maps $1
\rightarrow \calF$ are given by $\psi_{\alpha}$ and $\eta$). Then
$\calF'$ is strongly $(k-1)$-connected, so by the inductive
hypothesis (applied to $\emptyset \subseteq U \cap V_{\alpha}
\subseteq X$), there exists a closed subset $K \subset U \cap
V_{\alpha}$ having dimension $< -k$ in $X$, such that
$\psi_{\alpha}$ and $\eta$ are equivalent on $(U \cap V_{\alpha})
- K$. Since $\overline{K}$ has dimension $< -k$ in $X$, the
boundary $\bd K$ of $K$ has codimension $< -k-1$ in $X$. Let $W =
V_{\alpha} \cap (X - \overline{K})$. Then we can glue $\eta$ and
$\psi_{\alpha}$ to obtain a global section of $\calF$ over $W \cup
U$, which contains $V - \bd K = X - \bd K$.
\end{proof}

\begin{corollary}\label{hphp}
Let $\pi: X \rightarrow Y$ be a continuous map between Heyting
spaces of finite dimension. Suppose that $\pi$ has the property
that for any cocompact closed subset $K \subseteq X$ of dimension
$\leq n$, $\pi(K)$ is contained in a cocompact closed subset of
dimension $\leq n$. Then the functor $\pi_{\ast}$ carries strongly
$k$-connected stacks into strongly $k$-connected stacks.
\end{corollary}

\begin{proof}
This is clear from the characterization $(2)$ of Theorem
\ref{vanishing}.
\end{proof}

\begin{corollary}\label{gra}
Let $X$ be a Heyting space of finite dimension, and let $\calF$ be
a strongly $k$-connected stack on $X$. Then $\calF(X)$ is
$k$-connected.
\end{corollary}

\begin{proof}
Apply Corollary \ref{hphp} in the case where $Y$ is a point.
\end{proof}

\begin{corollary}
If $X$ is a Heyting space of Heyting dimension $\leq n$, then the
$\infty$-topos $X^{\infty}$ has homotopy dimension $\leq n$.
\end{corollary}

\begin{remark}
Since every compact open subset of $X$ is also a Heyting space of
Heyting dimension $\leq n$, we can deduce also that the
$\infty$-topos $\stacks X$ is locally of homotopy dimension $\leq
n$, and therefore is $t$-complete.
\end{remark}

In particular, we obtain Grothendieck's vanishing theorem:

\begin{corollary}
Let $X$ be a Noetherian topological space of Krull dimension $\leq
n$. Then $X$ has cohomological dimension $\leq n$.
\end{corollary}

\begin{example}
Let $V$ be a real algebraic variety (defined over the real
numbers, say). Then the lattice of open subsets of $V$ that can be
defined by polynomial equations and inequalities is a Heyting
algebra, and the spectrum of this Heyting algebra is a Heyting
space $X$ having dimension at most equal to the dimension of $V$.
The results of this section therefore apply to $X$.

More generally, let $T$ be an o-minimal theory (see for example
\cite{lou}), and let $S_n$ denote the set of complete $n$-types of
$T$. We endow $S_n$ with the following topology generated by the
sets $U_{\phi} = \{p: \phi \in p\}$, where $\phi$ ranges over
formula with $n$ free variables such that $T$ proves that the set
of solutions to $\phi$ is an open set. Then $S_n$ is a Heyting
space of dimension $\leq n$.
\end{example}

\begin{remark}
The methods of this section can be adapted to slightly more
general situations, such as the Nisnevich topology on a Noetherian
scheme of finite Krull dimension. It follows that the
$\infty$-topoi associated to such sites have finite homotopy
dimension and hence our theory is equivalent to the Joyal-Jardine
theory.
\end{remark}

\section{Appendix}
\subsection{Homotopy Limits and Colimits}\label{appendixdiagram}

In this appendix, we will summarize (with sketches of proofs) some
basic facts regarding limits and colimits in $\infty$-categories.
These are analogous to well-known facts concerning limits and
colimits in ordinary categories. We will limit the discussion to
colimits; the corresponding statements for limits may be obtained
by passing to opposite categories.

In this paper, we make use of colimits where the diagrams are
indexed by $\infty$-categories. Many authors use only ordinary
categories to index their colimits, which seems at first to give a
less general notion. In fact, this is not the case:

\begin{proposition}\label{puffy}
Let $\kappa$ be an infinite regular cardinal. Then all
$\kappa$-small homotopy colimits can be constructed using:
\begin{itemize}
\item Coequalizers and $\kappa$-small sums.
\item Pushouts and $\kappa$-small sums.
\end{itemize}
\end{proposition}

\begin{proof}
The coequalizer of a pair of morphisms $p,q: X \rightarrow Y$
coincides with the pushout of $X$ and $Y$ over $X \coproduct X$
(which maps to $Y$ via $p \coproduct q$). Given morphisms $X
\rightarrow Y$ and $X \rightarrow Z$, the pushout $Y \coproduct_X
Z$ may be constructed as the coequalizer of the induced maps $X
\rightarrow Y \coproduct Z$. Thus, coequalizers and coproducts can
be constructed in terms of one another (and finite sums); it
therefore suffices to prove the assertion regarding pushouts.

We next note that from finite sums and pushouts, one can construct
$K \otimes X$ where $K$ is any finite cell complex. We will need
this in the case where $K$ is the sphere $S^n$. If $n=-1$, then $K
\otimes S^n$ is the initial object (the empty sum). For $n \geq
0$, we note that $S^n \otimes K$ is the pushout of $K$ with itself
over $S^{n-1} \otimes K$.

Let $\calC$ be an $\infty$-category, and let $F: \calD \rightarrow
\calC$ be a diagram in $\calC$ with $\calD$ $\kappa$-presented. We
may view $\calD$ as obtained by a process of ``attaching
$n$-cells" as $n$ goes from $0$ to $\infty$, where a $0$-cell is
an object of $\calD$, a $1$-cell is a morphism of $\calD$, and
more generally adjoining an $n$-cell to $\calD$ means adjoining an
$(n-1)$-cell to some mapping space $\Hom_{\calD}(x,y)$ (together
with a mess of other morphisms which are generated by new cell).
The hypothesis that $\calD$ is $\kappa$-generated means that we
may construct $\calD$ using fewer than $\kappa$ $n$-cells for each
$n$, and if $\kappa=\omega$ then we use only finitely many cell in
total.

If $\calD$ is obtained from a ``subcategory'' $\calD_0$ by
adjoining an $(n-1)$-cell to $\Hom(x,y)$, then the colimit of $F$
in $\calC$ may be obtained as a pushout $X \coproduct_{S^{n-2}
\otimes Fx} Fx$, where $X$ is the colimit of $F|\calD_0$. If
$\calD$ is obtained from $\calD_0$ by adjoining a new object $x$,
then the colimit of $F$ is obtained from the colimit of
$F|\calD_0$ by taking the sum with $x$. By induction on the number
of cells needed to construct $\calD$, we can complete the proof in
the case where $\kappa=\omega$.

For $\kappa > \omega$, we must work a little bit harder. First of
all, we note that if $\calD$ is obtained from $\calD_0$ by
adjoining any $\kappa$-small collection of cells simultaneously,
then we may construct the colimit of $F$ as a pushout $X
\coproduct_{Y} Z$ as above, where $X$ is the colimit of
$F|\calD_0$ and the pair $(Y,Z)$ is a sum of pairs of the form $(
S^{n} \otimes Fx, Fx)$ considered above. In general, $\calD$ can
be obtained by an infinite sequence of simultaneous
cell-attachments, so that the colimit of $F$ is may be written as
a direct limit of objects $X_n = \colim F|\calD_n$, and each $X_n$
may be constructed using pushouts and $\kappa$-small sums by the
argument sketched above.

It therefore suffices to show that we may construct the direct
limit of a sequence $X_0 \rightarrow X_1 \rightarrow \ldots$. But
this is easy: the direct limit of the $\{X_i\}$ can be written as
a coequalizer of the identity and the shift map
$$\coproduct_i X_i \rightarrow \coproduct_i X_i.$$
It is permissible to form this countably infinite coproduct since
$\kappa > \omega$.
\end{proof}

\begin{remark}
If $\kappa = \omega$, then finite sums may be constructed by
taking pushouts over an initial object. Consequently, all finite
colimits may be constructed using an initial object and pushouts.
\end{remark}

\subsection{Filtered Colimits}

We now sketch a proof that the general notion of a filtered
colimit is not really any more general than ``directed colimits'',
which the reader will surely find familiar.

More specifically, we sketch a proof of the following result,
which is used in the body of the paper:

\begin{proposition}
Let $\kappa$ be a regular cardinal, and let $\calD$ be a small
$\kappa$-filtered $\infty$-category. Then there exists a
$\kappa$-directed partially ordered set (which we may regard as an
ordinary category, and therefore an $\infty$-category) $\calD_0$,
and a functor $F: \calD_0 \rightarrow \calD$ having the following
property: for any functor $G: \calD \rightarrow \calC$, the
diagrams $G$ and $G \circ F$ have the same colimits.
\end{proposition}

\begin{proof}
We will sketch a construction of $\calD_0$ and leave the property
to the reader. The construction uses the Boardman-Vogt approach to
$\infty$-categories based on simplicial sets satisfying a weak Kan
condition (see \cite{quasicat} and \cite{weakKan}). According to
this approach, we may model $\calD$ by a simplicial set $D$
(satisfying certain extension conditions).

Suppose that $D' \subseteq D$ is a simplicial subset and $v \in
D'_0$. We shall say that $v$ is {\it final} in $D'$ if every map
$\bd \Delta^n \rightarrow D'$ which carries the last vertex of
$\Delta^n$ into $v$ admits an extension $\Delta^n \rightarrow D'$.
We call $D' \subseteq D$ {\it $\kappa$-small} if it has fewer than
$\kappa$ nondegenerate vertices.

We let $\calD_0$ denote the collection of all $\kappa$-small
simplicial subsets $D' \subseteq D$ which possess a final vertex.
The set $\calD_0$ is partially ordered by inclusion. We claim that
it is $\kappa$-directed: in other words, any subset $S \subseteq
\calD_0$ of size $< \kappa$ has an upper bound. It suffices to
show that {\em any} $\kappa$-small simplicial subset $D' \subseteq
D$ can be enlarged to a $\kappa$-small simplicial subset with a
final vertex. This is more or less equivalent to the assertion
that $\calD$ is $\kappa$-directed.

The functor $F: \calD_0 \rightarrow \calD$ is obtained by choosing
a final vertex from each $D' \in \calD_0$. Although the final
vertex is not necessarily unique, it is unique up to a
contractible space of choices (essentially by definition) so that
the functor $F$ is well-defined.
\end{proof}

\subsection{Enriched $\infty$-Categories and
$(\infty,2)$-Categories}

In this section, we discuss $(\infty,2)$-categorical limits and
colimits, which were mentioned several times in the body of the
paper. Before we can describe these, we need to decide what an
$(\infty,2)$-category {\em is}. Just as an ordinary $2$-category
may be described as a category ``enriched over categories'', we
can obtain the theory of $(\infty,2)$-categories as a special case
of enriched $\infty$-category theory.

Given any $\infty$-category $\calE$, one can define a notion of
``$\infty$-categories enriched over $\calE$''. A $\calE$-enriched
$\infty$-category $\calC$ consists of a collection of objects,
together with an object $\Hom_{\calC}(X,Y) \in \calE$ for every
pair of objects $X,Y \in \calC$. Finally, these objects of $\calE$
should be equipped with {\em coherently associative} composition
products

$$\Hom_{\calC}(X_0,X_1) \times \ldots \times
\Hom_{\calC}(X_{n-1},X_n) \rightarrow \Hom_{\calC}(X_0, X_n).$$

\begin{remark}
More generally, one can replace the Cartesian product $\times$ by
any coherently associative monoidal structure $\otimes$ on
$\calE$. If $\calE$ is an ordinary category, then we arrive at the
usual notion of an $\calE$-enriched category.
\end{remark}

Taking $\calE$ to be the $\infty$-category $\calS$ of spaces, we
recover the notion of an $(\infty,1)$-category. On the other hand,
if we take $\calE$ to be the $\infty$-category of
$(\infty,1)$-categories (regarded as an $(\infty,1)$-category by
ignoring noninvertible natural transformations), then we obtain
the notion of an $(\infty,2)$-category. Iterating this procedure
leads to a notion of $(\infty,n)$-category for all $n$. The
formalization of this process leads to the Simpson-Tamsamani
theory of higher categories: see \cite{tamsamani}.

\begin{remark}
The ``$\infty$-category of $\infty$-categories'' is much more
natural than its $1$-categorical analogue, the ``category of
categories''. In the latter case, one should really have a
$2$-category: it is unnatural to ask for two functors to be equal,
and we lose information by ignoring the natural transformations.
In the $\infty$-categorical case, much less information is lost
since we discard only the noninvertible natural transformations.
\end{remark}

Now, if $\calC$ is any $\infty$-category enriched over $\calE$,
then any object of $\calC$ represents a $\calE$-valued functor on
$\calC$. Thus we obtain an ``enriched Yoneda embedding'' $\calC
\rightarrow \calE^{\calC^{op}}$, where the target denotes the
$\calE$-enriched $\infty$-category of $\calE$-valued presheaves on
$\calC$ (we consider only presheaves which are compatible with the
enrichment; that is, functors $\calC^{op} \rightarrow \calE$
between $\calE$-enriched categories).

Given another (small) $\infty$-category $\calD$ which is enriched
over $\calE$, and a functor $F: \calD \rightarrow \calE$, we
obtain a pullback functor $F^{\ast}: \calE^{\calC^{op}}
\rightarrow \calE^{\calD^{op}}$ between the $\infty$-categories of
$\calE$-valued presheaves. Composing this with Yoneda embedding,
we obtain a functor $T: \calC \rightarrow \calE^{\calD^{op}}$.
Finding colimits in $\calC$ amounts to finding a left-adjoint to
$T$. That is, given any $\calE$-valued presheaf $\calF$ on
$\calD$, the colimit of $(F, \calF)$ is an object $X \in \calC$
equipped with a map $\calF \rightarrow TX$ which is universal, in
the sense that it induces equivalences
$$\Hom_{\calC}(X,Y) \rightarrow \Hom_{\calE^{\calD^{op}}}(\calF,
TY)$$ where both sides are regarded as living in the
$\infty$-category $\calE$.

\begin{remark}
The above notion of colimit is apparently more general than the
notion we have considered in the body of the paper for the case
where $\calE = \SSet$. This generality is only apparent, however:
we can replace a pair $(\calD, \calF)$ by the $\infty$-category
$\calD_{\calF}$ ``fibered in spaces'' over $\calD$; then a colimit
of $(F, \calF)$ in the sense just defined is the same thing as a
colimit of $F|\calD_{\calF}$. This does not work in the
$\calE$-enriched case because there is no analogous ``fibered
category'' construction, so we must work with the more general
notion.
\end{remark}

Using a variant of the argument of the proof of Proposition
\ref{puffy}, one can show that all $\calE$-colimits in $\calC$ can
be constructed using pushouts, sums, and tensor products $\calE
\otimes \calC \rightarrow \calC$.

\subsection{Free Algebras}

Granting a good theory of $(\infty,2)$-categories, we can discuss
{\it monoidal $\infty$-categories}. As in classical category
theory, these may be thought of either as $(\infty,2)$-categories
with a single (specified) object, or as $(\infty,1)$-categories
equipped with a coherently associative product which we shall
denote by $\otimes$. We shall adopt the second point of view.

A {\it monoidal functor} $F: \calC \rightarrow \calC'$ between two
monoidal $\infty$-categories is a functor which is compatible with
the associative products $\otimes$ on $\calC$ and $\calC'$, up to
equivalence. (Alternatively, from the $(\infty,2)$-categorical
point of view, we may view a monoidal functor as simply a functor
$F$ between the underlying $2$-categories {\em together with} an
identification $\ast \simeq F(\ast)$, where $\ast$ is used to
denote the specified objects.)

Let $\calC$ be a monoidal $\infty$-category. Suppose that $\calC$
is presentable and $A \in \calC$. Using Theorem
\ref{representable}, we see that the following conditions are
equivalent:

\begin{itemize}
\item The functor $X \mapsto A \otimes X$ commutes with the
formation of colimits.

\item For each $B \in \calC$, there exists an object $\Hom(A,B)
\in \calC$ and a natural equivalence
$$\Hom_{\calC}( \bigdot, \Hom(A,B) ) \simeq \Hom_{\calC}(A \otimes
\bigdot, B) $$
\end{itemize}

If these equivalent conditions and their duals (involving products
where $A$ appears on the {\em right}) are satisfied for each $A
\in \calC$, then we shall say that $\otimes$ is {\it
colimit-preserving}.

If $\calC$ is any monoidal $\infty$-category, one may define an
$\infty$-category $M(\calC)$ of {\it monoid objects} of $\calC$.
The objects of $M(\calC)$ are objects $A \in \calC$ together with
multiplications $A^{\otimes n} \rightarrow A$ ($n \geq 0$),
together with various homotopies and higher homotopies which
explicate the idea that $A$ should be associative up to coherent
homotopy.

The purpose of this appendix is to sketch a proof of the following
fact which is needed in the body of the paper:

\begin{proposition}\label{free}
\begin{itemize}
\item  Let $\calC$ be a presentable monoidal $\infty$-category,
such that the monoidal structure $\otimes$ is colimit-preserving.
Let $A \in M(\calC)$, $N \in \calC$, and let $f: A \rightarrow N$
be a morphism in $\calC$. Then there exists an object $B(f) \in
M(\calC)$ and an equivalence
$$ \Hom_{M(\calC)}(B(f), \bigdot) \simeq \Hom_{M(\calC)}(A, \bigdot)
\times_{ \Hom_{\calC}(A, \bigdot)} \Hom_{\calC}(M, \bigdot).$$

\item Suppose that $\calC$ and $\calC'$ are both presentable
monoidal $\infty$-categories with colimit-preserving monoidal
structures. Let $F: \calC \rightarrow \calC'$ be a monoidal,
colimit-preserving functor. Let $A \in M(\calC)$, $N \in \calC$,
and $f: A \rightarrow M$ a morphism in $\calC$. Then the natural
map $B(Ff) \rightarrow F B(f)$ is an equivalence.
\end{itemize}
\end{proposition}

\begin{proof}
To prove the first assertion, we will sketch a construction of
$B(f)$. Moreover, in our construction, we will make use only of
the monoidal structure on $\calC$ and colimits in $\calC$. The
second part will then follow immediately.

For the construction, we will need to make use of a monoidal
$\infty$-category which is, in some sense, freely generated by a
monoid object, a second object, and a map between them. This free
$\infty$-category turns out to be an ordinary category, which is
easy to construct explicitly. We will proceed by giving an
explicit description of this free category, but we note that this
description is largely irrelevant to the proof, for which only the
existence is needed.

Let $\calD$ denote the ordinary category whose objects are finite,
linearly ordered sets with a distinguished subset of marked
points. If $J,J' \in \calD$, then the morphisms from $J$ to $J'$
in $\calD$ are maps $f: J \rightarrow J'$ satisfying the following
conditions:

\begin{itemize}
\item If $j \leq k$, then $f(j) \leq f(k)$.

\item If $j \in J$ is marked, then $f(j)$ is marked and
$f^{-1}f(j) = \{j\}$.
\end{itemize}

We define a functor $G: \calD \rightarrow \calC$ as follows. Given
an object $J = \{ j_1, \ldots, j_n \} \in \calD$, with marked
points $\{ j_{i_1}, \ldots, j_{i_k} \}$, we set $$G(J) = A_1
\otimes \ldots \otimes A_{j_{i_1}-1} \otimes M_{j_1} \otimes
\ldots \otimes M_{j_{i_k}} \otimes A_{j_{i_k}+1} \otimes \ldots
\otimes A_n.$$

In words, we simply take a large tensor product, indexed by $J$,
of copies of $A$ (for the unmarked points) and copies of $M$ (for
the marked points). Here the subscripts are introduced to index
several copies of $A$ and $M$ that are used in the tensor product.
The definition of $G$ on morphisms is encoded in the monoid
structure on $A$ and the map $f: A \rightarrow M$.

Note that $\calD$ has the structure of a monoidal category, with
the monoidal structure given by concatenation. With respect to
this monoidal structure, $G$ is a monoidal functor.

Let $G^n$ denote the induced functor $G \otimes \ldots \otimes G:
\calD^n \rightarrow \calC$, and let $B(f)^{\otimes n}$ denote the
colimit of $G^n$, regarded as a diagram in $G$ (in the exceptional
case where $n = 0$, we instead define $B(f)^{\otimes 0}$ to be the
unit in $\calC$). Using the fact that $\otimes$ is
colimit-preserving in $\calC$, we deduce that $B(f)^{\otimes n}$
is actually a tensor product of $n$ copies of $B(f) =
B(f)^{\otimes 1}$. Moreover, using the fact that $G$ is a monoidal
functor, we deduce that $G = G^n \circ \delta$ (where $\delta$
denotes the diagonal $\calD \rightarrow \calD^n$), so that we get
a family of maps $B(f)^{\otimes n} \rightarrow B(f)$; further
analysis of the situation gives coherence data for this family of
maps, so that $B(f)$ has the structure of a monoid in $\calC$.

Repeating the above discussion using the subcategory $\calD_0
\subseteq \calD$ consisting of those finite linearly ordered sets
with {\em no} marked points, we can construct another monoid
object of $\calC$. Since $\calD_0$ has a final object, the colimit
is trivial and we recover the monoid $A$ from this construction.
Moreover, the inclusion $\calD_0 \subseteq \calD$ gives a map $A
\rightarrow B(f)$ of monoids.

On the other hand, consider the object $J = \{ \ast \}$ consisting
of a single marked point. By construction we have a map $M = G(J)
\rightarrow B(f)$, through which the monoid morphism $A
\rightarrow B(f)$ factors. To complete the proof, it suffices to
show that $B(f)$ is universal with respect to this property. This
(ultimately) comes down the universal property of $\calD$.
\end{proof}

\end{document}